\documentclass[11pt,a4paper]{article}
\usepackage{amssymb}
\usepackage{amsmath}
\usepackage{amsthm}
\usepackage{pstricks}
\usepackage{graphbox}
\usepackage{color}
\usepackage{mathrsfs}
\usepackage{hyperref}

\newtheorem{thm}{Theorem}[section]
\newtheorem{lem}[thm]{Lemma}
\newtheorem{rem}[thm]{Remark}

\usepackage{kpfonts}
\newcommand*{\vv}[1]{\vec{\mkern0mu#1}}

\newcommand{\bR}{{\mathbb R}}
\newcommand{\bN}{{\mathbb N}}

\newcommand{\MD}{\mathcal{D}}
\newcommand{\MI}{\mathcal{I}}
\newcommand{\MR}{\mathcal{R}}
\newcommand{\vol}{\operatorname{vol}}
\def\widebar{\overline}

\newcommand{\dH}{{\rm d}\mathscr{H}}
\newcommand{\rd}{\;{\rm d}}
\newcommand{\Id}{{\rm Id}}
\newcommand{\id}{{\rm id}}
\newcommand{\dd}[1]{\frac{\rm d}{{\rm d}#1}}
\newcommand{\ddt}{\dd{t}}
\newcommand{\nn}{\nonumber}
\newcommand{\ttau}{\Delta t}

\newcommand{\VOh}{V^h(\Omega^h)}
\newcommand{\VGh}{V^h(\Gamma^m)}
\newcommand{\Vhpartial}{V^{h}_\partial(\Gamma^m)}
\newcommand{\Wh}{W^h(\Gamma^m)}
\newcommand{\VG}{V(\Gamma)}
\newcommand{\VO}{V(\Omega)}
\newcommand{\W}{{W}(\Gamma)}
\newcommand{\Wt}{{W}(\Gamma(t))}

\newcommand{\Vpartial}{V_\partial(\Gamma)}
\newcommand{\Vpartialt}{V_\partial(\Gamma(t))}

\newcommand{\tG}{{\widetilde{G}}}
\newcommand{\bl}{{\ell}}
{{\upshape\bfseries AMS subject classifications. }\ignorespaces}{}
\newenvironment{keywords}{{\upshape\bfseries Key words. }\ignorespaces}{}

\begin{document}
%\begin{frontmatter}
\title{%\texttt{\jobname} \textbf{DRAFT} \quad {\rm \mydate} \\
A structure-preserving finite element approximation of surface diffusion for curve networks and surface clusters
}

\author{Weizhu Bao\footnotemark[1] \and 
        Harald Garcke\footnotemark[2] \and 
        Robert N\"urnberg\footnotemark[3] \and Quan Zhao\footnotemark[2] }

\renewcommand{\thefootnote}{\fnsymbol{footnote}}
\footnotetext[1]{Department of Mathematics, National University of Singapore, 119076, Singapore}
\footnotetext[2]{Fakult{\"a}t f{\"u}r Mathematik, Universit{\"a}t Regensburg, 
93040 Regensburg, Germany}
\footnotetext[3]{Dipartimento di Mathematica, Universit\`a di Trento,
38123 Trento, Italy}

\date{}

\maketitle

\begin{abstract}
We consider the evolution of curve networks in two dimensions (2d) and surface clusters in three dimensions (3d). The motion of the interfaces is described by surface diffusion, with boundary conditions at the triple junction points\slash lines, where three interfaces meet, and at the boundary points\slash lines, where an interface meets a fixed planar boundary. We propose a parametric finite element method based on a suitable variational formulation. The constructed method is semi-implicit and can be shown to satisfy the volume conservation of each enclosed bubble and the unconditional energy-stability, thus preserving the two fundamental geometric structures of the flow. Besides, the method has very good properties with respect to the distribution of mesh points, thus no mesh smoothing or regularization technique is required. A generalization of the introduced scheme to the case of anisotropic surface energies and non-neutral external boundaries is also considered. Numerical results are presented for the evolution of two-dimensional curve networks and three-dimensional surface clusters in the cases of both isotropic and anisotropic surface energies. 
\end{abstract} 

\begin{keywords} Surface diffusion, curve networks, surface clusters, 
triple junctions, volume conservation, unconditional stability, anisotropy
\end{keywords}

%\end{frontmatter}

\renewcommand{\thefootnote}{\arabic{footnote}}

\setcounter{equation}{0}
\section{Introduction} \label{sec:intro}

\setlength\parindent{24pt}

%% curve networks and surface clusters and applications
 A droplet or soap bubble tends to form a spherical geometry in order to minimize the surface area with a prescribed volume. The soap bubble cluster is a generalization to minimizing the surface area for a number of enclosed regions with prescribed volumes. Such minimizing problems have received a lot of attention in the literature, with many questions remaining open. For example, natural conjectures are that the standard $k$-bubble is the unique global minimizer among all bubbles separating $k$ different volumes, where the surfaces making up these minimizers are spherical, i.e, they are either flat or part of a sphere. A definition of  standard $k$-bubbles and a proof of the existence and uniqueness of standard bubble clusters of given volumes can be found in \cite{amilibia2001}. However, in general it is not known that they minimize surface area when the volumes are given and whether other minimizers exist. In 2d, this was proved for double bubbles ($k =2$) \cite{Foisy1993} and triple bubbles ($k=3$) \cite{Wichiramala04}, and recently Paolini and Tortorelli proved it for the quadruple planar bubble ($k=4$) enclosing equal areas \cite{PaoliniT20}.  In 3d, the double bubble conjecture was proved in \cite{HMRR}, but it is still unknown for triple and quadruple bubbles. In addition, numerical approximations have shown that for bubbles with $k\geq 6$ enclosed regions, parts of the boundaries of locally stable clusters could be non-spherical \cite{sullivan1996}. The readers are referred to \cite{Taylor76,amilibia2001,Morgan2007,Morgan1998wulff, Wecht2000double} and the references therein for more details on this topic.  

The surface diffusion flow has applications in materials science, and
geometrically can be studied as a way to obtain perimeter and surface area minimizers for given prescribed volumes, often called soap bubble clusters.
In this work, we will study the numerical approximation of the surface diffusion of curve networks in 2d and surface clusters in 3d with the help of parametric finite elements, paying particular attention to the volume-preserving aspect. The networks and clusters we consider will feature both so-called triple junction points\slash lines, where three interfaces meet, as well as boundary points\slash lines, where a boundary component of an interface is constrained to lie in a fixed external plane. Moreover, in 3d four triple junction lines can meet at a quadruple junction point. For ease of presentation, from now on we will often use the 3d naming convections for interfaces, triple junctions and boundaries, referring to these as surfaces, triple junction lines and boundary lines also in the 2d situation.

For a single, closed evolving hypersurface $(\Gamma(t))_{t\geq0}$ in $\bR^d$, the motion by surface 
diffusion is given by
\begin{equation}\label{eq:SDclosed}
\mathcal{V} = - \Delta_s \varkappa,
\end{equation}
where $\mathcal{V}$ is the velocity of $\Gamma(t)$ in the direction of the unit normal
$\vec\nu$, $\Delta_s = \nabla_s \cdot \nabla_s$ is the Laplace-Beltrami
operator and $\varkappa = - \nabla_s \cdot \vec\nu$ denotes the mean curvature 
of $\Gamma$. The geometric evolution law in \eqref{eq:SDclosed} was first introduced  by Mullins 
\cite{Mullins57} to describe mass diffusion within interfaces in polycrystalline materials. Later Davi and Gurtin \cite{DaviG90} presented a derivation of the law
using principles from
% in the context of.
 rational thermodynamics. In fact, motion by surface diffusion has wide applications in materials science and solid-state physics, such as thermal grooving, void evolution in microelectronic circuits, epitaxial crystal growth, and solid-state dewetting; see e.g.\ \cite{Mullins57,LiZG99,BowerC98,AverbuchIR03,JiangZB20}. Theoretical results on existence, uniqueness and stability for surface diffusion of a single surface can be found in e.g. \cite{ElliottG97a,EscherMS98, Giga98}.

Geometrically the law \eqref{eq:SDclosed} can be viewed as a volume preserving
gradient flow for the surface area functional. In materials science and other applications,
anisotropic surface energies often play an important role. These energies take into account that 
the surface energy density may depend on the local orientation of the 
interface. The relevant evolution law is then anisotropic surface diffusion, 
defined by \eqref{eq:SDclosed} with $\varkappa$ replaced by 
the weighted mean curvature $\varkappa_\gamma = - \nabla_s \cdot \gamma'(\vec\nu)$, where $\gamma'(\vec\nu)$
denotes the so-called Cahn--Hoffmann vector \cite{CahnH74}. Here $\gamma: \bR^d\setminus\{\vec 0\} \to \bR_{>0}$ is a one-homogeneous
extension of the map $\vec\nu \mapsto \gamma(\vec\nu)$, and $\gamma'$ denotes
its gradient in $\bR^d$. For more details on anisotropic surface energies  
we refer to \cite{Giga06,DeckelnickDE05} and the references therein.  

In practical applications, clusters of surfaces with triple junction lines
may appear, see e.g.\ \cite{Mullins58,Cahn91,Pan03,grain}.
A model for surface diffusion of a network of curves has been
introduced in \cite{GarckeNC00} for $d=2$ and generalized to arbitrary space dimensions in \cite{Barrett10cluster,DepnerG13}. 
Well-posedness was shown in \cite{AbelsAG15} for $d=2$ and in \cite{GarckeG20} for higher space dimensions. We will present the precise mathematical formulation of this evolution law in Section~\ref{sec:mf} below. 
%\footnote{Harald: In \cite{GarckeNC00} only the case $d=2$ is considered. Is
%there a reference for $d=3$?} 
In \cite{Abels2019}, it was proved that the standard planar double bubbles in $\bR^2$ are stable under surface diffusion, and the result was then generalized to the high-dimensional double bubbles in \cite{DepnerG13,Garcke2021}.

%%% literature review on numerical approximations
We now give a short overview on existing work for the
numerical approximation of surface diffusion. In the absence of 
triple junctions, we focus on methods that employ parametric finite 
elements. Here the isotropic case has been considered in 
\cite{Bansch05,Barrett07,Barrett08JCP,Zhao20,Kovacs2020,Zhao2021},
while the more general anisotropic situation has been considered in
\cite{HausserV06a,Barrett07Ani,Barrett08Ani,Bao17,Zhao2020parametric,%
BaoZ20preprint,Li2020energy}. 
We note that in \cite{Barrett07} the second and third authors of this
paper, together with John W.~Barrett, introduced a novel variational
formulation of surface diffusion that upon discretization leads to a benevolent
tangential motion that guarantees nice mesh properties in practice.
We refer to the recent review article \cite{Barrett20} for more details on this
idea, including its application to the approximation of Willmore flow, 
(snow) crystal growth, two-phase flow and fluidic biomembranes.
However, the original motivation for the variational formulation pursued 
in \cite{Barrett07} 
was the numerical approximation of geometric evolution equations for curve
networks. In fact, for a well-posed formulation it is crucial to allow movement
of the triple junction points, which in turn requires a freedom in tangential
direction for the parameterizations used to describe the individual curves.
This novel approximation of curve networks was first used in
\cite{Barrett07} for surface diffusion, and then extended to more general
geometric evolution equations in \cite{Barrett07b}. The anisotropic case for
curve networks was studied in \cite{Barrett07Ani,Barrett2011},
while the method was extended to the 
evolution of surface clusters in \cite{Barrett10cluster,Barrett10}. 

For the numerical approximation of geometric evolution laws for curve networks
and surface clusters, and more generally for numerical methods to obtain
perimeter and surface area minimizing partitionings given prescribed volumes, several different
approaches are possible. The parametric finite element methods discussed so far
fall into the category of sharp interface front tracking methods. Other
examples of front tracking methods for curve networks and surface clusters 
with triple junctions include the well-known Surface Evolver by Brakke 
\cite{Brakke92,CoxGVM-PP03,CoxG04,KraynikRS04,CoxMG13}, as well as the works
\cite{bronsard1995,Thaddey99,Neubauer2002,PanW08}. An alternative sharp
interface approach is the level set method, which has been used in e.g.\
\cite{MerrimanBO94,Ruuth1998,ZhaoMOW98,SmithSC02}. On the other hand,
the phase field method, which is a diffuse interface approach, has been
employed in e.g.\ \cite{grain,GNSW,GNSSW,Robert09}.

% Our main objectives
Very recently, the first and fourth authors of this paper presented two
novel ideas for the parametric finite element approximation for the surface
diffusion of a single surface. Firstly, in \cite{Zhao2021}, building upon ideas developed in \cite{Jiang2021}, they proposed
a method with time-integrated discrete normals that enable an exact volume 
conservation for the fully discrete solutions. Secondly, in
\cite{BaoZ20preprint} they introduced an unconditionally stable method for the
situation where a surface with boundary is attached to a non-neutral external
substrate. It is the aim of this paper to combine the ideas on the numerical
approximation of surface clusters from
\cite{Barrett07,Barrett07Ani,Barrett10cluster,Barrett10}, from now on simply
referred to as ``BGN'' or ``the BGN scheme'', with the two novel ideas from
\cite{Zhao2021,BaoZ20preprint}, in order to obtain a 
structure-preserving parametric finite element method (SP-PFEM) for the 
evolution under surface diffusion of surface clusters. In particular, by using 
suitably weighted approximations of the surface normals, and similarly 
suitably weighted effective velocity vectors along the boundary lines, where
surfaces are constrained to remain attached to fixed external planes,
we are able to devise a fully discrete numerical method that
\renewcommand{\labelenumi}{(\alph{enumi})}
\begin{enumerate}
\item conserves the volume for each enclosed bubble in the cluster exactly,
\item is unconditionally stable, including in the case of attachments to 
non-neutral planar external boundaries.
\end{enumerate}
Both of the above aspects are new in the literature. In addition, on utilizing
the techniques from \cite{Barrett08Ani}, we extend our approximation to 
the anisotropic case, when the surface energy densities depend on the local
orientation of the surfaces.

%%organization of the paper
The rest of the paper is organized as follows. 
In Section~\ref{sec:mf} we describe the mathematical problem in detail and
discuss the energy decaying and volume preserving aspect of the surface
diffusion flow for surface clusters.
In Section~\ref{sec:fem} we review the weak formulation for the considered geometric equation and then introduce a parametric finite element method.  The properties of unconditional stability and volume conservation are shown for the discretized scheme. In Section~\ref{sec:ani} we  generalize the introduced scheme to the case of anisotropic surface energies. We then discuss the extension of the introduced scheme to the non-neutral external boundaries in Section~\ref{sec:exteneb}. In Section~\ref{sec:num} extensive numerical results are presented to show the applicability of the scheme. Finally, the paper is concluded in Section~\ref{sec:con}.

\setcounter{equation}{0} 
\section{Mathematical formulation} \label{sec:mf}

We follow the notations in \cite{Barrett10cluster} and specify the geometric evolution equations as follows. The evolving surface cluster is assumed to consist of $I_S$ hypersurfaces in $\bR^d$ ($d=2,3$) with $I_T$ triple junctions lines and $I_B$ boundary lines, which are denoted by
\begin{align}
\Gamma(t)&:=\left(\Gamma_1(t),~\ldots,~\Gamma_{I_S}(t)\right), \quad I_S\in\mathbb{N},\quad I_S\ge 1,\nn\\
\mathcal{T}(t)&:=\left(\mathcal{T}_1(t),~\ldots,~\mathcal{T}_{I_T}(t)\right),\quad I_T\in\mathbb{N},\quad I_T\ge 0,\nn\\
\mathcal{B}(t)&:=\left(\mathcal{B}_1(t),~\ldots,\mathcal{B}_{I_B}(t)\right),\quad I_B\in\mathbb{N},\quad I_B\geq 0.\nn
\end{align}
We introduce parameterizations of $\Gamma(t)$ using a collection of reference 
domains $\Omega:=\big(\Omega_1,~\ldots,$ $\Omega_{I_S}\big)$, which in
order to simplify the presentation we assume to be flat domains 
$\Omega_i\subset\bR^{d-1}$, $i=1,\ldots,I_S$. The generalization to the
case where the $\Omega_i$ themselves are allowed to be hypersurfaces in
$\bR^d$ is easily possible, and such a description is needed, for
example, for the trivial cluster consisting of a single closed surface.
However, for ease of notation we assume that the parameterizations 
$\vec x$ of the cluster are such that
%We introduce parameterizations of $\Gamma(t)$ using a collection of reference domains $\Omega:=\left(\Omega_1,~\ldots,~\Omega_{I_S}\right)$, with $\Omega_i\subset\bR^{d-1}$, such that
%
\begin{equation} \label{eq:para}
\vec x = \left(\vec x_1,~\ldots,~\vec x_{I_S}\right),\quad{\rm and}\quad \vec x_i: \Omega_i\times[0,T]\to \bR^d
\text{ with } \Gamma_i(t)=\vec x_i(\Omega_i,t),\quad i = 1,\ldots, I_S.
\end{equation}
For simplicity, throughout this paper we denote $\Gamma(t)=\vec x(\Omega,t)$. The velocity $\mathcal{\vv V} = (\mathcal{\vv V}_1,\dots,\mathcal{\vv V}_{I_S})$ induced by the parameterization $\vec x$ in \eqref{eq:para} is defined by
\begin{align}
\label{eq:ve}
\mathcal{\vv V}_i(\vec x_i(\vec q, t), t)=\partial_t\vec x_i(\vec q, t)\quad\forall\vec q\in \Omega_i,\quad  i = 1,\ldots, I_S.
\end{align}The motion of the surface $\Gamma_i(t)$ is given by surface diffusion 
\begin{subequations}
\label{eqn:iso}
\begin{equation}
\label{eq:isov}
\mathcal{V}_i  = -\Delta _s\mathcal{\varkappa} _i, \qquad i = 1,\ldots, I_S,
\end{equation}
where $\mathcal{V}_i=\mathcal{\vv V}_i\cdot\vec\nu_i$ denotes the velocity of $\Gamma_i(t)$ in the direction of the unit normal $\vec \nu_i$. In addition, $\varkappa_i$ is proportional to the sum of the principal curvatures, which is given by \cite{DeckelnickDE05}
\begin{align}
\label{eq:isok}
\varkappa_i\,\vec\nu_i = \sigma_i\,\Delta_s\vec\id, \qquad i = 1,\ldots, I_S,
\end{align}
\end{subequations}
where $\sigma_i$ is a positive constant representing the surface energy density of $\Gamma_i(t)$ and $\vec\id$ is the identity function in $\bR^d$.

\begin{figure}[tph]
\centering
\includegraphics[width=0.9\textwidth]{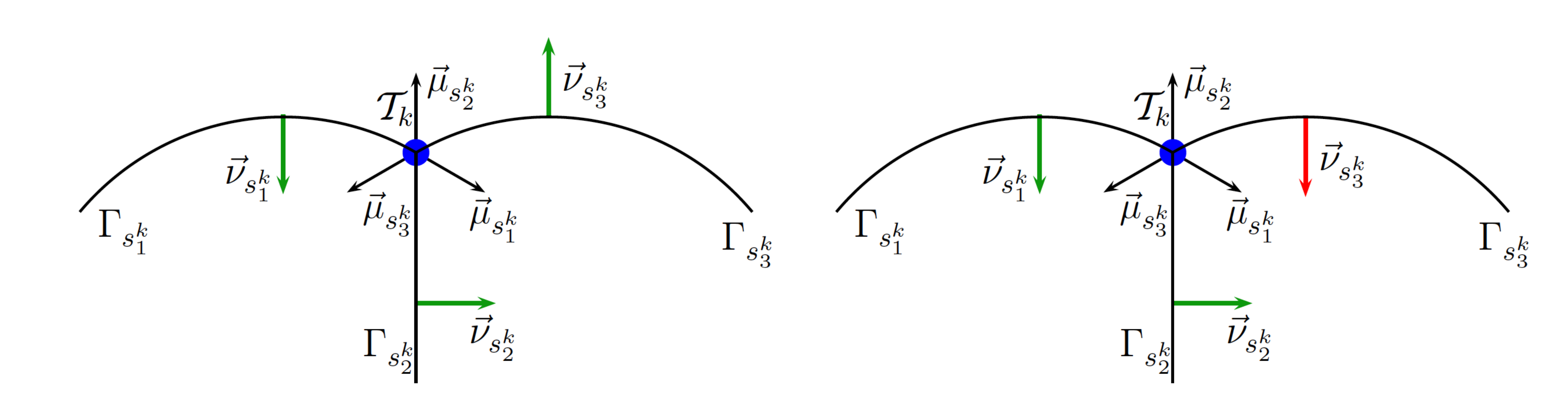}
\caption{Sketch of the local orientation of 
$(\Gamma_{s^k_1},\Gamma_{s^k_2},\Gamma_{s^k_3})$ at the triple junction line
$\mathcal{T}_k$ ({\blue\bf blue}). 
Depicted above is a plane that is perpendicular to $\mathcal{T}_k$. Left panel: 
$o^k:=(o^k_1,o^k_2,o^k_3)$ can be chosen
as $o^k=(1,1,1)$. Right panel: we require
$o^k=\pm(1,1,-1)$.}
\label{fig:orient}
\end{figure}

For the above geometric flows, we need to impose boundary conditions at the triple junction lines and boundary lines. We denote by
\begin{equation}
\partial_i\Omega = \bigcup_{j=1}^{I_P^i}\partial_j\Omega_i,\quad I_p^i\in\mathbb{N},\quad I_p^i\ge 1, \quad i = 1,\ldots, I_S\nn
\end{equation}
a partition of the boundary of $\Omega_i$. For each triple junction line $\mathcal{T}_k$, we set 
\begin{subequations}
\label{eqn:tj}
\begin{equation}
\mathcal{T}_k(t) := \vec{x}_{s^k_1}(\partial_{p^k_1}\Omega_{s^k_1},t) =
\vec{x}_{s^k_2}(\partial_{p^k_2}\Omega_{s^k_2},t) =
\vec{x}_{s^k_3}(\partial_{p^k_3}\Omega_{s^k_3},t)\,,
\quad k=1,\ldots, I_T
\,, \label{eq:tj_a}
\end{equation}
where $1\leq s_1^k<s_2^k<s_3^k\leq I_S$ and $1\leq p_j^k\leq I_P^{s_j^k}$,
$j=1,\ldots,3$. 
As a result, we can define $\mathcal{T}_k$ via the three pairs $\bigl((s_j^k,~p_j^k)\bigr)_{j=1}^3$, $k=1,\ldots, I_T$. Let $\vec\mu_i$ denote the conormal of $\Gamma_i(t)$, i.e., it is the outward unit normal to $\partial\Gamma_i(t)$ 
that lies within the tangent plane of $\Gamma_i(t)$. Then we have the following conditions  on $\mathcal{T}_k$ for $k=1,\ldots, I_T$
\begin{align}
\label{eq:tj_b}
&\sum_{j=1}^3\sigma_{s_j^k}\vec\mu_{s^k_j} = \vec 0,\\
\label{eq:tj_c}
 &o^k_1\,\vec\mu_{s^k_1} \cdot \nabla_s\,\varkappa_{s^k_1} =o^k_2\,
\vec\mu_{s^k_2} \cdot \nabla_s\,\varkappa_{s^k_2} =
o^k_3\,\vec\mu_{s^k_3} \cdot \nabla_s\,\varkappa_{s^k_3}, \\
&\sum_{j=1}^3 o^k_j\,\varkappa_{s^k_j} = 0,
\label{eq:tj_d}
\end{align}
\end{subequations}
where $o^k = \left(o_1^k,~o_2^k,~o_3^k\right)$ with $o_j^k\in\{1,-1\}$ representing the orientation of a triple junction point at $\mathcal{T}_k$ such that $\left(o_j^k\,\vec\nu_{s_j^k},~\vec\mu_{s_j^k}\right)$, $1\leq j\leq 3$, have the same orientation in the plane orthogonal to $\mathcal{T}_k$ at that point (see Fig.~\ref{fig:orient}). The equations \eqref{eq:tj_b} are force balance conditions at $\mathcal{T}_k$, which lead to the well-known $120^\circ$ angle condition at the triple junction lines when $\sigma_i$ are equal for $i=1,\ldots, I_S$.  Moreover, \eqref{eq:tj_c} and \eqref{eq:tj_d} can be interpreted as the flux balance condition and the chemical potential continuity condition, respectively. 

 We assume that part of the surfaces $\Gamma_i$, $i=1,\ldots, I_S$, are constrained to lie on the external planar surfaces $\{\mathcal{D}_k\}_{k=1}^{I_B}$.  Denote by
\begin{align}
\mathcal{B}_k(t):=\vec x_{s_k}(\partial_{p_k}\Omega_{s_k}, t)\subset\mathcal{D}_k,\quad k = 1,\ldots, I_B,\quad 1\leq s_k\leq I_S,\quad 1\leq p_k\leq I_p^{s_k},
\label{eq:attachment}
\end{align}
where $\mathcal{D}_k$ is a planar surface and its intersection with $\Gamma_{s_k}(t)$ produces the boundary line $\mathcal{B}_k(t)$. We assume for simplicity that no triple junction line $\mathcal{T}_k(t)$ is constrained to lie on the boundary, i.e., 
\begin{equation*}
 \bigcup_{k=1}^{I_B} \{ (s_k, p_k) \} \cap 
 \bigcup_{k=1}^{I_T} \bigcup_{j=1}^3 \{(s^k_j, p^k_j)\} = \emptyset.
\end{equation*}
For $1\leq k\leq I_B$, let $\vec n_k$ be the unit normal to $\mathcal{D}_k$, and pointing towards the clusters. Then we have the following conditions on $\mathcal{B}_k$ for $k=1,\ldots, I_B$
\begin{subequations}\label{eqn:bd}
\begin{align}
\label{eq:bda}
\vec n_k\cdot \mathcal{\vv V}_{s_k}=0,\\
\label{eq:bdb}
\vec n_k\cdot\vec\nu_{s_k}=0,\\
\vec \mu_{s_k}\cdot\nabla_s\varkappa_{s_k}=0.
\label{eq:bdc}
\end{align}
\end{subequations}
We note \eqref{eq:bda} together with the initial condition $\vec x_{s_k}(\partial_{p_k}\Omega_{s_k},0)\subset\mathcal{D}_k$ implies \eqref{eq:attachment} directly. 
Condition \eqref{eq:bdb} can be interpreted as a contact angle condition, which leads to a $90^\circ$ contact angle between $\Gamma_{s_k}$ and $\mathcal{D}_k$,
while \eqref{eq:bdc} is a zero-flux condition in order that the volume conservation is satisfied.

The relevant energy of the cluster is given by the weighted sum of the surface areas
\begin{align}
\label{eq:energy}
A(\Gamma(t)):=\sum_{i=1}^{I_S}\sigma_i\,|\Gamma_i(t)|=\sum_{i=1}^{I_S}\sigma_i\int_{\Gamma_i(t)}1\,\dH^{d-1}, 
\end{align}
where $\mathscr{H}^{d-1}$ denotes the $(d-1)$-dimensional Hausdorff measure in 
$\bR^d$, and similarly for $\mathscr{H}^{d-2}$.
%, and $|\Gamma_i(t)|$ represents the surface area of $\Gamma_i(t)$. 
In the cluster there are several bubbles, or volume regions, enclosed either by the surfaces or by the surfaces together with the external planar boundaries $\mathcal{D}_k$. For ease of presentation, we enumerate these regions by 
$\MR_1[\Gamma(t)],\ldots, \MR_{I_R}[\Gamma(t)]$ with corresponding index sets and orientations
\begin{align} \label{eq:IoI}
\MI_\Gamma^\ell \subset \{ 1,\ldots,I_S \}, \quad
o^{\MR_\ell} \in \{-1,1\}^{I_S}, \quad
\MI_\MD^\ell \subset \{1,\ldots,I_B\}, \quad \ell = 1,\ldots,I_R, \quad I_R \in \bN,\quad I_R\geq 1,
\end{align}
and denote by $\MR_\ell[\Gamma(t)]$ the region enclosed by the surfaces $\{\Gamma_i(t)\}_{i\in\MI^\ell_\Gamma}$, $\{\MD_k\}_{k\in \MI^\ell_\MD}$ and possibly an additional fixed hypersurface to create a finite volume. 
Here the orientations are chosen such that $o^{\MR_\ell}_i\vec\nu_i$ is the outer normal to $\MR_\ell[\Gamma(t)]$ on $\Gamma_i(t)$.
The geometric evolution equations in \eqref{eqn:iso} with the boundary conditions in \eqref{eqn:tj} and \eqref{eqn:bd} can be interpreted as a volume-preserving gradient flow. In other words, the dynamic system satisfies two geometric properties: (i) dissipation of the energy and (ii) conservation of the volume of each enclosed bubble.
In fact, it follows from a transport theorem, \eqref{eq:isok}, \eqref{eq:isov}, 
\eqref{eq:tj_a}, \eqref{eq:tj_b}, \eqref{eq:bda}, \eqref{eq:bdb}, 
\eqref{eq:tj_c}, \eqref{eq:tj_d} and \eqref{eq:bdc} that
\begin{subequations} \label{eqn:evlaws}
\begin{align}
\label{eq:energyd}
\ddt A(\Gamma(t)) & 
=-\sum_{i=1}^{I_S} \sigma_i 
\int_{\Gamma_i(t)} \frac1{\sigma_i}\varkappa_i\mathcal{V}_i \dH^{d-1}
+ \sum_{i=1}^{I_S} \sigma_i \int_{\partial\Gamma_i(t)} \vec{\mathcal{V}}_i \cdot
\vec\mu_i \dH^{d-2} %\nonumber \\ &
=\sum_{i=1}^{I_S}\int_{\Gamma_i(t)} \varkappa_i \Delta_s
\varkappa_i \dH^{d-1} \nonumber \\ &
=-\sum_{i=1}^{I_S}\int_{\Gamma_i(t)}|\nabla_s\varkappa_i|^2\,\dH^{d-1}
+ \sum_{i=1}^{I_S}\int_{\partial\Gamma_i(t)} \varkappa_i \nabla_s \varkappa_i 
\cdot\vec\mu_i \dH^{d-2} \nonumber \\ &
=-\sum_{i=1}^{I_S}\int_{\Gamma_i(t)}|\nabla_s\varkappa_i|^2\,\dH^{d-1}
\leq 0.
\end{align}
Moreover, it follows from the Reynolds transport theorem for any 
$\ell = 1,\ldots, I_R$ that
\begin{align}
\ddt \vol(\MR_\ell[\Gamma(t)]) &
= \sum_{i\in\MI^\ell_\Gamma}\int_{\Gamma_i(t)}o^{\MR_\ell}_i\mathcal{V}_i
\dH^{d-1}
= -\sum_{i\in\MI^\ell_\Gamma}o^{\MR_\ell}_i \int_{\Gamma_i(t)}\Delta_s
\varkappa_i \dH^{d-1}\nonumber\\ &
= -\sum_{i\in\MI^\ell_\Gamma}o^{\MR_\ell}_i \int_{\partial\Gamma_i(t)}\nabla_s
\varkappa_i \cdot \vec\mu_i  \dH^{d-2}
=0,
\label{eq:volumec}
\end{align}
\end{subequations}
where in the last line we have noted \eqref{eq:bdc} for the boundary lines,
and that all other boundary contributions correspond to surfaces meeting
pairwise at triple junction lines, with the chosen orientations meaning that
\eqref{eq:tj_c} implies pairwise cancellation, see also Fig.~\ref{fig:orient}
and the end of the proof of Theorem~\ref{thm:thE} below.

It is the main aim of this work to devise a fully discrete numerical method
that mimics the two fundamental structures of the flow in \eqref{eqn:evlaws} 
on the discrete level.

\setcounter{equation}{0} 
\section{Finite element approximation} \label{sec:fem}

In this section, we first revisit the BGN weak formulation for the considered geometric flow and then present a structure-preserving parametric finite element method for it.

\subsection{The weak formulation}
Let
\begin{align*}
\VO & := \big\{ (\vec\chi_1,\ldots,\vec\chi_{_{I_S}}) \in \mathop{\times}_{i=1}^{I_S} [H^1(\Omega_i)]^d :
\vec{\chi}_{s^k_1}(\partial_{p^k_1}\Omega_{s^k_1}) =
\vec{\chi}_{s^k_2}(\partial_{p^k_2}\Omega_{s^k_2}) =
\vec{\chi}_{s^k_3}(\partial_{p^k_3}\Omega_{s^k_3})\,,\, k = 1,\ldots, I_T\big\}.\nn
\end{align*}
Now any $\vec x \in \VO$ parameterizes a surface cluster $\Gamma = \vec
x(\Omega)$. Given such a cluster, we introduce the function spaces
\begin{align}
\W &:= \big\{(\chi_1,\ldots,\chi_{_{I_S}}) \in \mathop{\times}_{i=1}^{I_S} H^1(\Gamma_i) : 
\sum_{j=1}^3 o^k_j\,{\chi}_{s^k_{j}}  = 0\ 
\mbox{ on $\mathcal{T}_k$},\; k = 1,\ldots, I_T\big\} \,,\nn\\
\VG &:= \big\{(\vec\chi_1,\ldots,\vec\chi_{_{I_S}}) \in \mathop{\times}_{i=1}^{I_S} [H^1(\Gamma_i)]^d: \vec{\chi}_{s^k_{1}} = \vec{\chi}_{s^k_{2}} =
\vec{\chi}_{s^k_{3}} \ \mbox{ on $\mathcal{T}_k$},\; k = 1,\ldots, I_T\big\}\,, \nn \\
\Vpartial &:=\bigl\{(\vec\chi_1,~\ldots,\vec\chi_{_{I_S}})\in\VG: \vec \chi_{s_k}\cdot\vec n_k = 0 \ \mbox{ on $\mathcal{B}_k$},\; k = 1,\ldots, I_B\big\},\nn
\end{align}
and the $L^2$ inner product over $\Gamma$ as
\begin{align}
\left<u,~v\right>_{\Gamma}:=\sum_{i=1}^{I_S}\int_{\Gamma_i}u_i\cdot v_i\,\dH^{d-1},
\end{align}
where we allow $u, v$ to be scalar, vector or tensor valued functions.

We then introduce the weak formulation for the considered flow, i.e.\ \eqref{eqn:iso} with boundary conditions \eqref{eqn:tj} and \eqref{eqn:bd}, as follows. 
Let $\vec x(\cdot,0) \in \VO$, and $\vec x_{s_k}(\partial_{p_k}\Omega_{s_k},0)\subset\mathcal{D}_k$, $k=1,\ldots, I_B$. 
For $t>0$, we find $\vec x(\cdot,t)\in\VO$ such that $(\mathcal{\vv V}(\cdot, t),~\varkappa(\cdot, t))\in\Vpartialt\times \Wt$, for $\Gamma(t) = \vec x(\Omega,t)$, with
\begin{subequations}
\label{eqn:weak}
\begin{align}
\label{eq:weaka}
\big<\mathcal{\vv V}\cdot\vec\nu,~\chi\big>_{\Gamma(t)}-\big<\nabla_s\varkappa,~\nabla_s\chi\big>_{\Gamma(t)}=0\quad\forall\chi\in \Wt,\\[0.4em]
\big<\varkappa\,\vec\nu,~\vec\eta\big>_{\Gamma(t)}+\big<\sigma\,\nabla_s\vec\id, ~\nabla_s\vec\eta\big>_{\Gamma(t)}=0\quad\forall\vec\eta \in\Vpartialt.
\label{eq:weakb}
\end{align}
\end{subequations}
 Here \eqref{eq:weaka} is obtained by multiplying \eqref{eq:isov} with $\chi_i$, integrating over $\Gamma_i$, summing up for $i=1,\ldots, I_S$, using integration by parts and the boundary conditions \eqref{eq:tj_c}, \eqref{eq:bdc}. Similarly, using test functions $\vec\eta\in\Vpartial$ to multiply \eqref{eq:isok}, we can obtain \eqref{eq:weakb} by noting the boundary conditions \eqref{eq:tj_b} and \eqref{eq:bdb}.

\subsection{The discretization}
\label{sse:dis}
 For $i=1 ,\ldots, I_S$, let $\Omega_i^h = \cup_{j=1}^{J_i} \overline{\sigma^i_j}$ be a triangulation approximating 
$\overline\Omega_i \subset \bR^{d-1}$, where $\{\sigma^i_j\}_{j=1}^{J_i}$ is a family of mutually disjoint open 
$(d-1)$-simplices with vertices $\{\vec{q}^i_k\}_{k=1}^{K_i}$. Denote by $\partial_j\Omega_i^h$ an approximation of $\partial_j\Omega_i$, $j=1,\ldots, I^i_P$, $i=1,\ldots, I_S$. Then we
assume that the endpoints of $\partial_j\Omega_i^h$ and
$\partial_j\Omega_i$ coincide and that  the
triangulations of $\Omega^h$ ``match up'' at their boundaries at triple junction lines, i.e.,
\begin{align}
Z_k := \# \{ \{\vec{q}^{s^k_1}_l\}_{l=1}^{K_{s^k_1}} 
\cap \partial_{p^k_1}\Omega_{s^k_1}^h \} =
\# \{ \{\vec{q}^{s^k_2}_l\}_{l=1}^{K_{s^k_2}} 
\cap \partial_{p^k_2}\Omega_{s^k_2}^h \} =
\# \{ \{\vec{q}^{s^k_3}_l\}_{l=1}^{K_{s^k_3}} 
\cap \partial_{p^k_3}\Omega_{s^k_3}^h \},\quad k = 1,\ldots, I_T.\nn
\end{align}
In addition, for the discrete boundary parts $\partial_{p^k_j}\Omega_{s^k_j}^h$, we let
\begin{equation}
\vec\rho^k_j : \{1,\ldots, Z_k \} \to
\big\{ \{\vec{q}^{s^k_j}_l\}_{l=1}^{K_{s^k_j}} 
\cap \partial_{p^k_j}\Omega_{s^k_j}^h \big\}\,,\quad j = 1,\ldots, 3 \,,\quad 1,\ldots, I_T,
\label{eq:rho}
\end{equation}
be a bijective map such that
$(\vec\rho^k_j(1),\ldots,\vec\rho^k_j(Z_k))$ is an ordered sequence of vertices. Then we define the natural discrete analogue of $\VO$ by
\begin{align}  
\VOh &=  \Big\{(\vec \chi_1,~\ldots,~\vec \chi_{I_S})
\in \mathop{\times}_{i=1}^{I_S} [C^0(\overline\Omega_i^h)]^d : \vec\chi_i
\!\mid_{\sigma_j^i}
\mbox{ is linear}\ \forall\ j=1,\ldots, J_i,\ i = 1,\ldots, I_S; \nn\\
&\hspace{0.2cm}\vec{\chi}_{s^k_{1}} (\vec{\rho}^k_1(l)) =
\vec{\chi}_{s^k_{2}} (\vec{\rho}^k_2(l)) =
\vec{\chi}_{s^k_{3}} (\vec{\rho}^k_3(l))\,,\ 
l = 1,\ldots, Z_k,\, k = 1,\ldots, I_{T}\Big \}.
\end{align}

Let $M$ be a positive integer and $\bigcup_{m=0}^{M-1}[t_m,~t_{m+1}]$ be a partition of the time domain $[0,~T]$ such that $0=t_0<t_1<\ldots<t_M=T$ with possibly variable time steps $\ttau_m := t_{m+1} -
t_{m}$.  Denote by $\Gamma^m=\vec {\mathcal X}^m(\Omega^h)$, for $\vec {\mathcal X}^m \in\VOh$, the discrete approximation of the cluster $\Gamma(t_m)$, with $\Gamma_i^m=\vec {\mathcal X}_i^m(\Omega_i^h)$, $i=1,\ldots, I_S$. This introduces a sequence of polyhedral surfaces in $\bR^d$. Let $\Gamma_i^m=\bigcup_{i=1}^{J_i}\overline{\sigma_j^{m,i}}=\bigcup_{j=1}^{J_i}\vec {\mathcal X}^m_i(\overline{\sigma_j^i})$,  where $\{\sigma_j^{m,i}\}_{j=1}^{J_i}$ are mutually disjoint open $(d-1)$-simplices with vertices $\{\vec q_k^{m,i}\}$ defined by $\vec{q}^{m,i}_k := \vec{{\mathcal X}}^m_i(\vec{q}^i_k)$. As a discrete analogue of $\mathcal{T}_k(t_m)$, the triple junction $\mathcal{T}_k^m$ of the polyhedral surface cluster $\Gamma^m$ is defined by the ordered sequence of vertices 
\begin{align}
(\vec{{\mathcal X}}^m_{s_1^k} (\vec\rho_1^k(1)), \ldots,\linebreak
\vec{{\mathcal X}}^m_{s_1^k} (\vec\rho_1^k(Z_k))),\quad k = 1,\ldots, I_T. \nn
\end{align}
Similarly, the boundaries $\mathcal{B}_k^m$ are given by an appropriately defined
ordering of the vertices $\{\vec {\mathcal X}^m(\vec q):\vec q\in \{\vec q_l^{s_k}\}_{k=1}^{K_{s_k}}\cap \partial_{p_k}\Omega^h_{s_k}\}$.

 We  define the function spaces $\widehat{W}^h(\Gamma^m) := \{\chi \in \mathop{\times}_{i=1}^{I_S}
C^0(\Gamma_i^m): \chi_i
\!\mid_{\sigma_j^{m,i}}
\mbox{ is linear}\ \forall\ j=1,\ldots, J_i,\ i = 1,\ldots, I_S\}$ and 
$\widehat{V}^h(\Gamma^m):=\{\vec\chi \in \mathop{\times}_{i=1}^{I_S} [C^0(\Gamma_i^m)]^d: \vec\chi_i
\!\mid_{\sigma_j^{m,i}}
\mbox{ is linear}\ \forall\ j=1,\ldots, J_i,\ i = 1,\ldots, I_S\}.$
Then the natural discrete analogues of $\VG$, $\W$ and $\Vpartial$ are given by
\begin{subequations}
\begin{align}
\Wh&:=\Big\{\chi\in\widehat{W}^h(\Gamma^m) : \sum_{j=1}^3 o_j^k\chi_{s_j^k}=0 \mbox{ on $\mathcal{T}_k^m$}, 
\; k = 1,\ldots, I_{T}\Big\},\\
\VGh&:=\Big\{\vec\chi\in\widehat{V}^h(\Gamma^m) : \vec\chi_{s_1^k}=\vec\chi_{s_2^k}=\vec\chi_{s_3^k}\mbox{ on $\mathcal{T}_k^m$}, \; k = 1,\ldots, I_{T}\Big\},\\
\Vhpartial&:=\Big\{\vec\chi\in\VGh:\;\vec n_k\cdot\vec\chi_{s_k}(\vec q)=0\;\forall\vec q\in\mathcal{B}_k^m, k= 1,\ldots, I_B \Big\}.
\end{align}
\end{subequations}

In addition, let
$\left\{\vec q_{j_k}^{m,i}\right\}_{k=0}^{d-1}$ be the vertices of $\sigma_j^{m,i}$, and ordered with the same orientation for all $\sigma_j^{m,i}$, $j=1,\ldots, J_i$. For simplicity, we denote $\sigma_j^{m,i}=\Delta\left\{\vec q_{j_k}^{m,i}\right\}_{k=0}^{d-1}$. Then we introduce the unit normal $\vec{\nu}^m_i$ to $\Gamma^m_i$; that is,
\begin{equation}\label{eq:vG}
\vec{\nu}^m_{i,j} := \vec{\nu}^m_i \mid_{\sigma^{m,i}_j} :=
\frac{\vec A\{\sigma_j^{m,i}\}}{
|\vec A\{\sigma_j^{m,i}\}|}\quad\mbox{ with}\quad \vec A\{\sigma_j^{m,i}\}=( \vec{q}^{m,i}_{j_1} - \vec{q}^{m,i}_{j_0} ) \wedge \ldots \wedge
( \vec{q}^{m,i}_{j_{d-1}} - \vec{q}^{m,i}_{j_0}),
\end{equation}
where $\wedge$ is the wedge product and $\vec A\{\sigma_j^{m,i}\}$ is the orientation vector of $\sigma_j^{m,i}$.  To approximate the inner product $\langle\cdot,\cdot\rangle_{\Gamma(t_m)}$, we introduce the inner products 
$\langle\cdot,\cdot\rangle_{\Gamma^m}$ and  $\langle\cdot,\cdot\rangle_{\Gamma^m}^h$ over
the current polyhedral surface cluster $\Gamma^m$ via 
\begin{subequations}
\begin{align}
\label{eq:erule}
\langle u, v\rangle_{\Gamma^m} &:= \sum_{i=1}^{I_S} 
\int_{\Gamma^m_i} u_i\cdot v_i \dH^{d-1},\\
\langle u, v \rangle^h_{\Gamma^m} &:= \sum_{i=1}^{I_S} 
\frac{1}{d}\sum_{j=1}^{J_i} |\sigma^{m,i}_j|
\sum_{k=0}^{d-1} 
\underset{\sigma^{m,i}_j\ni \vec{p}\to \vec{q}^{m,i}_{j_k}}{\lim}\, 
(u_i\cdot v_i)(
\vec{p}),
\label{eq:tprule}
\end{align}
\end{subequations}
where $u,v$ are piecewise continuous, with possible jumps
across the edges of $\{\sigma^{m,i}_j\}_{j=1}^{J_i}$, $i=1,\ldots, I_S$, $\{\vec{q}^{m,i}_{j_k}\}_{k=0}^{d-1}$ are the vertices of $\sigma^{m,i}_j$, and $|\sigma^{m,i}_j| = \frac{1}{(d-1)!}\,|\vec A\{\sigma_j^{m,i}\}|$ 
is the measure of $\sigma^{m,i}_j$. 

In what follows, given the cluster $\Gamma^m$ we will devise a system of 
equations for $\vec{X}^{m+1} \in \VGh$, which then defines the new cluster
$\Gamma^{m+1} = \vec X^{m+1}(\Gamma^m)$. 
Based on the ideas in \cite{Jiang2021, Zhao2021}, it is our aim to propose 
a finite element approximation of the weak formulation in \eqref{eqn:weak} 
in order that the energy dissipation law \eqref{eq:energyd} and the volume 
conservation law \eqref{eq:volumec} are still satisfied on the discrete level. 
To this end, we need to introduce appropriately weighted surface normals that 
approximate $\vec\nu_i$. Precisely, we first introduce a family of polyhedral 
surfaces via a linear interpolation between $\Gamma^m$ and $\Gamma^{m+1}$ 
defined by
\begin{equation}
\label{eq:gammah}
\Gamma^h_i(t)=\frac{t_{m+1} - t}{\Delta t_m}\Gamma_i^m + \frac{t - t_m}{\Delta t_m}\,\Gamma_i^{m+1},\quad t\in[t_m,~t_{m+1}],\quad i = 1,\ldots, I_S.
\end{equation}
Denote by $\Gamma_i^h(t)=\bigcup_{j=1}^{J_i} \overline{\sigma_j^{h,i}(t)}$ the polyhedral surfaces, where $\{\sigma_j^{h,i}(t)\}_{j=1}^{J_i}$ are the mutually disjoint $(d-1)$-simplices with vertices $\{\vec q_k^{h,i}(t)\}_{k=1}^{K_i}$, and  
\begin{equation}
\label{eq:qh}
\vec q_k^{h,i}(t) = \frac{t_{m+1} - t}{\Delta t_m}\vec q_k^{m,i} + \frac{t - t_m}{\Delta t_m}\,\vec q_k^{m+1,i},\quad t\in[t_m,~t_{m+1}],\quad k = 1,\ldots, K_i.
\end{equation}
We then define the time-weighted approximation $\vec\nu^{m+\frac{1}{2}}\in \displaystyle \mathop{\times}_{i=1}^{I_S} [L^\infty(\Gamma_i^m)]^d$  such that
\begin{align}
\label{eq:weightv}
\vec\nu_i^{m+\frac{1}{2}}|_{\sigma_j^{m,i}}=\vec\nu_{i,j}^{m+\frac{1}{2}}&:=\frac{1}{\Delta t_m\,|\vec A\{\sigma_j^{m,i}\}|}\int_{t_m}^{t_{m+1}}\vec A\{\sigma_j^{h,i}(t)\}\,\rd t,\quad j = 1,\ldots, J_i,\quad i = 1,\ldots, I_S.
\end{align}
In a similar manner as in \cite{Zhao2021} , we have the following lemma for the discrete quantities defined in \eqref{eq:weightv}.

\begin{lem}\label{lem:vc}
Let $\vec X^{m+1}\in\VGh$ with $\vec X^{m+1} - \vec\id\!\mid_{\Gamma^m}
\in\Vhpartial$. 
Then it holds
\begin{align}
\label{eq:Vlem}
\vol(\MR_\ell[\Gamma^{m+1}]) - \vol(\MR_\ell[\Gamma^m]) = 
\big<(\vec X^{m+1} - \vec\id)\cdot\vec\nu^{m+\frac{1}{2}},~\chi\big>_{\Gamma^m}^h,\quad \ell = 1,\ldots, I_R,
\end{align}
where $\chi=(\chi_1,~\ldots,~\chi_{_{I_S}})$ is given by
\begin{equation} 
\label{eq:xior}
\chi_i = \begin{cases}
o^{\MR_\ell}_i &{\text{if}}\; i\in\MI_\Gamma^\ell,\\[0.4em]
0 &{\text{if}}\; i\notin\MI_\Gamma^\ell, 
\end{cases}
\end{equation}
with $o^{\MR_\ell}$ defined as in \eqref{eq:IoI}.
\end{lem}
\begin{proof}
For $t\in[t_m,~t_{m+1}]$ and $\Gamma^h(t)=(\Gamma_1^h(t),\ldots,\Gamma_{I_S}^h(t))$ defined in \eqref{eq:gammah}, denote $\Gamma^h(t):=\vec X^h(\Gamma^m,~t)$ with $\vec X^h(t)\in\VGh$. Then we have  
\begin{align}
\label{eq:xh}
\vec X^h_i(\vec q, t) =\frac{t_{m+1}-t}{\Delta t_m}\vec q + \frac{t - t_m}{\Delta t_m}\,\vec X_i^{m+1}(\vec q),\quad\forall\vec q\in\Gamma^m,\quad t\in[t_m,~t_{m+1}],\quad i = 1,\ldots, I_S.
\end{align}
Denote by $\vec\nu^h(t)=(\vec\nu_1^h,~\ldots,\vec\nu_{I_S}^h)$ the unit normal 
to $\Gamma^h(t)$. We now apply the Reynolds transport theorem to the region
$\MR_\ell[\Gamma^h(t)]$, for $t\in[t_m,~t_{m+1}]$.
Here the boundaries $\{\MD_k\}_{k\in \MI^\ell_\MD}$ 
do not move, and so do not
contribute to the change in volume. Hence, similarly to \cite{Zhao2021} and
\cite{mullins}, we obtain that 
\begin{align}
\ddt\vol(\MR_\ell[\Gamma^h(t)])&=\sum_{i\in \MI^\ell_\Gamma}\int_{\Gamma_i^h(t)}\,o^{\MR_\ell}_i\,\vec\nu_i^h\cdot(\partial_t\vec X^h_i)\circ (\vec X^h_i)^{-1} 
\,\dH^{d-1}\nn\\
&=\sum_{i\in\MI^\ell_\Gamma}o^{\MR_\ell}_i\sum_{j=1}^{J_i} \int_{\sigma_j^{m,i}} \frac{\vec X^{m+1}_i - \vec\id}{\Delta t_m}\cdot\frac{\vec A\{\sigma_j^{h,i}(t)\}}{|\vec A\{\sigma_j^{h,i}(t)\}|}\,\frac{|\vec A\{\sigma_j^{h,i}(t)\}|}{|\vec A\{\sigma_j^{m,i}\}|}\,\dH^{d-1},
\label{eq:dv1}
\end{align}
where in the first equality we have dropped the integrals over subsets of $\{\mathcal{D}_k\}_{k\in {\MI^\ell_\MD}}$ as they are zero.
Integrating \eqref{eq:dv1} from $t_m$ to $t_{m+1}$ with respect to $t$, we arrive at
\begin{align}
&{\rm vol}(\MR_\ell[\Gamma^{m+1}]) - {\rm vol}(\MR_\ell[\Gamma^m]) \nn\\
&\qquad=\int_{t_m}^{t_{m+1}}\sum_{i\in\MI^\ell_\Gamma}o^{\MR_\ell}_i\sum_{j=1}^{J_i} \int_{\sigma_j^{m,i}}\,\frac{\vec X^{m+1}_i - \vec\id}{\Delta t_m}\cdot\frac{\vec A\{\sigma_j^{h,i}(t)\}}{|\vec A\{\sigma_j^{m,i}\}|}\,\dH^{d-1}\rd t\nn\\
&\qquad=\sum_{i\in\MI^\ell_\Gamma}o^{\MR_\ell}_i\sum_{j=1}^{J_i} \int_{\sigma_j^{m,i}}\,\left(\vec X^{m+1}_i - \vec\id\right)\cdot\frac{1}{\Delta t_m\,|\vec A\{\sigma_j^{m,i}\}|}\,\int_{t_m}^{t_{m+1}}\vec A\{\sigma_j^{h,i}(t)\}\,\rd t\,\dH^{d-1}\nn\\
&\qquad=\sum_{i\in \MI_\Gamma^\ell}o^{\MR_\ell}_i\int_{\Gamma_i^m}\left(\vec X^{m+1}_i - \vec\id\right)\cdot\vec\nu_i^{m+\frac{1}{2}}\,\dH^{d-1},
\label{eq:dv2}
\end{align}
where we have invoked \eqref{eq:weightv} for the last equality. This implies \eqref{eq:Vlem} on recalling \eqref{eq:tprule}.
\end{proof}

\begin{rem} \label{rem:poly}
We note that in \eqref{eq:weightv}, $\vec A\{\sigma_j^{h,i}(t)\}$ is a polynomial of degree $d-1$ for the variable $t$, recall \eqref{eq:vG} and \eqref{eq:qh}. Therefore, in the case of $d=2$, applying the trapezoidal rule to \eqref{eq:weightv} yields
\begin{equation*}
\vec\nu_{i,j}^{m+\frac{1}{2}} = \frac{\vec A\{\sigma_j^{m,i}\} + \vec A\{\sigma_j^{m+1,i}\}}{2|\vec A\{\sigma_j^{m,i}\}|},\quad j = 1,\ldots, J_i,\quad i = 1,\ldots, I_S,\quad m = 0,\ldots, M-1,
\end{equation*}
which gives \cite[(2.10)]{Zhao2021}. While in the case of $d=3$, we can apply Simpson's quadrature rule and obtain
\begin{equation*}
\vec\nu_{i,j}^{m+\frac{1}{2}} = \frac{\vec A\{\sigma_j^{m,i}\} + 4\vec A\{\sigma_j^{m+\frac{1}{2},i}\} + \vec A\{\sigma_j^{m+1,i}\}}{6\,|\vec A\{\sigma_j^{m,i}\}|}\quad\mbox{with}\quad \sigma_j^{m+\frac{1}{2},i} = \Delta\{\frac{\vec q_{j_k}^{m,i} + \vec q_{j_k}^{m+1,i}}{2}\}_{k=0}^{2}.
\end{equation*}
This gives a form similar to \cite[(3.12)]{Zhao2021}. 
\end{rem}

We now propose the following structure-preserving discretization for the weak formulation in \eqref{eqn:weak}. Let $\vec X^0\in\VOh$, and $\vec X^0_{s_k}(\partial_{p_k}\Omega_{s_k}^h)\subset\mathcal{D}_k$, $k=1,\ldots, I_B$. For $m=0,\ldots,M-1$, 
find $(\vec X^{m+1}, \kappa^{m+1}) \in \VGh \times \Wh$,
with $\vec X^{m+1} - \vec\id\!\mid_{\Gamma^m} \in \Vhpartial$, 
such that
\begin{subequations}
\label{eqn:fulld}
\begin{align}
\label{eq:fulld1}
\frac{1}{\Delta t_m}\big<\vec X^{m+1}-\vec\id,~\vec\nu^{m+\frac{1}{2}}\,\chi\big>_{\Gamma^m}^h -\big<\nabla_s\kappa^{m+1},~\nabla_s\chi\big>_{\Gamma^m} = 0\quad\forall\chi\in\Wh,\\
\big<\kappa^{m+1}\,\nu^{m+\frac{1}{2}},~\vec\eta\big>_{\Gamma^m}^h + \big<\sigma\,\nabla_s\vec X^{m+1},~\nabla_s\vec\eta\big>_{\Gamma^m} = 0\quad\forall\vec\eta\in\Vhpartial.
\label{eq:fulld2}
\end{align}
\end{subequations}
We note the method \eqref{eqn:fulld} is very similar to the BGN scheme, see e.g.\ \cite[(4.7)]{Barrett10cluster}. The difference is that here in the first terms of \eqref{eq:fulld1} and \eqref{eq:fulld2} we employ the semi-implicit approximation of the unit normal from \eqref{eq:weightv} instead of the explicit approximation with $\vec\nu^m$, which results in a nonlinear set of equations, compared to the linear scheme from BGN. These treatments will lead to a volume-preserving and unconditionally stable method. Furthermore, the method has very good properties with respect to the distribution of mesh points. In other words, for a semi-discrete approximation, it generally leads to the equidistribution of mesh points in 2d and conformal polyhedral surfaces in 3d, which has been studied in detail in \cite{Barrett07,Barrett08JCP}, see also \cite{Barrett20}. The discretized method gives rise to a system of nonlinear polynomial equations, 
recall Remark~\ref{rem:poly},
and in practice can be solved e.g.\ with a Picard-type iterative method, see Remark~\ref{rem:iterativemethod} below.

\begin{rem}
Formally the method \eqref{eqn:fulld} is first order in temporal discretization and second order in spatial discretization, which was numerically confirmed in \cite{Zhao2021} for surface diffusion of a single surface. However, the mathematical analysis of the error and convergence for the type of BGN schemes is an open problem and still very challenging due to the introduced tangential movements of the vertices and the complexity of the differential equations.
\end{rem}

\subsection{Volume conservation and stability}
We have the following theorem for the discretization \eqref{eqn:fulld}, which mimics the energy dissipation  and volume conservation laws in \eqref{eqn:evlaws} on the discrete level.
\begin{thm}[stability and volume conservation] \label{thm:thE}
Let $(\vec X^{m+1},~\kappa^{m+1})$ be a solution to \eqref{eqn:fulld}. 
Then it holds that
\begin{align}
\label{eq:thE}
A(\Gamma^{m+1})+\Delta t_m\big<\nabla_s\kappa^{m+1},~\nabla_s\kappa^{m+1}\big>_{\Gamma^m}\leq A(\Gamma^{m}).
\end{align}
Moreover, it holds that
\begin{align}
\label{eq:thV}
{\rm vol}(\MR_\ell[\Gamma^{m+1}]) = {\rm vol}(\MR_\ell[\Gamma^m]),\quad \ell = 1,\ldots, I_R.
\end{align}
\end{thm}
\begin{proof}
Setting $\chi = \Delta t_m\kappa^{m+1}$ in \eqref{eq:fulld1} and $\vec\eta=\vec X^{m+1}-\vec\id\!\mid_{\Gamma^m}$ in \eqref{eq:fulld2}, and combing the two equations, yields
\begin{align}
\label{eq:the12}
\big<\sigma\,\nabla_s\vec X^{m+1},~\nabla_s(\vec X^{m+1}-\vec\id)\big>_{\Gamma^m} + \Delta t_m\big<\nabla_s\kappa^{m+1},~\nabla_s\kappa^{m+1}\big>_{\Gamma^m} = 0.
\end{align}
It follows directly from \cite[Lemma~57]{Barrett20} that
\begin{align}
\sigma_i\,\int_{\Gamma_i^m}\nabla_s\vec X^{m+1}_i:\nabla_s(\vec X^{m+1}_i - \vec\id)\,\dH^{d-1}\geq \sigma_i(|\Gamma_i^{m+1}| - |\Gamma_i^m|)\quad\forall i = 1,\ldots, I_S,
\label{eq:the1}
\end{align}
which immediately implies \eqref{eq:thE} by inserting \eqref{eq:the1} into \eqref{eq:the12}.

Moreover, in \eqref{eq:fulld1} we set $\chi=(\chi_1,~\ldots,~\chi_{I_S})$ with $\chi_i$ satisfying \eqref{eq:xior}. This gives
\begin{equation}
\big<\vec X^{m+1}-\vec\id,~\vec\nu^{m+\frac{1}{2}}\chi\big>_{\Gamma^m}^h=0,
\end{equation}
which implies \eqref{eq:thV} by noting Lemma \ref{lem:vc}. What remains to be done is to show that the chosen test function satisfies $\chi\in\Wh$. 
For an arbitrary triple junction line,
if $\mathcal{T}_k^m \cap \widebar{\MR_\ell[\Gamma^m]} = \emptyset$, then
$\chi_{s^k_1} = \chi_{s^k_2} = \chi_{s^k_3} = 0$ and there is nothing to show.
Otherwise, we assume without loss of generality that $\chi_{s^k_1} = 0$ and
$\{s_2^k,s_3^k\}\subset\MI_{\Gamma}^\ell$. 
As shown in Fig.~\ref{fig:orient}, in order that 
$o^{\MR_\ell}_{s_2^k}\vec\nu_{s_2^k}^m$ and
$o^{\MR_\ell}_{s_3^k}\vec\nu_{s_3^k}^m$ 
are the outer normal to the considered region $\MR_\ell[\Gamma^m]$, 
on the left panel we require $o^{\MR_\ell}_{s_2^k}=-1$ and 
$o^{\MR_\ell}_{s_3^k}=1$, while on the right panel $o^{\MR_\ell}_{s_2^k}=-1$ 
and $o^{\MR_\ell}_{s_3^k}=-1$. 
In both cases $\sum_{j=1}^3 o_j^k \chi_{s_j^k}=o_2^k o^{\MR_\ell}_{s_2^k}
+ o_3^k o^{\MR_\ell}_{s_3^k} = 0$ holds, and thus $\chi\in\Wh$.
\end{proof}

\begin{rem} In the case of curved boundaries $\mathcal{D}_k$, the attachment condition \eqref{eq:attachment} will only be approximately satisfied. Usually an orthogonal projection of $\vec X^{m+1}$ onto $\mathcal{D}_k$ can be employed so that the attachment condition is exactly satisfied.  But the price is that the numerical solutions will lose the properties of volume conservation and unconditional stability. Therefore, we restrict our attention to the case of planar external boundaries in this work.   
\end{rem}

\setcounter{equation}{0} 
\section{Anisotropic surface energies} \label{sec:ani}

\subsection{Mathematical formulations}
\label{sec:amathfor}

In materials science, the surface energy of a material often exhibits strong dependence on its crystallographic orientations. This yields the anisotropy and could influence the kinetic evolution of the material. To this end, we assume the anisotropic surface energy density for the cluster $\Gamma(t)=(\Gamma_1(t),~\ldots,\Gamma_{I_S}(t))$ is given by $\gamma$. In particular, we restrict ourselves to the surface energy of the form that was introduced in \cite{Barrett08Ani}:
\begin{align}
%\label{eq:anif}
\gamma(\vec p ) = \left(\sum_{\ell=1}^{L}[\gamma_\bl(\vec p)]^{r}\right)^{\frac{1}{r}}\quad\mbox{with}\quad \gamma_\bl(\vec p):=\sqrt{\vec p\cdot G_\bl\vec p},\quad r\in[1,\infty), \quad\forall\vec p\in\bR^d\backslash\{\vec 0\},\nn
\end{align}
where $G_\bl\in\bR^{d\times d}$, $\ell=1,\ldots, L$, are symmetric and positive definite. 
Building on the techniques in \cite{Barrett08Ani}, the restriction to this
class of anisotropies will allow us to establish an analogue of
Theorem~\ref{thm:thE} for the anisotropic generalization of the scheme
\eqref{eqn:fulld}. 
Direct calculation yields the gradient of $\gamma(\vec p)$ as
\begin{align}
\label{eq:anipy}
\gamma'(\vec p) =\sum_{\ell=1}^{L}\left[\frac{\gamma_\bl(\vec p)}{\gamma(\vec p)}\right]^{r-1}\gamma_\bl'(\vec p)\quad\mbox{with}\quad \gamma_\bl^\prime(\vec p) = \frac{1}{\gamma_\bl(\vec p)}G_\bl\vec p.
\end{align}
%We observe that $\gamma_i(\vec p)$, $i=1,\ldots, I_S$ are positive and absolutely homogeneous of degree one, and they satisfy 
%\begin{align}
%\gamma_i(\lambda\vec p) = |\lambda|\gamma_i(\vec p)%\quad\forall\lambda\in \bR\quad\Rightarrow\quad \gamma_i^\prime(\vec p)\cdot\vec p = \gamma_i(\vec p)\quad\forall\vec p\in\bR^3\backslash\{\vec 0\}.
%\end{align}
Some typical examples of $\gamma(\vec p)$ are the isotropic surface energy with $L=1, r=1, G_1=\Id\in\bR^{d\times d}$, which gives $\gamma(\vec p ) = |\vec p|$, as well as $L=d$ with
\begin{align}
\label{eq:cuspgamma}
\gamma(\vec p) = \left(\sum_{\ell=1}^d\left[(1-\epsilon^2) p_\ell^2 + \epsilon^2 |\vec p|^2\right]^{\frac{r}{2}}\right)^\frac{1}{r},\qquad \vec p = (p_1,\dots,p_d)^T.
\end{align}
In the case of $r=1$, \eqref{eq:cuspgamma} can be regarded as a smooth regularization of the $l^1$-norm $\gamma(\vec p)=\sum_{\ell=1}^d|p_\ell|$, while for $r\gg1$ and $\epsilon\ll1$ it approximates an octahedral anisotropy in the case $d=3$. For more choices of $L$, $r$, $G_\bl$ and their corresponding Wulff shapes,  readers can refer to Refs.~\cite{Barrett07Ani,Barrett08Ani,Barrett10cluster} and the references therein. 

We now generalize the gradient flow in \eqref{eqn:iso}, with boundary conditions \eqref{eqn:tj} and \eqref{eqn:bd}, to the case of anisotropic surface energies. The motion of $\Gamma_i(t)$ is given by the anisotropic surface diffusion
\begin{subequations}
\label{eqn:anivk}
\begin{align}
\label{eq:aniv}
\mathcal{V}_i = -\Delta_s\varkappa_{\gamma,i},\quad i = 1,\ldots, I_S,
\end{align}
where $\varkappa_{\gamma,i}$ for $i=1,\ldots, I_S$ are the weighted mean curvatures and are defined via the Cahn-Hoffman vector $\vec\nu_{\gamma,i}$ \cite{Hoffman72, CahnH74}: 
\begin{align}
\label{eq:anik}
\varkappa_{\gamma,i}=-\nabla_s\cdot\vec\nu_{\gamma,i}\quad\mbox{with}\quad \vec\nu_{\gamma,i}=\gamma'(\vec\nu_i).
\end{align}
\end{subequations}

We next consider the boundary conditions for the anisotropic system. At the triple junction lines $\mathcal{T}_k$, $k=1,\ldots, I_T$, we still have the attachment conditions \eqref{eq:tj_a}. The anisotropic variants of \eqref{eq:tj_b} - \eqref{eq:tj_d} are then given by \cite{Hoffman72,GarckeNS98,GarckeNC00,Taylor99} 
\begin{subequations}
\label{eqn:atj}
\begin{align}
&\sum_{j=1}^3 \biggl[\gamma(\vec\nu_{s^k_j})\,\vec \mu_{s^k_j} -
(\gamma'(\vec \nu_{s^k_j})\cdot\vec \mu_{s^k_j})\,\vec{\nu}_{s^k_j}\biggr] = 
\vec 0,
\label{eq:atj_b} \\
&o^k_1\,\vec\mu_{s^k_1} \,
\nabla_s\,\varkappa_{\gamma,s^k_1} =
o^k_2\,\vec\mu_{s^k_2} \,\nabla_s\,\varkappa_{\gamma,s^k_2} =
o^k_3\,\vec\mu_{s^k_3} \,
\nabla_s\,\varkappa_{\gamma,s^k_3}, \label{eq:atj_c} \\
&\sum_{j=1}^3 o^k_j\,\varkappa_{\gamma,s^k_j} = 0\,.
\label{eq:atj_d}
\end{align}
\end{subequations}
At the boundary lines $\mathcal{B}_k$, $k=1,\ldots, I_B$,  we still require \eqref{eq:bda} to hold so that the boundary lines remain attached to the external planes. The generalizations of \eqref{eq:bdb}, \eqref{eq:bdc} are given by 
\begin{subequations}
\label{eqn:abd}
\begin{align}
\label{eq:abdb}
\vec n_k\cdot\gamma^\prime (\vec\nu_{s_k})=0,\\
\vec \mu_{s_k}\cdot\nabla_s\varkappa_{\gamma,s_k}=0.
\label{eq:abdc}
\end{align}
\end{subequations}
Here \eqref{eq:abdb} is the contact angle condition, which gives rise to a $90^\circ$ angle between $\gamma^\prime(\vec\nu_{s_k})$ and $\vec n_k$, and \eqref{eq:abdc} is the no-flux boundary condition.

\begin{rem}
For ease of presentation, we consider a single anisotropy $\gamma(\vec p)$ for all the surfaces $\Gamma_i(t)$, $i=1,\ldots, I_S$. Extending the model and the finite element approximation to individual anisotropies $\gamma^{(i)}(\vec p)$,
$i=1,\ldots, I_S$, is straightforward, see e.g.\
\cite{Barrett10cluster, Barrett10}.
We note that in this case choosing $\gamma^{(i)}(\vec p) = \sigma_i\,|\vec p|$
collapses to the isotropic case discussed in Section~\ref{sec:mf}, since then 
$\vec\nu_{\gamma,i} =\gamma'(\vec\nu_i)= \sigma_i\,\vec\nu_i$ and 
$\varkappa_{\gamma,i}=\varkappa_i$ on recalling \eqref{eq:isok}.
\end{rem}

The geometric evolution equations in \eqref{eqn:anivk}, together with the boundary conditions \eqref{eq:tj_a}, \eqref{eqn:atj}, \eqref{eq:bda} and \eqref{eqn:abd}, form a complete model for the evolution of the cluster $\Gamma(t)$ in the case of anisotropic surface energies. The relevant energy is defined by
\begin{align}
A_\gamma(\Gamma(t)):=\sum_{i=1}^{I_S}\int_{\Gamma_i(t)}\gamma(\vec\nu_i)\,\dH^{d-1}.
\end{align}
Analogously to the isotropic case \eqref{eqn:evlaws}, the dynamic system obeys the energy dissipation and volume conservation laws
\begin{subequations} 
\begin{align}
\label{eq:aed}
&\ddt A_\gamma(\Gamma(t)) = -\sum_{i=1}^{I_S}\int_{\Gamma_i(t)}|\nabla_s\varkappa_{\gamma,i}|^2\,\dH^{d-1}\leq 0,\\
&\ddt \vol(\MR_\ell[\Gamma(t)]) =0,\quad \ell = 1,\ldots, I_R.
\label{eq:avc}
\end{align}
\end{subequations}

 To formulate the weak BGN formulation, we introduce some necessary notations from \cite{Barrett08Ani} in the following. For a symmetric positive matrix $G_\bl$, we set $\tG_\bl  = [\rm det\,G_\bl]^{\frac{1}{d-1}}\,[G_\bl]^{-1}$ and define the $\tG_\bl$-inner product
\begin{align}
\bigl(\vec\eta,~\vec\chi\bigr)_{\tG_\bl}= \vec\eta\cdot\tG_\bl\vec\chi,\qquad\forall\vec\eta,~\vec\chi\in\bR^d.\nn
\end{align}
For a smooth scalar field $g$ over $\Gamma_i(t)$, we define the anisotropic surface gradient
\begin{align}
\nabla_s^{\tG_\bl}g = \sum_{j=1}^{d-1}\partial_{\vec t_{j}^{\bl}}g\,\vec t_{j}^{\bl}=\sum_{j=1}^{d-1}(\nabla_sg\cdot\vec t_{j}^{\bl})\,\vec t_{j}^{\bl},
\end{align}
where $\partial_{\vec t_{j}^{\bl}}g = \nabla_sg\cdot\vec t_{j}^{\bl}$ is the directional derivative, $\nabla_s$ is the usual surface gradient operator, and $\{\vec t_{j}^{\bl}\}_{j=1}^{d-1}$ forms an orthonormal basis with respect to the $\tG_\bl$-inner product for the tangent plane of $\Gamma_i(t)$ at the point of interest, i.e., 
\begin{equation*}
\vec t_{j}^{\bl}\cdot\vec\nu_i = 0,\qquad \left(\vec t_{j}^{\bl}, ~\vec t_{k}^{\bl}\right)_{\tG_\bl} = \delta_{jk}, \quad 1\leq j,k\leq d-1,\quad \ell = 1,\ldots, L.
\end{equation*}
Moreover, the anisotropic surface divergence and gradient of a smooth vector field $\vec g$ are given by
\begin{align}
\nabla_s^{\tG_\bl}\cdot\vec g = \sum_{j=1}^{d-1} (\partial_{\vec t_{j}^{\bl}}\vec g)\cdot(\tG_\bl\vec t_{j}^{\bl}),\qquad 
\nabla_s^{\tG_\bl}\vec g =  \sum_{j=1}^{d-1}(\partial_{\vec t_{j}^{\bl}}\vec g)\otimes(\tG_\bl\vec t_{j}^{\bl}),
\end{align}
where $\otimes$ is the stand tensor product for two vectors in $\bR^d$.

Now we present the generalization of \eqref{eqn:weak} to the case of anisotropic surface energies in the form of  \eqref{eq:anipy}. 
Let $\vec x(\cdot,0) \in \VO$, and $\vec x_{s_k}(\partial_{p_k}\Omega_{s_k},0)\subset\mathcal{D}_k$, $k=1,\ldots, I_B$. 
For $t>0$, we find $\vec x(\cdot,t)\in\VO$ such that $(\mathcal{\vv V}(\cdot, t),~\varkappa_\gamma(\cdot, t))\in\Vpartialt\times \Wt$, for $\Gamma(t) = \vec x(\Omega,t)$, with
\begin{subequations}
\label{eqn:aweak}
\begin{align}
\label{eq:aweaka}
\big<\mathcal{\vv V}\cdot\vec\nu,~\chi\big>_{\Gamma(t)}-\big<\nabla_s\varkappa_\gamma,~\nabla_s\chi\big>_{\Gamma(t)}=0\quad\forall\chi\in \Wt,\\[0.4em]
\big<\varkappa_\gamma\,\vec\nu,~\vec\eta\big>_{\Gamma(t)}+\big<\nabla_s^{\tG}\,\vec\id, ~\nabla_s^{\tG}\vec\eta\big>_{\gamma, \Gamma(t)}=0\quad\forall\vec\eta \in\Vpartialt,
\label{eq:aweakb}
\end{align}
\end{subequations}
where we define 
\begin{align}
 \big<\nabla_s^{\tG}\vec \eta, ~\nabla_s^{\tG}\vec\chi\big>_{\gamma, \Gamma(t)}=\sum_{i=1}^{I_S}\sum_{\ell=1}^{L} \int_{\Gamma_i(t)}\left[\frac{\gamma_\bl(\vec\nu_i)}{\gamma(\vec\nu_i)}\right]^{r-1}\left(\nabla_s^{\tG_\bl}\vec\eta,~\nabla_s^{\tG_\bl}\vec\chi\right)_{\tG_\bl}\gamma_\bl(\vec \nu_i)\,\dH^{d-1}.\nn
 \end{align}

\subsection{The generalized SP-PFEM}
Based on the weak formulation \eqref{eqn:aweak} and making use of the discretization in \S\ref{sse:dis}, we can generalize the method \eqref{eqn:fulld} to the case of anisotropic surface energies as follows. Let $\vec X^0\in\VOh$, and $\vec X^0_{s_k}(\partial_{p_k}\Omega_{s_k}^h)\subset\mathcal{D}_k$, $k=1,\ldots, I_B$. For $m=0,\ldots,M-1$, 
find $(\vec X^{m+1}, \kappa_\gamma^{m+1}) \in \VGh \times \Wh$,
with $\vec X^{m+1} - \vec\id\!\mid_{\Gamma^m} \in\Vhpartial$, such that
\begin{subequations}
\label{eqn:afulld}
\begin{align}
\label{eq:afulld1}
\frac{1}{\Delta t_m}\big<\vec{X}^{m+1}-\vec\id, \chi\,\vec\nu^{m+\frac{1}{2}}\big>_{\Gamma^m}^h
- \big<\nabla_s\,\kappa^{m+1}_\gamma, 
\nabla_s\, \chi \big>_{\Gamma^m}  = 0
\quad \forall \chi \in \Wh, \\
\big< \kappa^{m+1}_\gamma\,\vec\nu^{m+\frac{1}{2}}, \vec\eta\big>_{\Gamma^m}^h
+ \big< \nabla_s^{\tG}\, \vec{X}^{m+1}, \nabla_s^{\tG}\, \vec\eta \big>_{\gamma,\Gamma^m} 
 = 0
\quad \forall\ \vec\eta \in \Vhpartial,
\label{eq:afulld2}
\end{align}
\end{subequations}
where we define the discrete inner product $\langle \nabla_s^{\tG}\, \cdot, \nabla_s^{\tG}\, \cdot \rangle_{\gamma,m}$ via
\begin{equation}
\big< \nabla_s^{\tG}\, \vec\eta, ~\nabla_s^{\tG}\, \vec\chi \big>_{\gamma,\Gamma^m} :=
 \sum_{i=1}^{I_S} \sum_{\ell=1}^{L} \int_{\Gamma^m_i} \left[
\frac{\gamma_\bl(\vec\nu_i^{m+1})}{\gamma(\vec\nu^{m+1}_i)} 
\right]^{r-1}\!\!\!\!
\bigl(\nabla_s^{\tG_\bl}\,\vec\eta_i,\nabla_s^{\tG_\bl}\,
\vec\chi_i\bigr)_{\tG_i^{\ell}} \,\gamma_\bl(\vec\nu^m_i) \dH^{d-1} \,.
\label{eq:dip}
\end{equation}
The above scheme \eqref{eqn:afulld} is very similar to \cite[(4.9)]{Barrett10cluster} except that we apply a semi-implicit approximation of the unit normal in the first terms of \eqref{eq:afulld1} and \eqref{eq:afulld2}.  
That means in the case $r=1$ the scheme \eqref{eqn:afulld} introduces a
nonlinearity compared to the linear scheme \cite[(4.9)]{Barrett10cluster}.
But for $r\not=1$ the introduced nonlinearity is mild compared to the
dependence of \eqref{eq:dip} on the unit normal $\vec\nu^{m+1}$ on $\Gamma^{m+1}$, which is necessary in order to prove unconditional stability \cite{Barrett08Ani}.

We first present a lemma which will be used to prove the unconditional stability for the discretized scheme in \eqref{eqn:afulld}, and its proof can be found in \cite[Lemma 3.1]{Barrett08Ani}.
\begin{lem}\label{lem:ani} 
Let $\vec X^{m+1}\in\VGh$ with $\vec X^{m+1} - \vec\id\!\mid_{\Gamma^m}
\in\Vhpartial$. 
Then it holds
\begin{align}
\sum_{\ell=1}^{L} \int_{\Gamma^m_i} \left[
\frac{\gamma_{\ell}(\vec\nu_i^{m+1})}{\gamma(\vec \nu^{m+1}_i)} 
\right]^{r-1}\!\!\!\!
\bigl(\nabla_s^{\tG_\bl}\,\vec{X}_i^{m+1},~\nabla_s^{\tG_\bl}
(\vec X_i^{m+1} - \vec\id)\bigr)_{\tG_\bl} \,\gamma_\bl(\vec\nu^m_i) \dH^{d-1}\nn\\ \geq \int_{\Gamma_i^{m+1}}\gamma(\vec\nu_i^{m+1})\,\dH^{d-1} - \int_{\Gamma_i^m}\gamma(\vec\nu_i^m)\,\dH^{d-1},\nn
\end{align}
which yields the following inequality on recalling \eqref{eq:dip}
\begin{align}
\big<\nabla_s^{\tG}\,\vec X^{m+1}, ~\nabla_s^{\tG}\,(\vec X^{m+1} - \vec\id)\big>_{\gamma,\Gamma^m}\geq \sum_{i=1}^{I_S}\int_{\Gamma_i^{m+1}}\gamma(\vec\nu_i^{m+1})\,\dH^{d-1} - \sum_{i=1}^{I_S}\int_{\Gamma_i^{m}}\gamma(\vec\nu_i^{m})\,\dH^{d-1}.\nn
\end{align}
\end{lem}

For the discretized scheme in \eqref{eqn:afulld}, we can prove the unconditional energy decay and the conservation  of volume for each enclosed bubble.

\begin{thm}[stability and volume conservation]\label{thm:ani}Let $(\vec X^{m+1},~\kappa_\gamma^{m+1})$ be a solution to \eqref{eqn:afulld}, then it holds that
\begin{align}
\label{eq:athE}
A_\gamma(\Gamma^{m+1})+\Delta t_m\big<\nabla_s\kappa_\gamma^{m+1},~\nabla_s\kappa_\gamma^{m+1}\big>_{\Gamma^m}\leq A_\gamma(\Gamma^m).
\end{align}
Moreover, it holds that
\begin{align}
\label{eq:athV}
{\rm vol}(\MR_\ell[\Gamma^{m+1}]) = {\rm vol}(\MR_\ell[\Gamma^m]),\quad \ell = 1,\ldots, I_R.
\end{align}
\end{thm}
\begin{proof}
Setting $\chi = \Delta t_m\,\kappa_\gamma^{m+1}$ in \eqref{eq:afulld1} and $\eta = \vec X^{m+1} - \vec\id\!\mid_{\Gamma^m}$ in \eqref{eq:afulld2} and combining the two equations yields
\begin{align}
\Delta t_m\big<\nabla_s\,\kappa^{m+1}_\gamma, 
\nabla_s\, \kappa_\gamma^{m+1} \big>_{\Gamma^m} + \big<\nabla_s^{\tG}\,\vec X^{m+1}, ~\nabla_s^{\tG}\,(\vec X^{m+1} - \vec\id)\big>_{\gamma,\Gamma^m}=0.\nn
\end{align}
On recalling Lemma \ref{lem:ani}, we directly obtain the unconditional stability in \eqref{eq:athE} as claimed. 

Finally, in \eqref{eq:afulld2}, we choose $\chi$ with $\chi_i$ satisfying \eqref{eq:xior}. This yields \eqref{eq:athV} by Lemma \ref{lem:vc}. 
\end{proof}

\begin{rem}\label{rem:iterativemethod}Like in the isotropic case, we can solve the nonlinear system resulting from \eqref{eqn:afulld} with a lagged Picard-type iteration as follows.  For each $p\geq 0$, find $(\vec X^{m+1, p+1}, \kappa_\gamma^{m+1, p+1}) \in \VGh \times \Wh$,
with $\vec X^{m+1, p+1} - \vec\id\!\mid_{\Gamma^m} \in \Vhpartial$, such that for all $(\chi,\vec\eta)\in\Wh\times\Vhpartial$ the following two equations hold
\begin{subequations}
\label{eqn:apicard}
\begin{align}
\label{eq:apicard1}
&\frac{1}{\Delta t_m}\big<\vec{X}^{m+1, p+1}-\vec\id, \chi\,\vec\nu^{m+\frac{1}{2},p}\big>_{\Gamma^m}^h
-\big<\nabla_s\,\kappa^{m+1,p+1}_\gamma, 
\nabla_s\, \chi \big>_{\Gamma^m}=0, \\
&\big< \kappa^{m+1, p+1}_\gamma\,\vec\nu^{m+\frac{1}{2}, p}, \vec\eta\big>_{\Gamma^m}^h
+ \sum_{i=1}^{I_S} \sum_{\ell=1}^{L} \int_{\Gamma^m_i} \left[
\frac{\gamma_\bl(\vec\nu_i^{m+1,p})}{\gamma(\vec\nu^{m+1, p}_i)} 
\right]^{r-1}\!\!\!\!
\left(\nabla_s^{\tG_\bl}\,\vec X_i^{m+1,p+1},\nabla_s^{\tG_\bl}\,
\vec\eta_i\right)_{\tG_\bl} \,\gamma_\bl(\vec\nu^m_i) \dH^{d-1}
 = 0,\label{eq:apicard2}
\end{align}
\end{subequations}
where we denote $\Gamma^{m+1,p} = \vec X^{m+1,p}(\Gamma^m)$, and $\vec\nu^{m+1,p}$ and $\vec\nu^{m+\frac{1}{2},p}$ are defined by using the similar formulas in \eqref{eq:vG} and \eqref{eq:weightv} except that $\Gamma^{m+1}$ is replaced by $\Gamma^{m+1,p}$ instead. In particular, we choose $\vec X^{m+1,0} = \vec\id\!\mid_{\Gamma^m}$. The resulting linear system from \eqref{eqn:apicard} can then be solved efficiently with the Schur complement approaches in BGN.
\end{rem}

\setcounter{equation}{0} 
\section{Extension to non-neutral external boundaries} \label{sec:exteneb}

So far, for ease of presentation, we have only considered  the simplified case when the contact energy densities, for the two phases separated by the interface at the external boundary, are the same, so that they have no contribution to the total energy of the system. As suggested by \eqref{eq:abdb}, this then leads to a $90^\circ$ angle between $\gamma^\prime(\vec\nu_{s_k})$ and $\vec n_k$. However, in practical physical applications, this is usually not the case and the contact energies play a non-negligible role in the evolution of the surface cluster. To this end, we consider the dynamic system in \S\ref{sec:amathfor} but replace the contact angle condition \eqref{eq:abdb}  with the following anisotropic Young's equation \cite{Barrett10}
\begin{align}
\vec n_k\cdot\gamma^\prime (\vec\nu_{s_k})=\varrho_k,\quad k = 1,\ldots, I_B,
\label{eq:aniy}
\end{align}
which gives rise to more general contact angles. Here $\varrho_k$ are given constants and represent the change of contact energy density in the direction of $-\vec \nu_{s_k}$, that the two phases separated by the surface $\Gamma_{s_k}$ have with the external boundary $\mathcal{D}_k$. A similar contact angle condition has also been derived in \cite{JiangZB20}. It is easy to see that \eqref{eq:aniy} yields an angle of $\arccos\frac{\varrho_k}{|\gamma^\prime(\vec\nu_{s_k})|}$ between $\gamma^\prime(\vec \nu_{s_k})$ and $\vec n_{k}$ when $|\rho_k|\leq |\gamma^\prime(\vec \nu_{s_k})|$. 
In particular, in the isotropic case we obtain a contact angle $\vartheta_k$ 
with $\cos\vartheta=\varrho_k$, for $\varrho_k\in[-1,1]$.

\begin{figure}[!htp]
\centering
\includegraphics[width=0.8\textwidth]{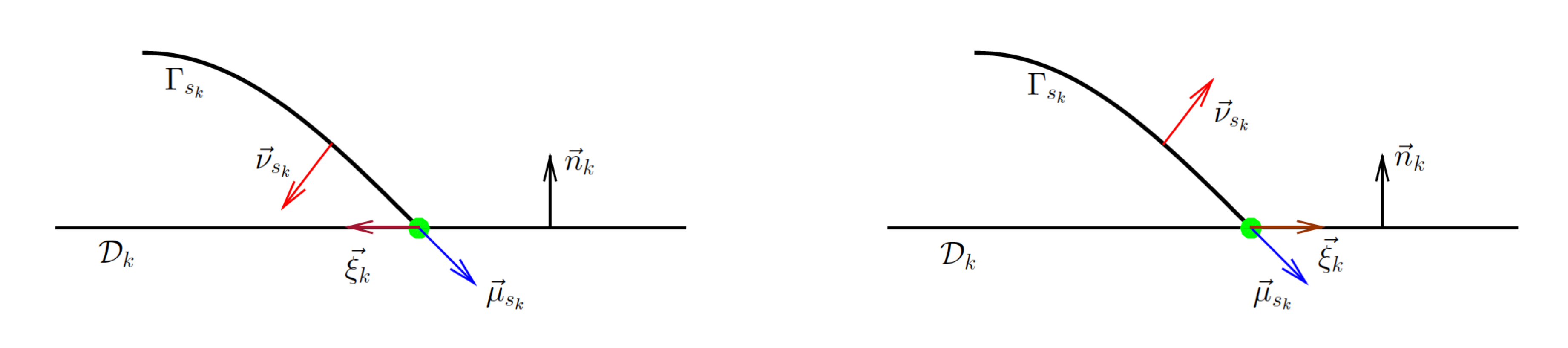}
\caption{Sketch of the structures at the boundary line $\mathcal{B}_k$ ({\green\bf green}) where $\Gamma_{s_k}$ meets the external planar boundary $\mathcal{D}_k$. Depicted above is a plane that is perpendicular to
$\mathcal{B}_k$. }
\label{fig:bl}
\end{figure}

We now discuss the contact energy contributions to the  system. At $\mathcal{B}_k$, we define 
\begin{align}
\vec\xi_k = (\vec n_k\cdot\vec \nu_{s_k})\,\vec\mu_{s_k} - (\vec n_k\cdot\vec \mu_{s_k})\,\vec\nu_{s_k},\quad k = 1,\ldots, I_B,
\label{eq:xik}
\end{align}
where we observe that $\vec\xi_k$ is normal to $\mathcal{B}_k$ and lies in the tangent plane of the surface $\mathcal{D}_k$. In particular, $\vec\xi_k$ is obtained through a $90^\circ$ rotation of $\vec n_k$ in the plane spanned by $\vec\nu_{s_k}$ and $\vec\mu_{s_k}$, and that $(\vec n_{k},\vec\xi_{k})$ have the same orientation with $(\vec\nu_{s_k},~\vec\mu_{s_k})$, as shown in Fig.~\ref{fig:bl}. Let $\mathbf{B}_R^d$ be a ball in $\bR^d$ with sufficiently large radius $R$, and for $k=1,\ldots, I_B$ we set $\mathscr{G}_k=\mathcal{D}_k\cap\mathbf{B}_R^d$. 
Then the boundary point/line $\mathcal{B}_k$ divides the segment/disk $\mathscr{G}_k$ into two parts by
\begin{align}
\overline{\mathscr{G}_k^+}\cap \overline{\mathscr{G}_k^-} = \mathcal{B}_k,\qquad \overline{\mathscr{G}_k^+}\cup\overline{\mathscr{G}_k^-}=\mathscr{G}_k,
\end{align}
where $\mathscr{G}_k^-$ is chosen such that $\vec\xi_k$ in \eqref{eq:xik} is the outer normal to $\mathscr{G}_k^-$ on $\mathcal{B}_k$. The relevant energy of the considered system is then given by
\begin{align}
E(\Gamma(t)) &= A_\gamma(\Gamma(t)) + A_\partial(\Gamma(t))= \sum_{i=1}^{I_S}\int_{\Gamma_i(t)}\gamma(\vec\nu_i)\,\dH^{d-1} + \sum_{k=1}^{I_B}\left(\widehat{\varrho^+_k}\,|\mathscr{G}_k^+(t)| + \widehat{\varrho^-_k}\,|\mathscr{G}_k^-(t)|\right),
\end{align}
where $A_\partial(\Gamma(t))$ represents the contact energies, $\widehat{\varrho_k^{\pm}}$ are the contact energy densities of the plane surfaces $\mathscr{G}_k^{\pm}(t)$ which satisfy the relation $\widehat{\varrho_k^+}- \widehat{\varrho_k^-}=\varrho_k$, and $|\mathscr{G}^\pm_k|$ represent the surface area of $\mathscr{G}_k^\pm$, respectively. Direct calculation yields the energy dissipation law (see \cite[Proposition 2.1]{Barrett10}): 
\begin{align}
\ddt E(\Gamma(t)) + \sum_{i=1}^{I_S}\int_{\Gamma_i(t)}|\nabla_s\varkappa_{\gamma,i}|^2\,\dH^{d-1}= 0.
\end{align}
In addition, we still have the volume conservation law \eqref{eq:avc}.

We then generalize the weak formulation in \eqref{eqn:aweak} to the case of non-neutral external boundaries. In order that \eqref{eq:aniy} can be weakly enforced, we add the following terms on the right hand side of \eqref{eq:aweakb}
\begin{align}
\label{eq:bterm}
\sum_{k=1}^{I_B}\varrho_k\int_{\mathcal{B}_k(t)}\vec\xi_k\cdot\vec\eta_{s_k}\,\dH^{d-2}.
\end{align}
Similarly, we generalize the discretized numerical method \eqref{eqn:afulld} as follows. On the right hand of \eqref{eq:afulld2}, we add 
\begin{align}
\sum_{k=1}^{I_B}\varrho_k\int_{\mathcal{B}_k^m}\vec\xi_k^{m+\frac{1}{2}}\cdot\vec\eta_{s_k}\,\dH^{d-2},
\end{align}
where $\mathcal{B}_k^m$ is the natural discrete analogue of $\mathcal{B}_k(t_m)$ and $\vec\xi_k^{m+\frac{1}{2}}$ is an appropriate approximation in order to  guarantee the unconditional stability for the generalized scheme.

 Following \cite{BaoZ20preprint}, we next discuss the treatment of $\vec\xi_k^{m+\frac{1}{2}}$ in detail. In the case of $d=2$, $\vec\xi_k^{m+\frac{1}{2}}$ can be simply determined from $\vec n_k$ via a $90^\circ$ rotation in $\bR^2$. While in the case of $d=3$,  we have $\vec\xi_k = \vec n_k\times(\vec\mu_{s_k}\times \vec\nu_{s_k})$ by \eqref{eq:xik}. We assume that  $\{\vec\lambda^{m,k}_\ell\}_{\ell=0}^{Y_k}$ is an ordered sequence of vertices of $\mathcal{B}_k^m$ according to the direction of $\vec\mu_{s_k}\times \vec\nu_{s_k}$ and denote 
\begin{align}
\mathcal{B}_k^m=\bigcup_{\ell=1}^{Y_k}\overline{L_{\ell,k}^m}=\bigcup_{\ell=1}^{Y_k}[\vec\lambda_{\ell-1}^{m,k},~\vec\lambda_\ell^{m,k}],\qquad \vec f\{L_{\ell,k}^m\}=\vec\lambda_{\ell}^{m,k} - \vec\lambda_{\ell-1}^{m,k},\nn
\end{align}
where $L_{\ell,k}^m$ is the $\ell$th line segment of $\mathcal{B}_k^m$ and $\vec f\{L_{\ell,k}^m\}$ represents its orientation vector. Based on \eqref{eq:qh}, we can naturally set $\mathcal{B}_k^h(t)=\bigcup_{\ell=1}^{Y_k}\overline{L_{\ell,k}^h(t)}=\bigcup_{\ell=1}^{Y_k}[\vec\lambda_{\ell-1}^{h,k}(t),~\vec\lambda_\ell^{h,k}(t)]$ as a linear interpolation between $\mathcal{B}_k^m$ and $\mathcal{B}_k^{m+1}$, and
\begin{align}
\vec \lambda_\ell^{h,k}(t) &=\frac{t_{m+1}-t}{\Delta t_m}\vec \lambda_\ell^{m,k} + \frac{t - t_m}{\Delta t_m}\vec \lambda_\ell^{m+1,k},\quad t\in[t_m,~t_{m+1}], \quad \ell = 0,\ldots, Y_k.\label{eq:BLpoints}
\end{align}
We then define $\vec\xi_k^{m+\frac{1}{2}}$ in an average sense via 
\begin{align}
\vec\xi_{k}^{m+\frac{1}{2}}|_{L_{\ell,k}^m}&=\vec\xi_{k,\ell}^{m+\frac{1}{2}} =\vec n_k\times\left(\frac{1}{\Delta t_m\,|\vec f\{L_{\ell,k}^m\}|}\int_{t_m}^{t_{m+1}}\,\vec f\{L_{\ell,k}^h(t)\}\rd t\right)\nn\\
&=\frac{1}{2|\vec f\{L_{\ell,k}^m\}|}\,\vec n_k\times \left(\vec f\{L_{\ell,k}^m\}+\vec f\{L_{\ell,k}^{m+1}\} \right),\quad k = 1,\ldots, I_B,\quad \ell = 1,\ldots, Y_k.
\label{eq:xiweight}
\end{align}

From \cite[Lemma~3.1]{BaoZ20preprint}, we have the following lemma for $\vec\xi^{m+\frac{1}{2}}_k$. For completeness, here we present a new proof in a similar manner as we did in the proof of Lemma~\ref{lem:vc}. 
\begin{lem}\label{lem:bd} 
Let $\vec X^{m+1}\in\VGh$ with $\vec X^{m+1} - \vec\id\!\mid_{\Gamma^m}
\in\Vhpartial$. Then it holds that
\begin{align} \label{eq:bd}
|\mathscr{G}_k^{m+1,\pm}| - |\mathscr{G}_k^{m,\pm}| = \mp\int_{\mathcal{B}_k^m}\vec \xi_k^{m+\frac{1}{2}}\cdot(\vec X^{m+1}_{s_k}-\vec\id)\,\dH^{d-2},\quad k = 1,\ldots, I_B,
\end{align}
where $\mathscr{G}_k^{m,\pm}$ are the natural discrete analogues of $\mathscr{G}_k^{\pm}(t_m)$. 
\end{lem} 
\begin{proof}
In the case of $d=2$ we have that $\mathcal{B}_k^m$ and 
$\vec X^{m+1}_{s_k}(\mathcal{B}_k^m)$ are points on the line $\mathcal{D}_k$, 
while $\mathscr{G}_k^{m,\pm}$ are line segments on $\mathcal{D}_k$, meaning
the result \eqref{eq:bd} is elementary.

In the case of $d=3$, we recall that the vertices of the polygonal curve $\mathcal{B}_k^h(t)$ are given by \eqref{eq:BLpoints}. It is natural to define
\begin{equation}
\vec\xi_k^h(t)|_{L_{\ell,k}^h(t)} = \vec n_k \times \frac{\vec f\{L_{\ell, k}^h(t)\}}{|\vec f\{L_{\ell,k}^h(t)\}|},\qquad \overline{\mathscr{G}_k^{h,-}(t)}\cap\overline{\mathscr{G}_k^{h,+}(t)}=\mathcal{B}_k^h(t),\quad \overline{\mathscr{G}_k^{h,-}(t)}\cup\overline{\mathscr{G}_k^{h,+}(t)} = \mathcal{D}_k\cap\mathbf{B}_R^3,\nn
\end{equation}
such that $\vec{\xi}_k^{h}(t)$ is the outer normal to $\mathscr{G}_k^{h,-}(t)$. Applying the Reynolds transport theorem to the two-dimensional domain $\mathscr{G}_k^{h,-}(t)$ gives
\begin{align}
\frac{\rd}{\rd t}|\mathscr{G}_k^{h,-}(t)|&= \int_{\mathcal{B}_k^h(t)}\vec\xi_k^h(t)\cdot(\partial_t\vec X_{s_k}^h) \circ(\vec X_{s_k}^h)^{-1}\,\dH^{1}\nn\\
& = \sum_{\ell = 1}^{Y_k}\int_{L_{\ell,k}^m}\left(\vec n_k \times \frac{\vec f\{L_{\ell, k}^h(t)\}}{|\vec f\{L_{\ell,k}^h(t)\}|}\right)\cdot\frac{\vec X_{s_k}^{m+1} - \vec\id}{\ttau_m}\,\frac{|\vec f\{L_{\ell,k}^h(t)\}|}{|\vec f\{L_{\ell,k}^m\}|}\,\dH^{1}\nn\\
&=\sum_{\ell = 1}^{Y_k}\int_{L_{\ell,k}^m}\left(\vec n_k\times\frac{\vec f\{L_{\ell, k}^h(t)\}}{\ttau_m\,|\vec f\{L_{\ell,k}^m\}|}\right)\cdot(\vec X_{s_k}^{m+1} - \vec\id)\,\dH^{1},
\label{eq:bdac}
\end{align}
where $\vec X_{s_k}^h$ is defined in \eqref{eq:xh}. Integrating \eqref{eq:bdac} from $t_m$ to $t_{m+1}$ with respect to $t$ yields 
\begin{equation}
|\mathscr{G}_k^{m+1,-}| - |\mathscr{G}_k^{m,-}| = \sum_{\ell = 1}^{Y_k}\int_{L_{\ell,k}^m}\vec\xi_{k,\ell}^{m+\frac{1}{2}}\cdot(\vec X_{s_k}^{m+1} - \vec\id)\dH^{1} = \int_{\mathcal{B}_k^m}\vec\xi_k^{m+\frac{1}{2}}\cdot(\vec X_{s_k}^{m+1} - \vec\id)\,\dH^{1}
\end{equation}
on recalling \eqref{eq:xiweight}. Using a similar approach to $\mathscr{G}_k^{h,+}(t)$ yields that 
\begin{equation}
|\mathscr{G}_k^{m+1,+}| - |\mathscr{G}_k^{m,+}| = -\int_{\mathcal{B}_k^m}\vec\xi_{k}^{m+\frac{1}{2}}\cdot(\vec X_{s_k}^{m+1} - \vec\id)\,\dH^{1}.
\end{equation}
Thus we obtain \eqref{eq:bd}.
\end{proof}

We then have the following theorem which generalizes Theorem \ref{thm:ani}.  

\begin{thm}
Let $(\vec X^{m+1},~\kappa_\gamma^{m+1})$ be a solution to \eqref{eqn:afulld} with \eqref{eq:bterm} added to the right hand side of \eqref{eq:afulld2}. 
Then it holds that
\begin{align}
\label{eq:athEb}
E(\Gamma^{m+1})+\Delta t_m\big<\nabla_s\kappa_\gamma^{m+1},~\nabla_s\kappa_\gamma^{m+1}\big>_{\Gamma^m}\leq E(\Gamma^m).
\end{align}
Moreover, it holds that
\begin{align}
\label{eq:athVb}
{\rm vol}(\MR_\ell[\Gamma^{m+1}]) = {\rm vol}(\MR_\ell[\Gamma^m]),\quad \ell = 1,\ldots, I_R.
\end{align}
\end{thm}
\begin{proof}
Setting $\chi = \Delta t_m\,\kappa_\gamma^{m+1}$ in \eqref{eq:afulld1} and $\eta = \vec X^{m+1} - \vec\id\!\mid_{\Gamma^m}$ in the adapted \eqref{eq:afulld2} and combining the two equations yields
\begin{align}
\Delta t_m\big<\nabla_s\,\kappa^{m+1}_\gamma, 
\nabla_s\, \kappa_\gamma^{m+1} \big>_{\Gamma^m} + \big<\nabla_s^{\tG}\,\vec X^{m+1}, ~\nabla_s^{\tG}\,(\vec X^{m+1} - \vec\id)\big>_{\gamma,\Gamma^m}
=\sum_{k=1}^{I_B}\varrho_k\int_{\mathcal{B}_k^m}\vec\xi_k^{m+\frac{1}{2}}\cdot(\vec X^{m+1}_{s_k}-\vec\id)\,\dH^{d-2}.\nn
\end{align}
By Lemma~\ref{lem:bd}, and on noting 
$\widehat{\varrho_k^+}- \widehat{\varrho_k^-}=\varrho_k$, we have
\begin{align}
\widehat{\varrho_k^+}(|\mathscr{G}_k^{m+1,+}|-|\mathscr{G}_k^{m,+}|) + \widehat{\varrho_k^-}(|\mathscr{G}_k^{m+1,-}|-|\mathscr{G}_k^{m,-}|) = - \sum_{k=1}^{I_B}\varrho_k\int_{\mathcal{B}_k^m}\vec\xi_k^{m+\frac{1}{2}}\cdot(\vec X^{m+1}_{s_k}-\vec\id)\,\dH^{d-2},\nn
\end{align}
which yields \eqref{eq:athEb} on recalling Lemma \ref{lem:ani}. Finally, \eqref{eq:athVb} follows directly by choosing $\chi$ in \eqref{eq:afulld1} with $\chi_i$ satisfying \eqref{eq:xior}. 
\end{proof}

\setcounter{equation}{0}
\section{Numerical results} \label{sec:num}

We implemented our fully discrete finite element approximations 
within the finite element toolbox ALBERTA, see \cite{Alberta}. The systems of
linear equations arising from the Picard-iteration are solved with the help
of the Schur complement approach from BGN, 
employing a preconditioned conjugate gradient solver with preconditioners based
on the sparse factorization package UMFPACK, see \cite{Davis04}.

Throughout this section we use uniform time steps $\ttau_m=\ttau$. 
We let $J = \sum_{i=1}^{I_S} J_i$ denote the total number of elements, 
and $K = \sum_{i=1}^{I_S} K_i$ the total number of vertices.
Unless otherwise stated, we use $\varrho_k=\varrho$ for $k=1,\ldots, I_B$,
with $\rho = 0$ by default. For many of the presented simulations 
we will put particular emphasis on the volume preserving aspect. 
Hence, for later use we 
define the relative volume error at time $t=t_m$ as 
\[
v_\Delta^m = \max_{\ell=1,\ldots,I_R} \left|\frac{\vol(\MR_\ell[\Gamma^m]) - \vol(\MR_\ell[\Gamma^0])}
{\vol(\MR_\ell[\Gamma^0])}\right|.
\]
We also define the mesh ratio 
\begin{equation} \label{eq:ratio}
r^m = \max_{i=1,\ldots,I_S}
\dfrac{\max_{j=1,\ldots,J_i} |\sigma^{m,i}_j|}
{\min_{j=1,\ldots,J_i} |\sigma^{m,i}_j|}.
\end{equation}
Throughout we use solid red lines for the introduced structure-preserving schemes, and dashed blue lines for the standard BGN scheme.
We stress that all the presented numerical simulations were performed without
any mesh smoothings or remeshings.

\subsection{Numerical results in 2d} 

\begin{figure}[tbh]
\centering
\includegraphics[width=0.90\textwidth]{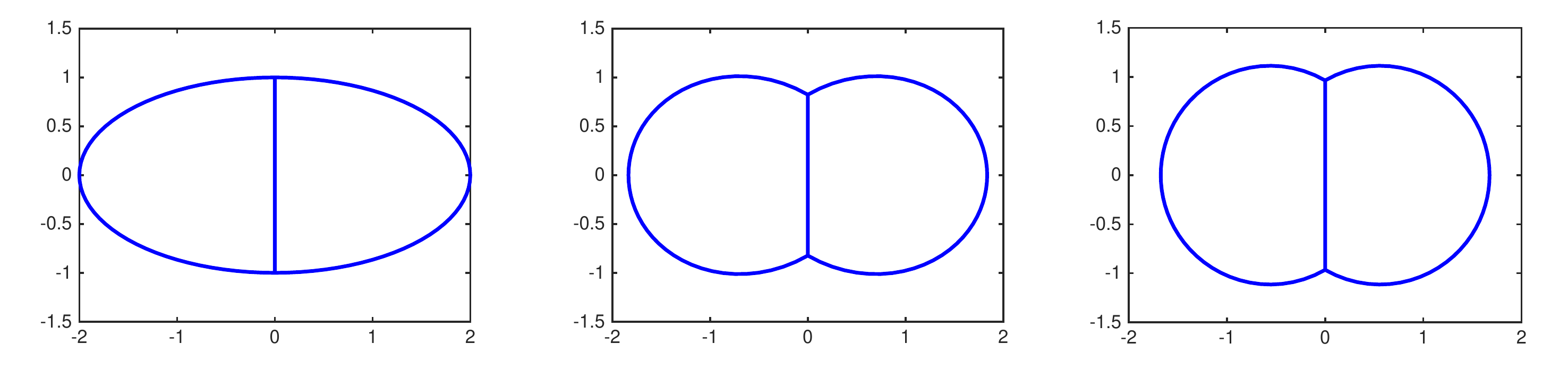}
\includegraphics[width=0.90\textwidth]{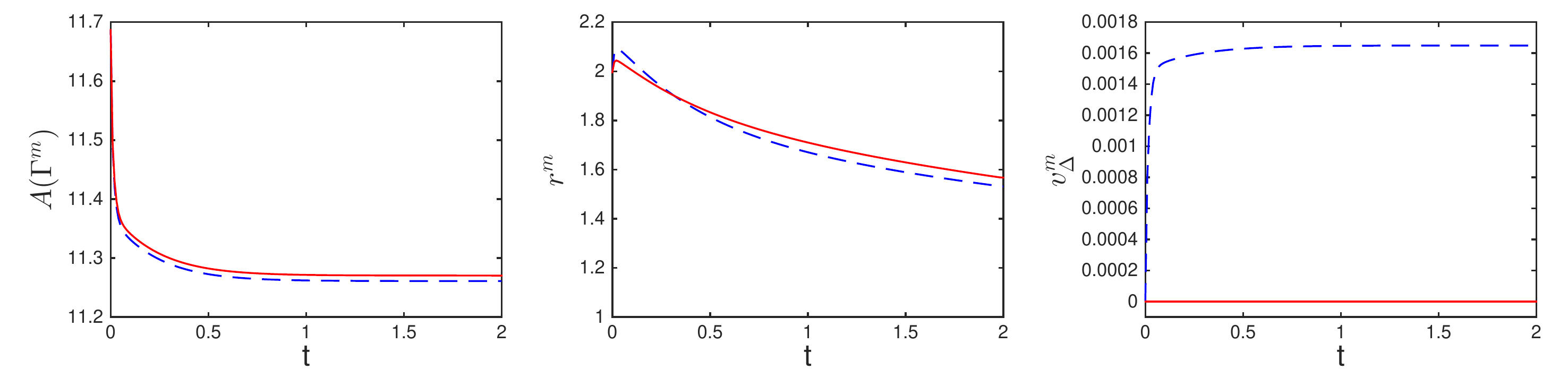}
\caption{
Evolution towards the 2d standard double bubble.
Plots of $\Gamma^m$ at times $t=0, 0.1, 2$.
We also show plots of the discrete energy $A(\Gamma^m)$, the ratio $r^m$ and 
the relative volume error $v_\Delta^m$ over time, where  $\sigma = (1,1,1)$, $K = 129$ and $\ttau =10^{-2}$.
}
% ~/hpc_cluster/data/alberta/tjtrue/2d.db
\label{fig:2ddbw1}
\end{figure}%

We start with the evolution of a curve network towards the well-known double bubble minimizer. The initial network is
given by two $2:1$ semi-ellipses and a straight line, 
meeting at two triple junction points. The discretization
parameters are chosen as $K = 129$ and $\ttau = 10^{-2}$. In the first simulation, we consider the standard double bubble with equal surface energy densities $\sigma = (1,1,1)$. The numerical results are shown in Fig.~\ref{fig:2ddbw1}, where we observe that triple junction angles approach $120^\circ$ in the steady state. Based on the observation, we also find that (i) the volume preservation for the introduced SP-PFEM is well satisfied, as expected, while for the BGN scheme more than $0.15\%$ volume loss is observed. (ii) the mesh ratios for both schemes remain at small values, which implies the good mesh qualities; and (iii) the energy dissipation shows a good agreement.

\begin{figure}[tbh]
\centering
\includegraphics[width=0.90\textwidth]{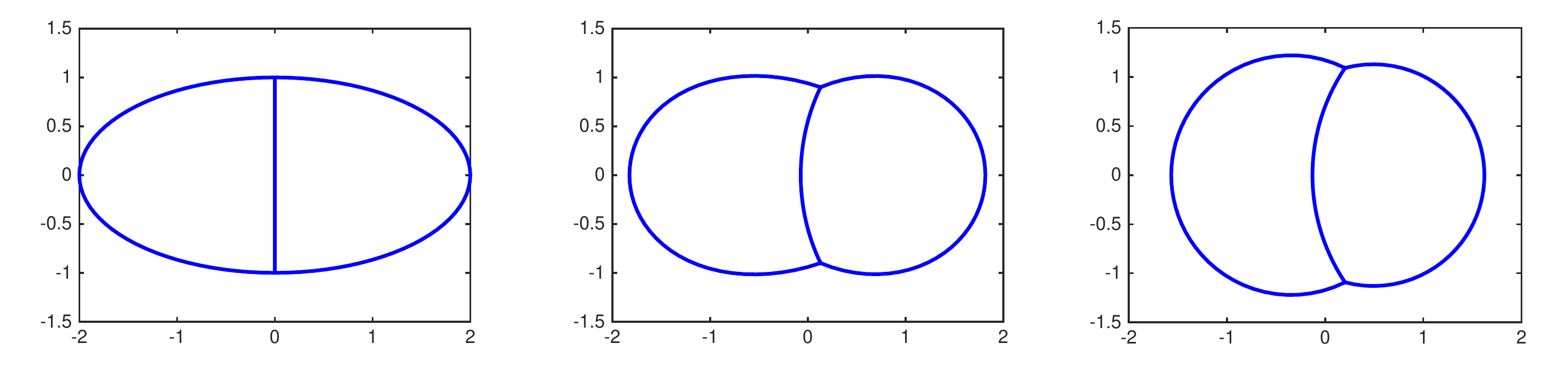}
\includegraphics[width=0.90\textwidth]{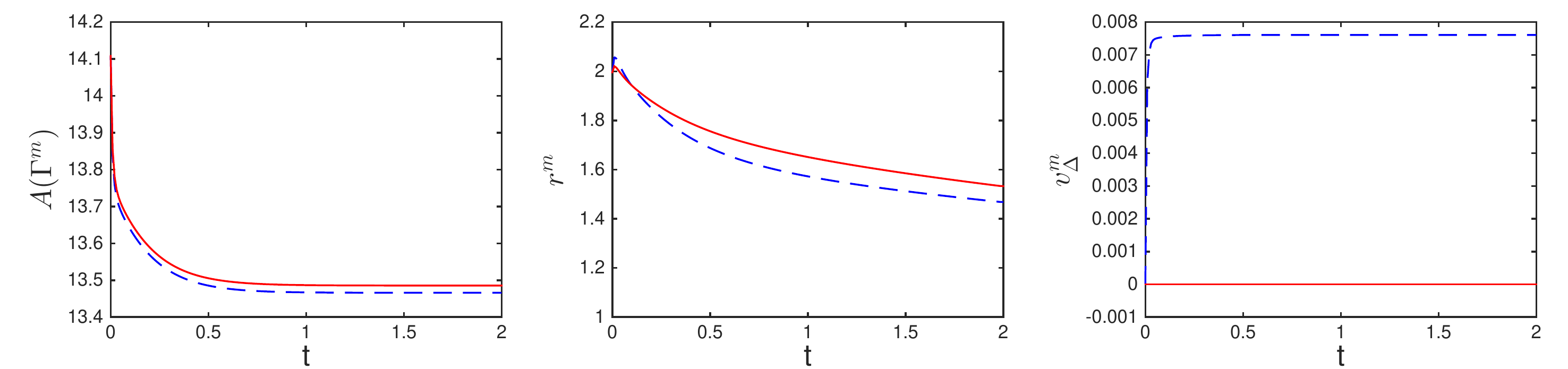}
\caption{
Evolution towards a 2d double bubble, with weightings $\sigma = (1,1,1.5)$.
Plots of $\Gamma^m$ at times $t=0, 0.1, 2$.
We also show plots of the discrete energy $A(\Gamma^m)$, the ratio $r^m$ and 
the relative volume error $v_\Delta^m$ over time, where $K = 129$ and $\ttau=10^{-2}$.
}
% ~/hpc_cluster/data/alberta/tjtrue/2d.db_weight1.5
\label{fig:2ddbw1.5}
\end{figure}

\begin{figure}[!tbh]
\center
\includegraphics[width=0.90\textwidth]{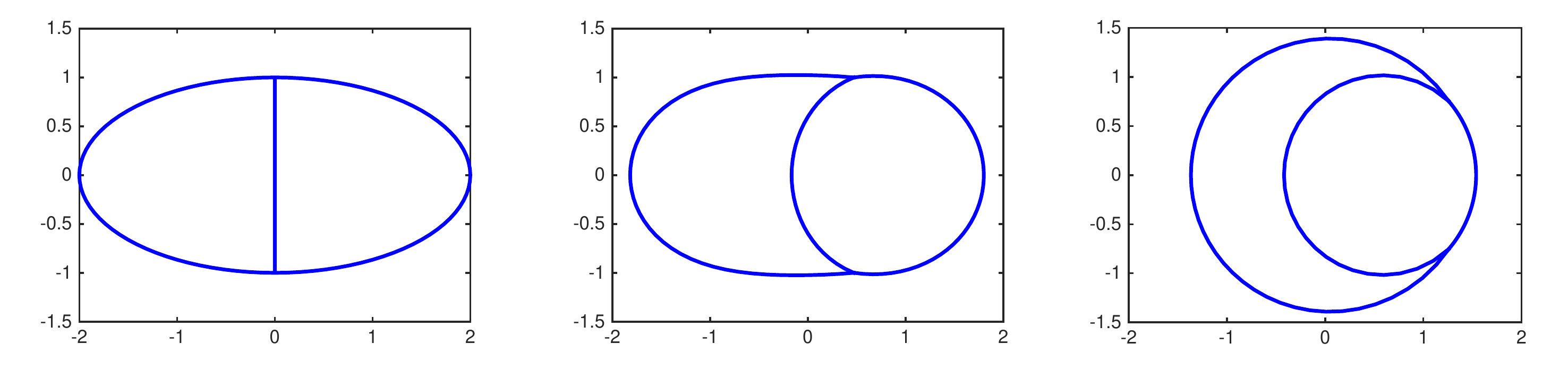}
\includegraphics[width=0.90\textwidth]{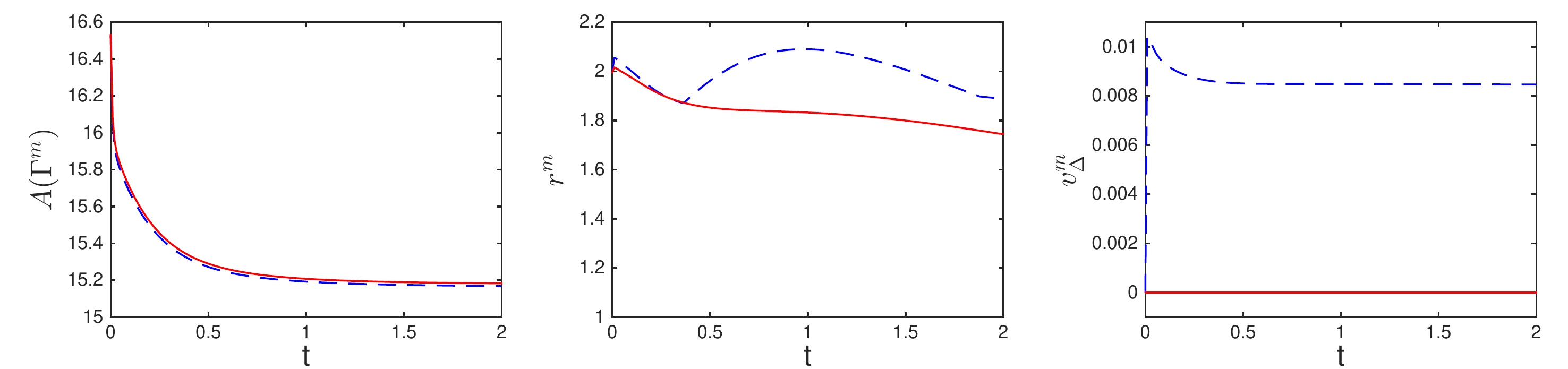}
\caption{
Evolution towards a degenerate 2d double bubble, with weightings $\sigma = (1,1,2)$.
Plots of $\Gamma^m$ at times $t=0, 0.1, 2$.
We also show plots of the discrete energy $A(\Gamma^m)$, the ratio $r^m$ and 
the relative volume error $v_\Delta^m$ over time, where $K = 129$ and $\ttau=10^{-2}$.
}
% ~/hpc_cluster/data/alberta/tjtrue/2d.db_weight2
\label{fig:2ddbw2}
\end{figure}

\begin{figure}[!htp]
\centering
\includegraphics[width=0.90\textwidth]{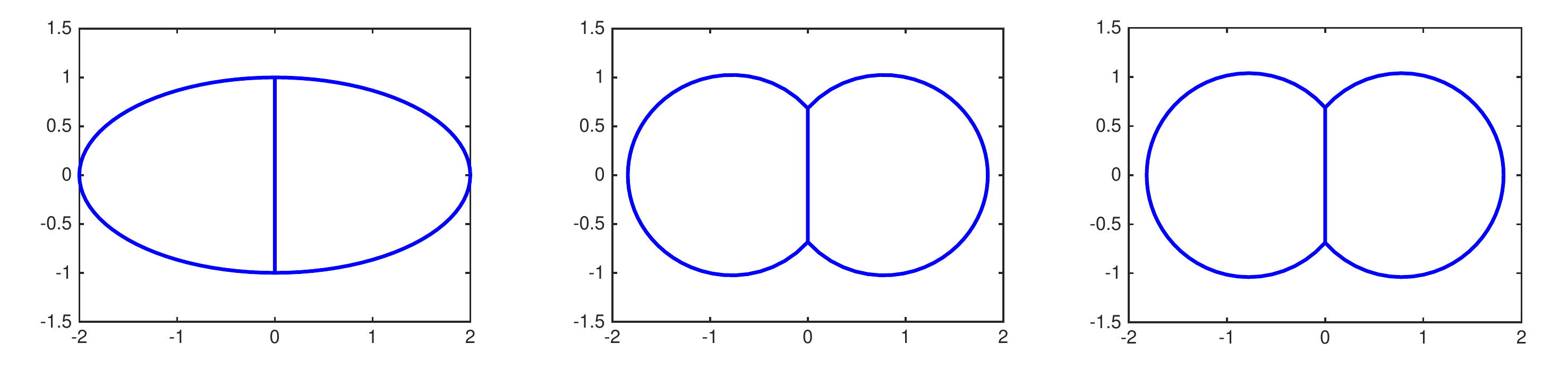}
\includegraphics[width=0.90\textwidth]{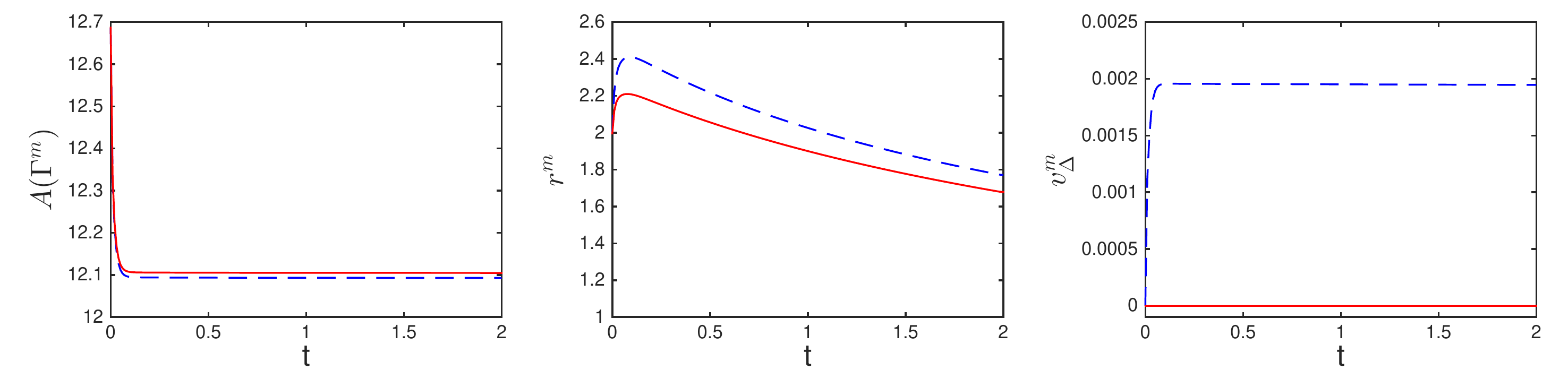}\hspace{0.4cm}
\caption{
Evolution towards a 2d double bubble, with weightings $\sigma = (1,1.5,1)$.
Plots of $\Gamma^m$ at times $t=0, 0.1, 2$.
We also show plots of the discrete energy $A(\Gamma^m)$, the ratio $r^m$ and 
the relative volume error $v_\Delta^m$ over time, where $K = 129$ and $\ttau=10^{-2}$.
}
% ~/hpc_cluster/data/alberta/tjtrue/2d.db_1weight1.5
\label{fig:2ddb1w1.5}
\end{figure}

\begin{figure}[!htp]
\centering
\includegraphics[width=0.90\textwidth]{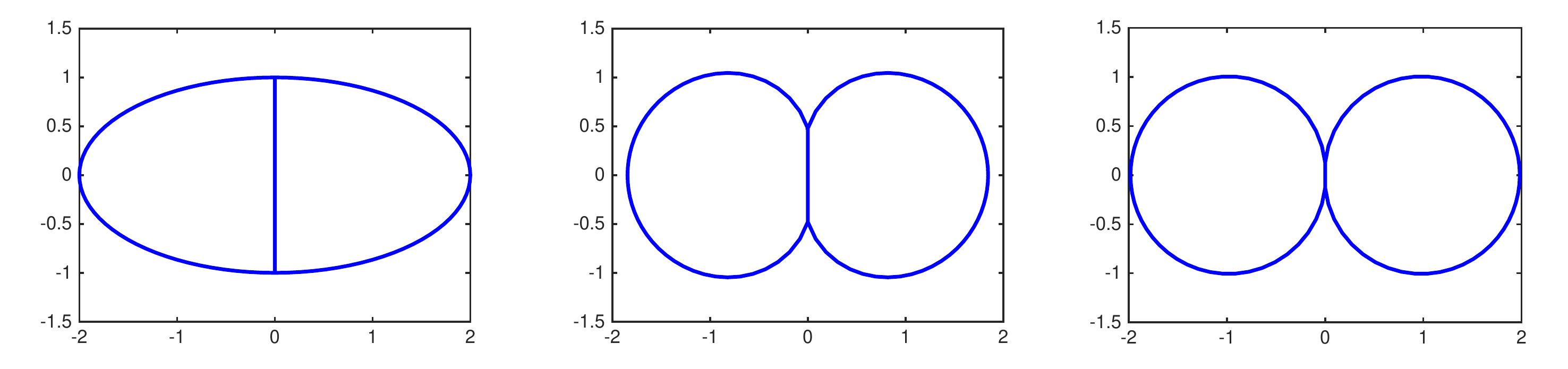}
\includegraphics[width=0.90\textwidth]{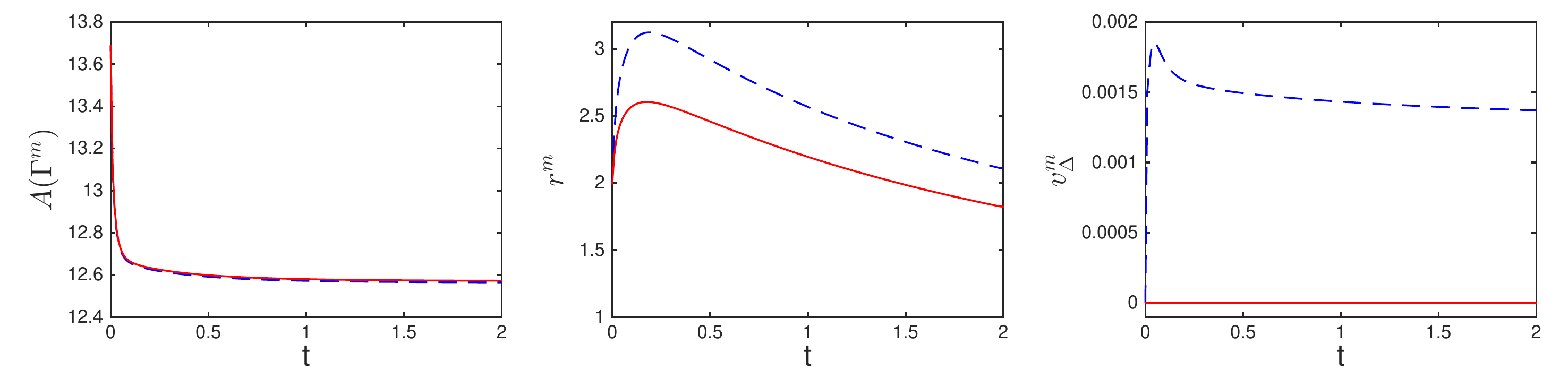}
\caption{
Evolution towards a degenerate 2d double bubble, with weightings $\sigma = (1,2,1)$.
Plots of $\Gamma^m$ at times $t=0, 0.1, 2$.
We also show plots of the discrete energy $A(\Gamma^m)$, the ratio $r^m$ and 
the relative volume error $v_\Delta^m$ over time, where $K = 129$ and $\ttau=10^{-2}$.
}
% ~/hpc_cluster/data/alberta/tjtrue/2d.db_1weight2
\label{fig:2ddb1w2}
\end{figure}

We then conduct experiments for the double bubble with different weightings of 
the surface energies, and the results are presented in 
Figs.~\ref{fig:2ddbw1.5}, \ref{fig:2ddbw2}, \ref{fig:2ddb1w1.5} 
and \ref{fig:2ddb1w2}. We observe that different weightings generally lead to different shapes of networks with different triple junction angles. For example, when $\sigma=(1,1,2)$, as time evolves, the triple junction angle between $\Gamma_1$ and $\Gamma_2$ approaches $0^\circ$ while the angles between $\Gamma_1$, $\Gamma_3$ and between $\Gamma_2$, $\Gamma_3$ tend to $180^\circ$, as shown in Fig.~\ref{fig:2ddbw2}. In fact, the third curve will finally shrink to a point, leading to a steady state of only two circular curves, as discussed in \cite{Barrett07}. Despite the different weightings being used, the energy dissipation and the volume conservation are satisfied, and the mesh quality is well preserved for the discrete numerical solutions in these experiments.

\begin{figure}[!htp]
\centering
\includegraphics[width=0.90\textwidth]{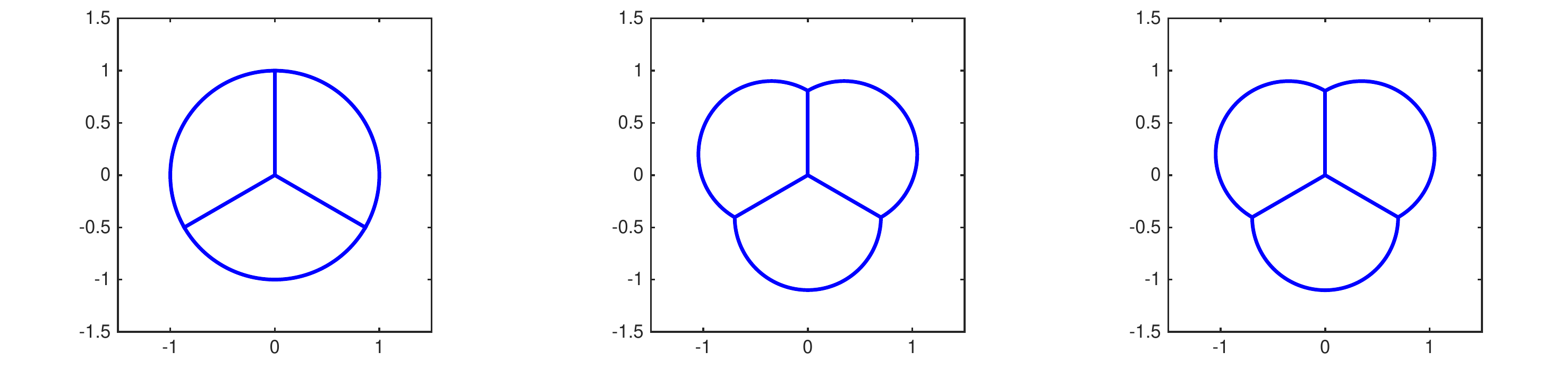}
\includegraphics[width=0.90\textwidth]{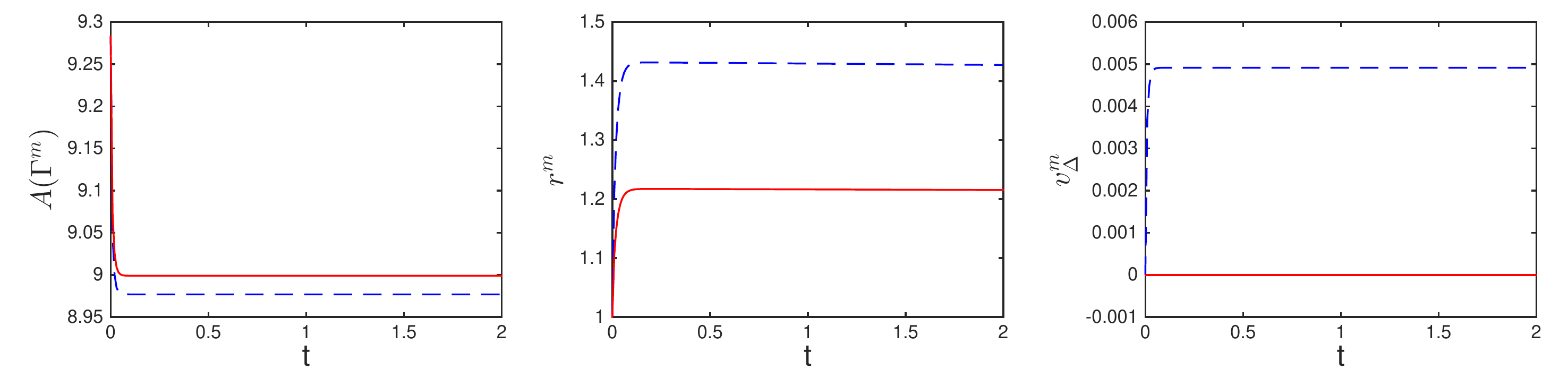}
\caption{
Evolution towards the 2d standard triple bubble.
Plots of $\Gamma^m$ at times $t=0, 0.1, 2$.
We also show plots of the discrete energy $A(\Gamma^m)$, the ratio $r^m$ and 
the relative volume error $v_\Delta^m$ over time, where $K = 1029$ and $\ttau=10^{-2}$.
}
% .apmc created in ~/c/triplej/triple_bubble/
% ~/hpc_cluster/data/alberta/tjtrue/2d.tb
\label{fig:2dtb}
\end{figure}

\begin{figure}[!htp]
\center
\includegraphics[width=0.90\textwidth]{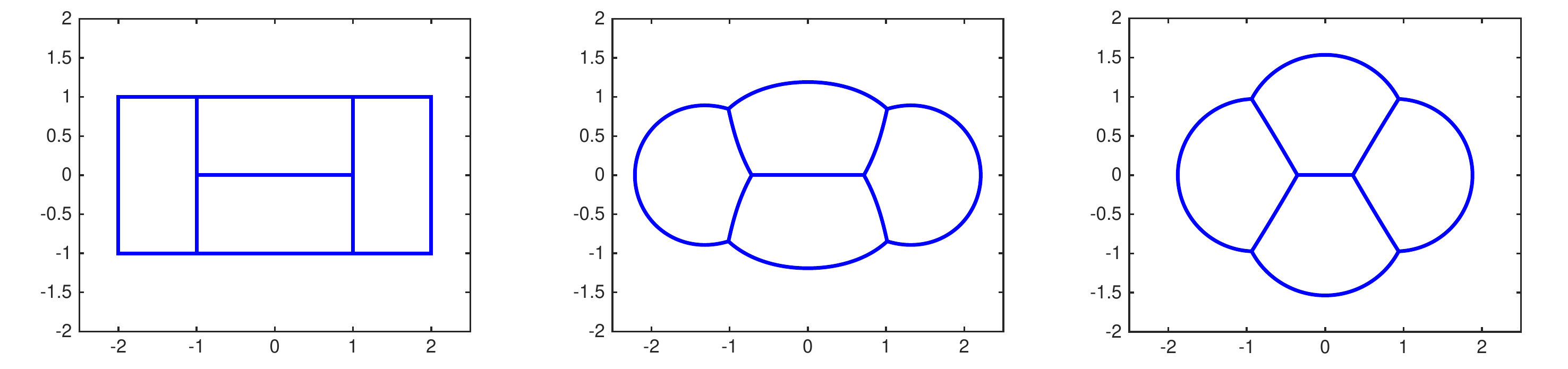}
\includegraphics[width=0.90\textwidth]{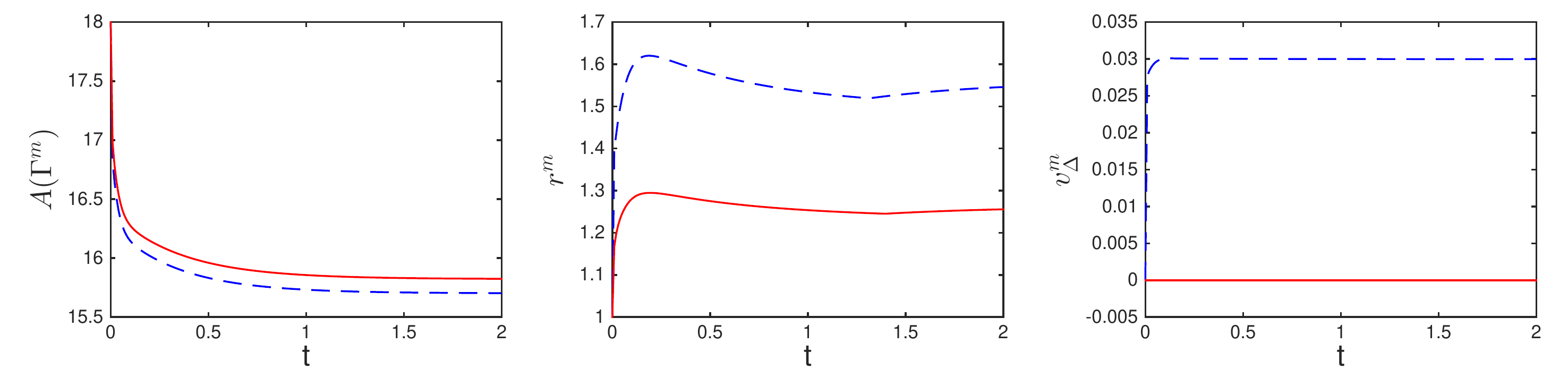}
\caption{
Evolution towards the 2d standard quadruple bubble.
Plots of $\Gamma^m$ at times $t=0, 0.1, 2$.
We also show plots of the discrete energy $A(\Gamma^m)$, the ratio $r^m$ and 
the relative volume error $v_\Delta^m$ over time, where $K = 1029$ and $\ttau=10^{-2}$.
}
% .apmc created in ~/c/triplej/quartic_bubble/
% ~/hpc_cluster/data/alberta/tjtrue/2d.qb
\label{fig:2dqb}
\end{figure}

\begin{figure}[!htp]
\center
\includegraphics[width=0.90\textwidth]{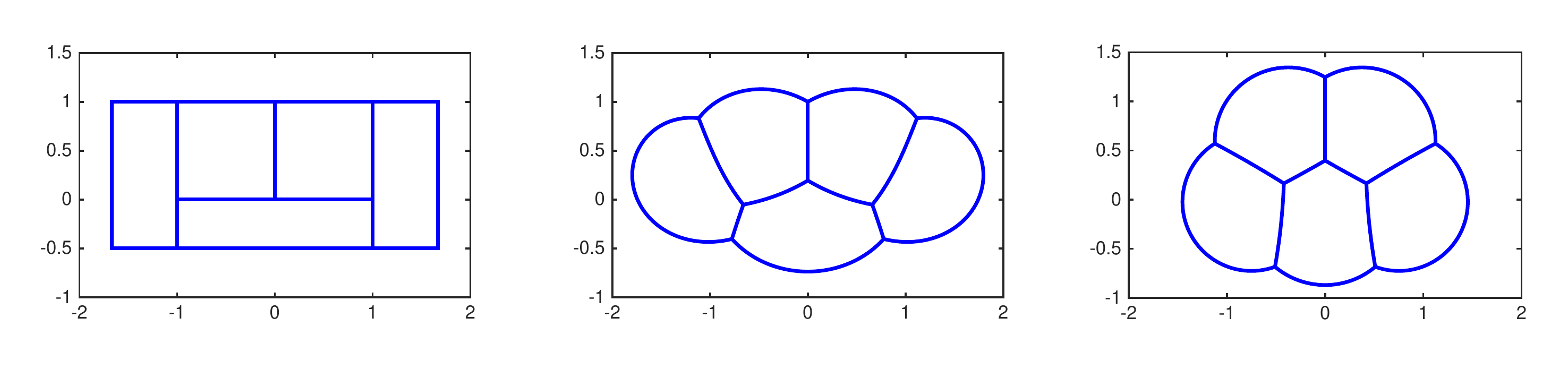}
\includegraphics[width=0.90\textwidth]{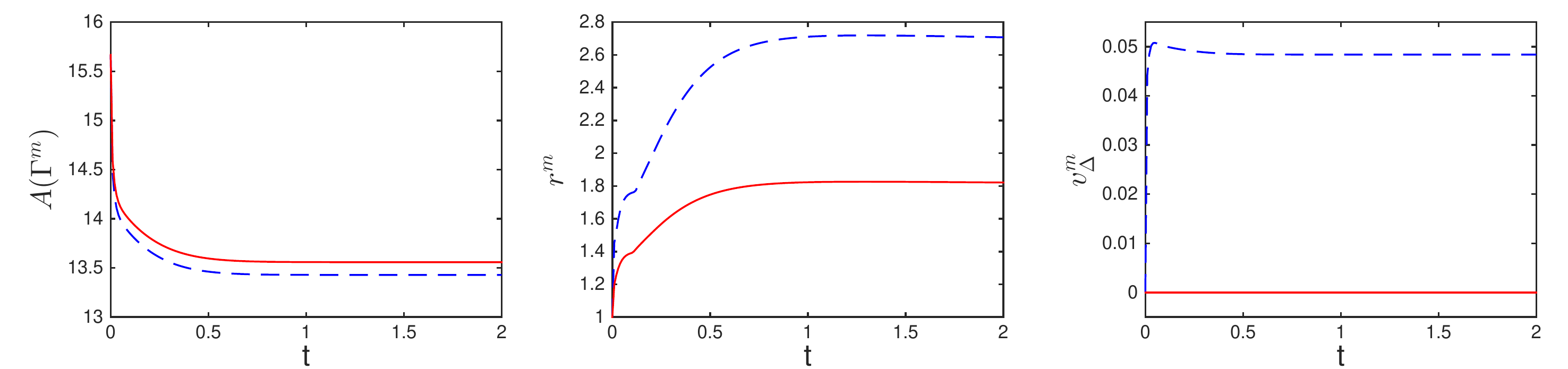}
\caption{
Evolution towards a possible 2d length minimizing quintuple bubble.
Plots of $\Gamma^m$ at times $t=0, 0.1, 2$.
We also show plots of the discrete energy $A(\Gamma^m)$, the ratio $r^m$ and 
the relative volume error $v_\Delta^m$ over time, where $K = 1032$ and $\ttau=10^{-2}$.
}
% .apmc created in ~/c/triplej/quintic_bubble/
% ~/hpc_cluster/data/alberta/tjtrue/2d.qnb
\label{fig:2dqnb}
\end{figure}

\begin{figure}[!htp]
\center
\includegraphics[width=0.90\textwidth]{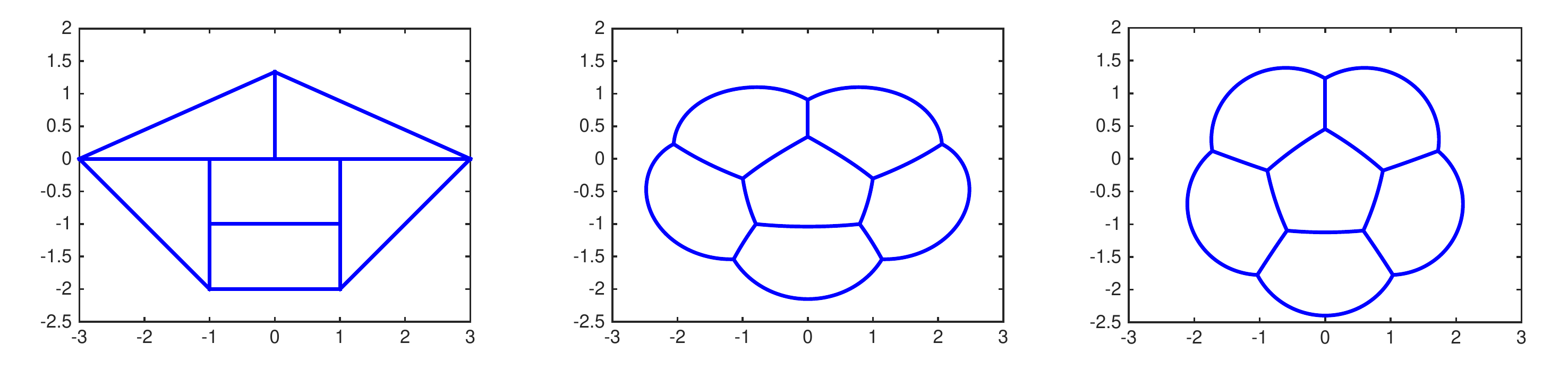}
\includegraphics[width=0.90\textwidth]{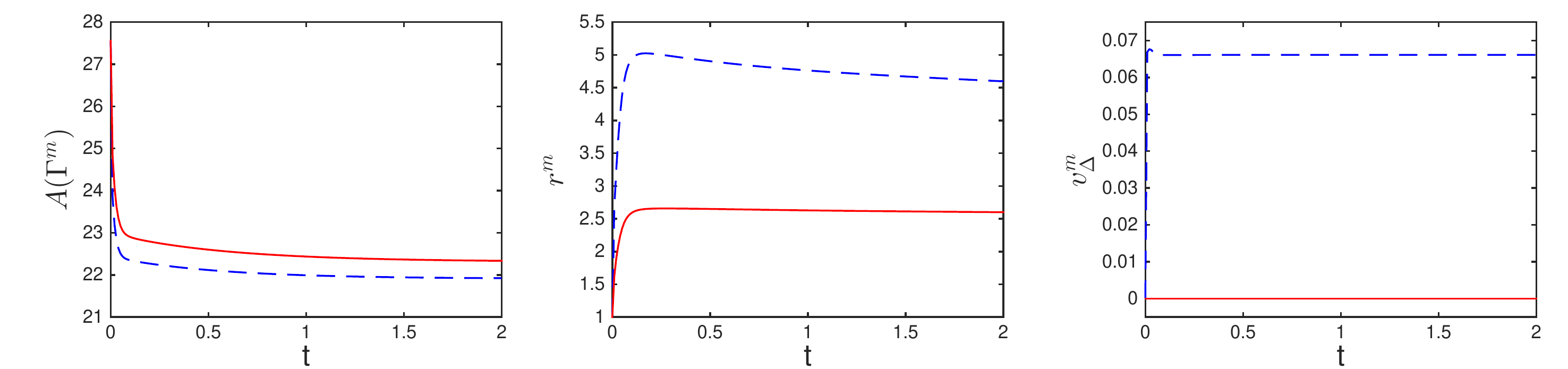}
\caption{
Evolution towards a possible 2d length minimizing sextuple bubble.
Plots of $\Gamma^m$ at times $t=0, 0.1, 2$.
We also show plots of the discrete energy $A(\Gamma^m)$, the ratio $r^m$ and 
the relative volume error $v_\Delta^m$ over time, where $K = 1025$ and $\ttau=10^{-2}$.
}
% .apmc created in ~/c/triplej/sextic_bubble/
% ~/hpc_cluster/data/alberta/tjtrue/2d.sb
\label{fig:2dsb}
\end{figure}

\begin{figure}[!htp]
\center
\includegraphics[width=0.90\textwidth]{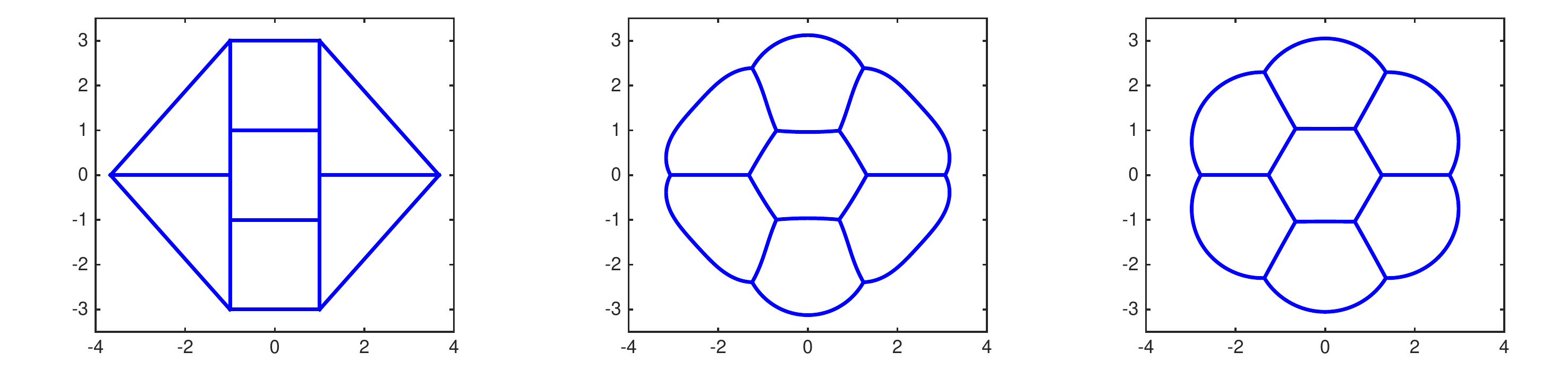}
\includegraphics[width=0.90\textwidth]{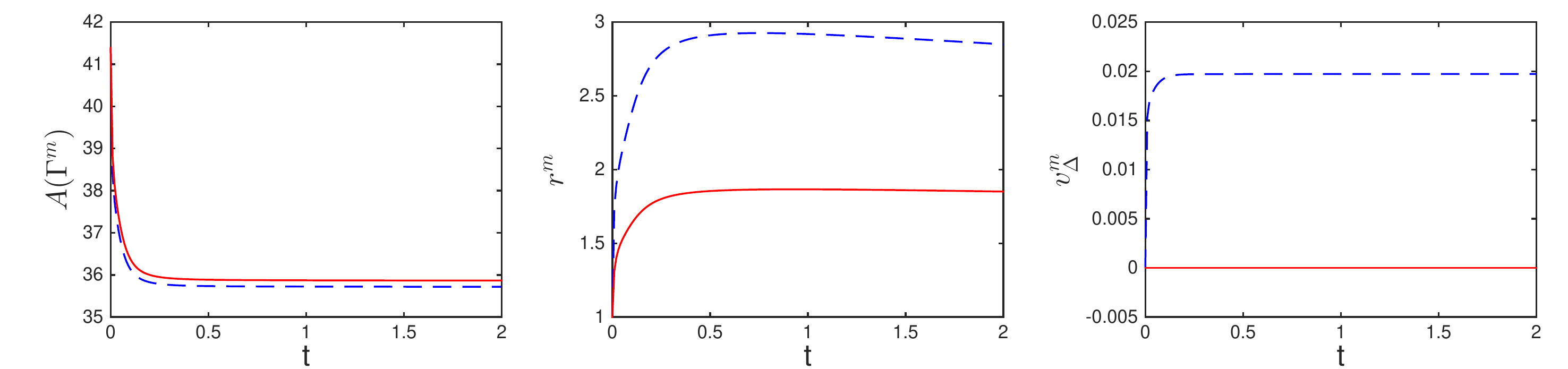}
\caption{
Evolution towards a possible 2d length minimizing septuple bubble.
Plots of $\Gamma^m$ at times $t=0, 0.1, 2$.
We also show plots of the discrete energy $A(\Gamma^m)$, the ratio $r^m$ and 
the relative volume error $v_\Delta^m$ over time, where $K = 1032$ and $\ttau=10^{-2}$.
}
% .apmc created in ~/c/triplej/septic_bubble/
% ~/hpc_cluster/data/alberta/tjtrue/2d.spb
\label{fig:2dspb}
\end{figure}

We next perform simulations for the standard triple, quadruple, quintuple, sextuple and 
septuple bubbles with equal surface energy densities, as shown in Figs.~\ref{fig:2dtb}, \ref{fig:2dqb}, \ref{fig:2dqnb},
\ref{fig:2dsb} and \ref{fig:2dspb}, respectively. We observe the energy is decreasing and the mesh ratio remains at small values for the numerical solutions during the simulation. In particular, in all these simulations the volume of the enclosed bubbles is preserved exactly for the introduced SP-PFEM. However, for the BGN scheme the observed relative volume loss can be up to $6.5\%$ during the evolution, as can be seen from the last subfigure in Fig. \ref{fig:2dsb}. These results demonstrate the reliability of our method.

\subsection{Anisotropic numerical results in 2d}

\begin{figure}[!htp]
\center
\includegraphics[width=0.9\textwidth]{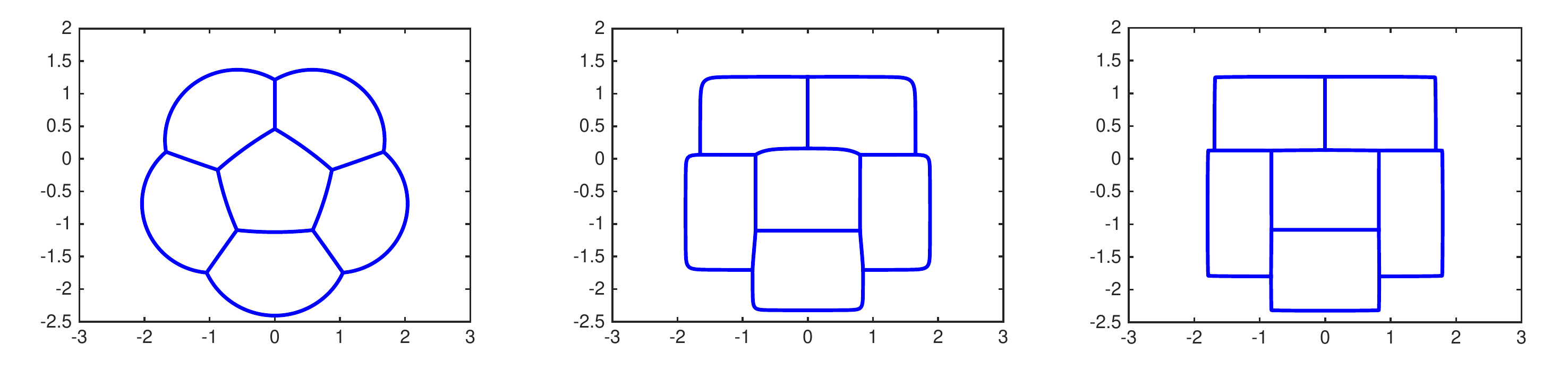}
\includegraphics[width=0.75\textwidth]{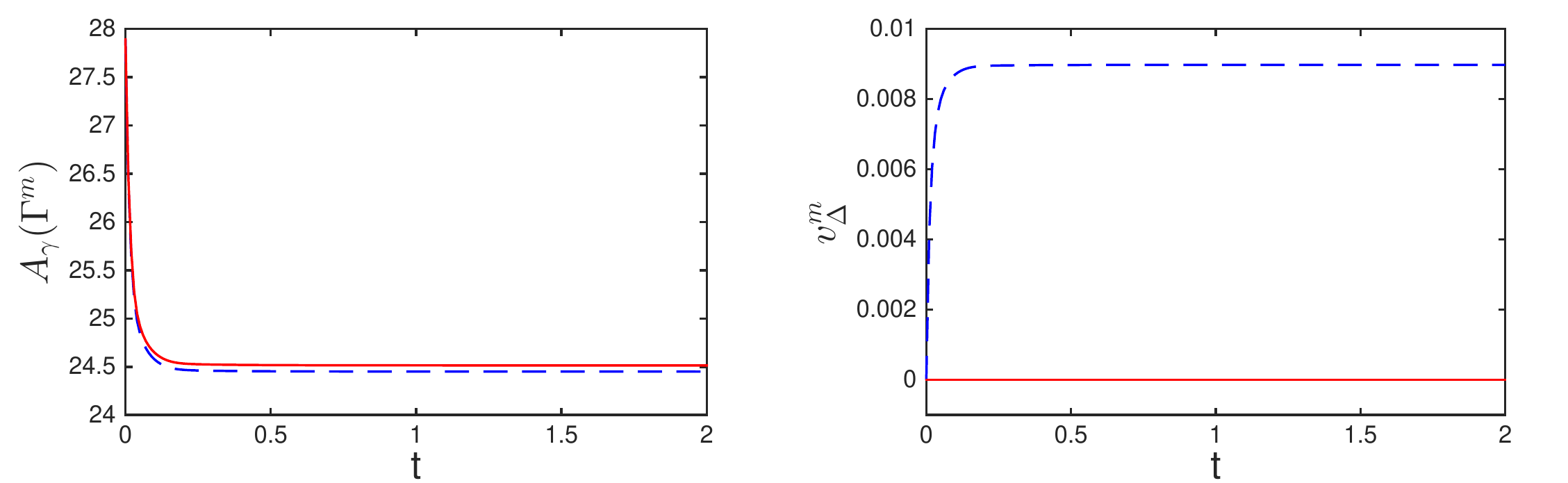}
\caption{
Evolution towards an anisotropic 2d sextuple bubble, for the 
anisotropy \eqref{eq:2dani} with $L=2$ and $\epsilon=0.01$.
Plots of $\Gamma^m$ at times $t=0, 0.1, 2$.
We also show plots of the discrete energy $A_\gamma(\Gamma^m)$ and the relative volume error $v_\Delta^m$ over time, where $K = 1025$ and $\ttau=10^{-2}$.
}
\label{fig:2dsbL2}
\end{figure}

\begin{figure}[!htp]
\center
\includegraphics[width=0.90\textwidth]{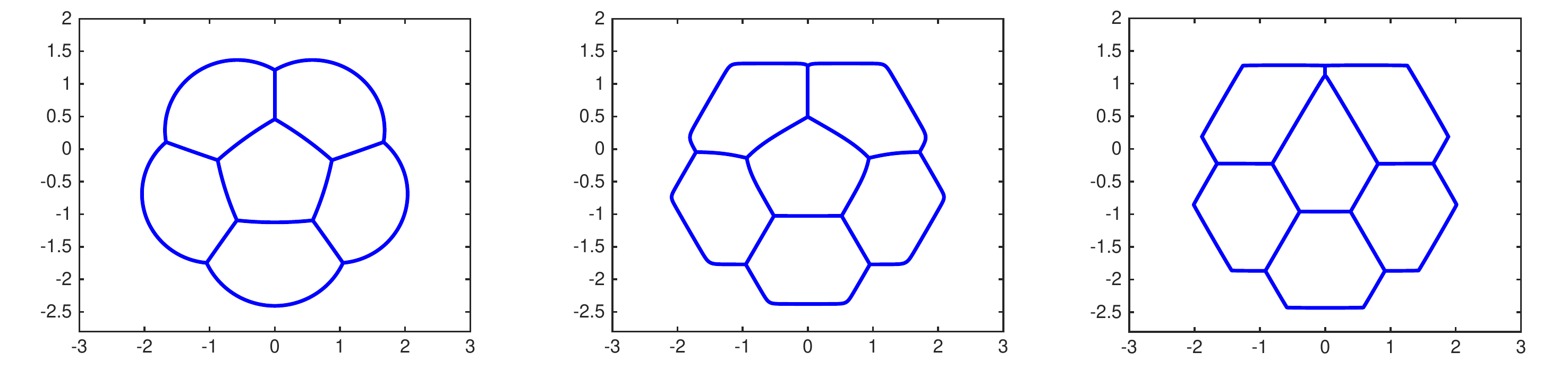}
\includegraphics[width=0.75\textwidth]{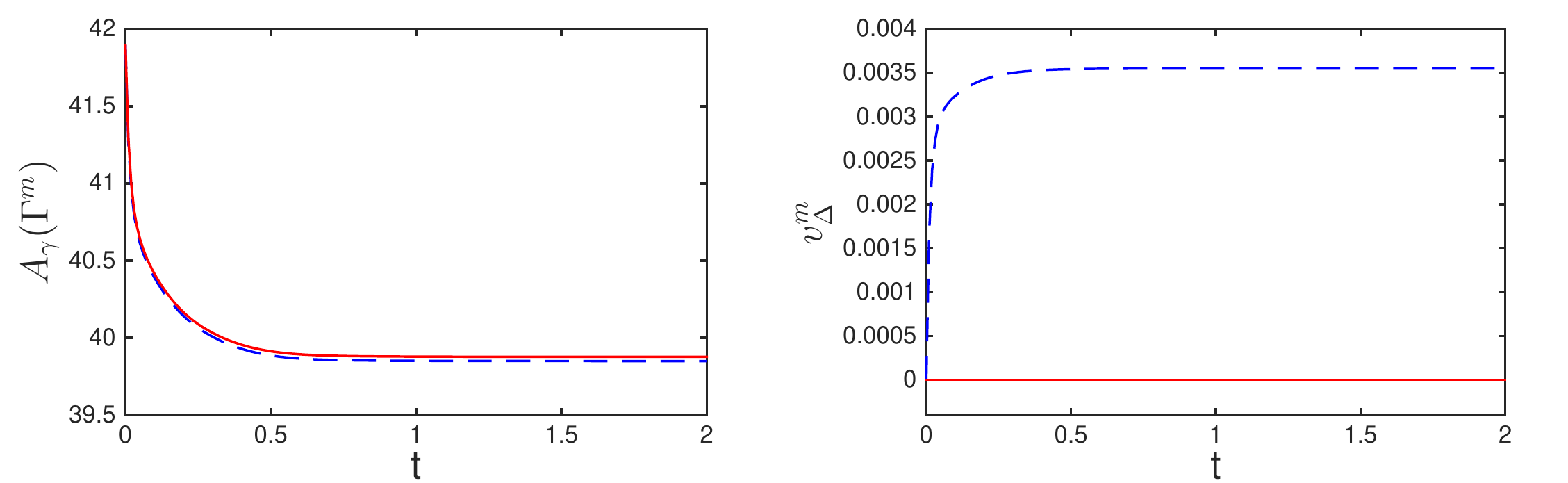}
\caption{
Evolution towards an anisotropic 2d sextuple bubble, for the 
anisotropy \eqref{eq:2dani} with $L=3$ and $\epsilon=0.01$.
Plots of $\Gamma^m$ at times $t=0, 0.1, 2$.
We also show plots of the discrete energy $A_\gamma(\Gamma^m)$ and the relative volume error $v_\Delta^m$ over time, where $K = 1025$ and $\ttau=10^{-2}$. 
}
\label{fig:2dsbL3}
\end{figure}

\begin{figure}[!htp]
\center
\includegraphics[width=0.90\textwidth]{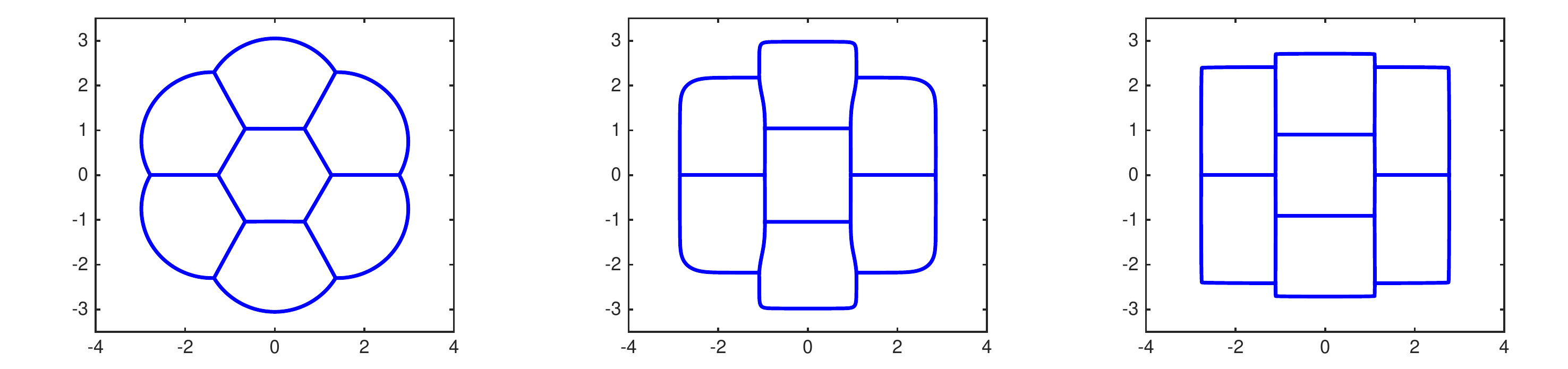}
\includegraphics[width=0.75\textwidth]{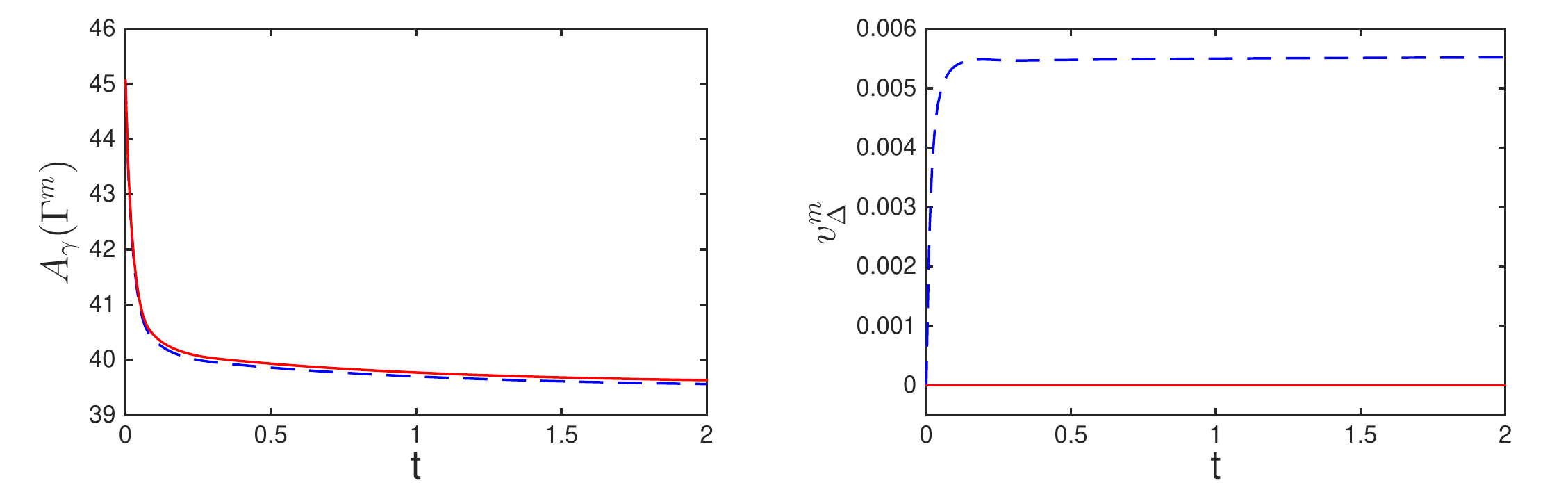}
\caption{
Evolution towards an anisotropic 2d septuple bubble, for the 
anisotropy \eqref{eq:2dani} with $L=2$ and $\epsilon=0.01$.
Plots of $\Gamma^m$ at times $t=0, 0.1, 2$.
We also show plots of the discrete energy $A_\gamma(\Gamma^m)$ and the relative volume error $v_\Delta^m$ over time, where $K = 1032$ and $\ttau=10^{-2}$. 
}
\label{fig:2dspbL2}
\end{figure}

\begin{figure}[!htp]
\center
\includegraphics[width=0.90\textwidth]{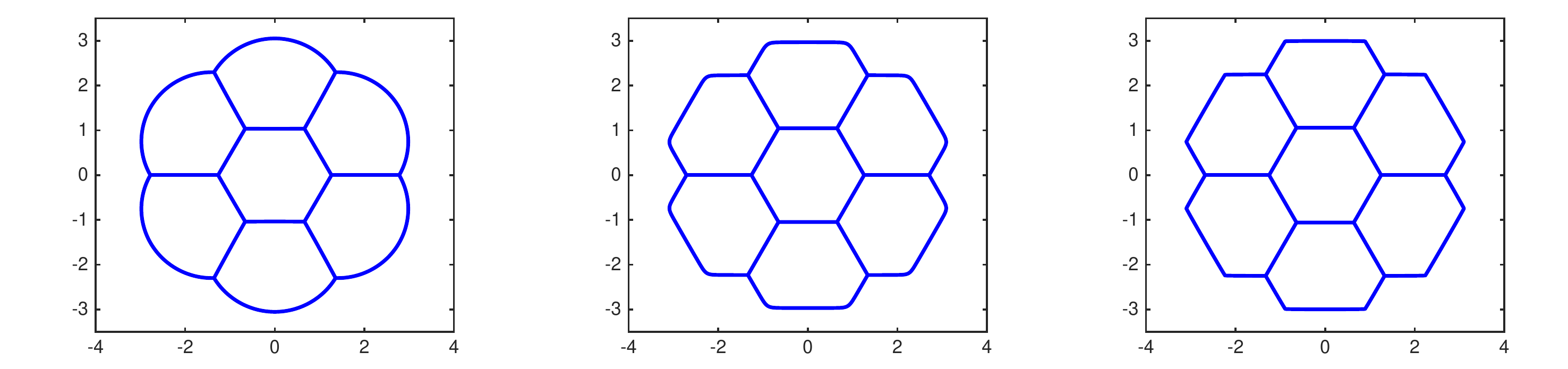}
\includegraphics[width=0.75\textwidth]{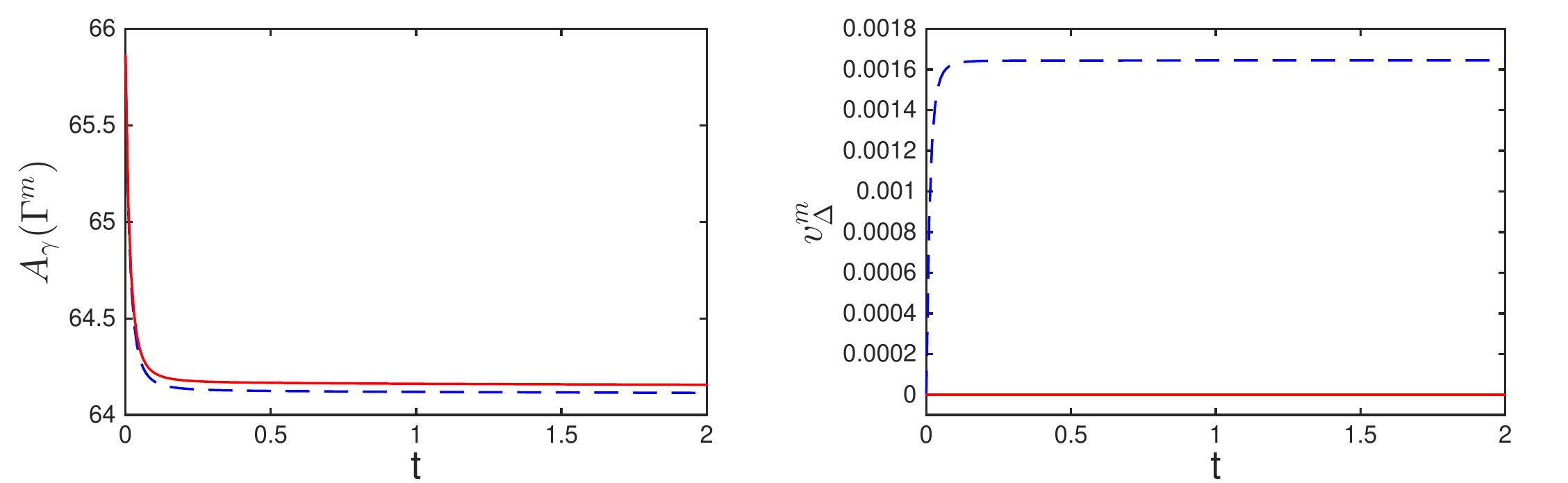}
\caption{
Evolution towards an anisotropic 2d septuple bubble, for the 
anisotropy \eqref{eq:2dani} with $L=3$ and $\epsilon=0.01$.
Plots of $\Gamma^m$ at times $t=0, 0.1, 2$.
We also show plots of the discrete energy $A_\gamma(\Gamma^m)$ and the relative volume error $v_\Delta^m$ over time, where $K = 1032$ and $\ttau=10^{-2}$.
}
\label{fig:2dspbL3}
\end{figure}

We simulate the evolution of curve networks with the anisotropy given by 
\begin{equation}
\gamma(\vec p) = \sum_{\ell = 1}^L \sqrt{\vec p \cdot R(-\tfrac{(\ell-1)\pi}{L})D(\epsilon) R(\tfrac{(\ell-1)\pi}{L})\vec p}\quad\mbox{with}\quad R(\theta) = \left(\begin{matrix}
\cos\theta & \sin\theta\\
-\sin\theta & \cos\theta
\end{matrix}\right), 
\label{eq:2dani} 
\end{equation}
where $D(\epsilon) = {\rm diag}(1,\epsilon^2)$ and $R(\theta)$ is a clockwise rotation matrix through the given angle $\theta$. 
Note that for $L=2$ the anisotropy \eqref{eq:2dani} is the same as
\eqref{eq:cuspgamma} for $d=2$ and $r=1$.
In the first simulation, we repeat the experiment from Fig.~\ref{fig:2dsb} for the anisotropy \eqref{eq:2dani} with $L=2$ and
$\epsilon=0.01$. The results are shown in Fig.~\ref{fig:2dsbL2}. 
Similarly, we show in Fig.~\ref{fig:2dsbL3} the corresponding evolution for the anisotropy \eqref{eq:2dani} with $L=3$ and $\epsilon=0.01$.  
In both cases it can be observed that the circular segments of the cluster in
the isotropic case now become facetted, with the orientations of the facets
aligned with the Wulff shape of the anisotropy.
We also repeat the experiment from Fig.~\ref{fig:2dspb} with the two considered anisotropies, and the numerical results are presented in Fig.~\ref{fig:2dspbL2} and
Fig.~\ref{fig:2dspbL3}, respectively. 
Once again, the previously smooth parts of the steady state clusters now become
facetted.
It is clearly observed that in all of these experiments the volume conservation and energy dissipation are well satisfied for the numerical solutions.
%, and different steady states have been obtained due to the anisotropies.  

\subsection{Numerical results in 3d} 

\begin{figure}[!htb]
\center
\includegraphics[angle=-0,width=0.3\textwidth]{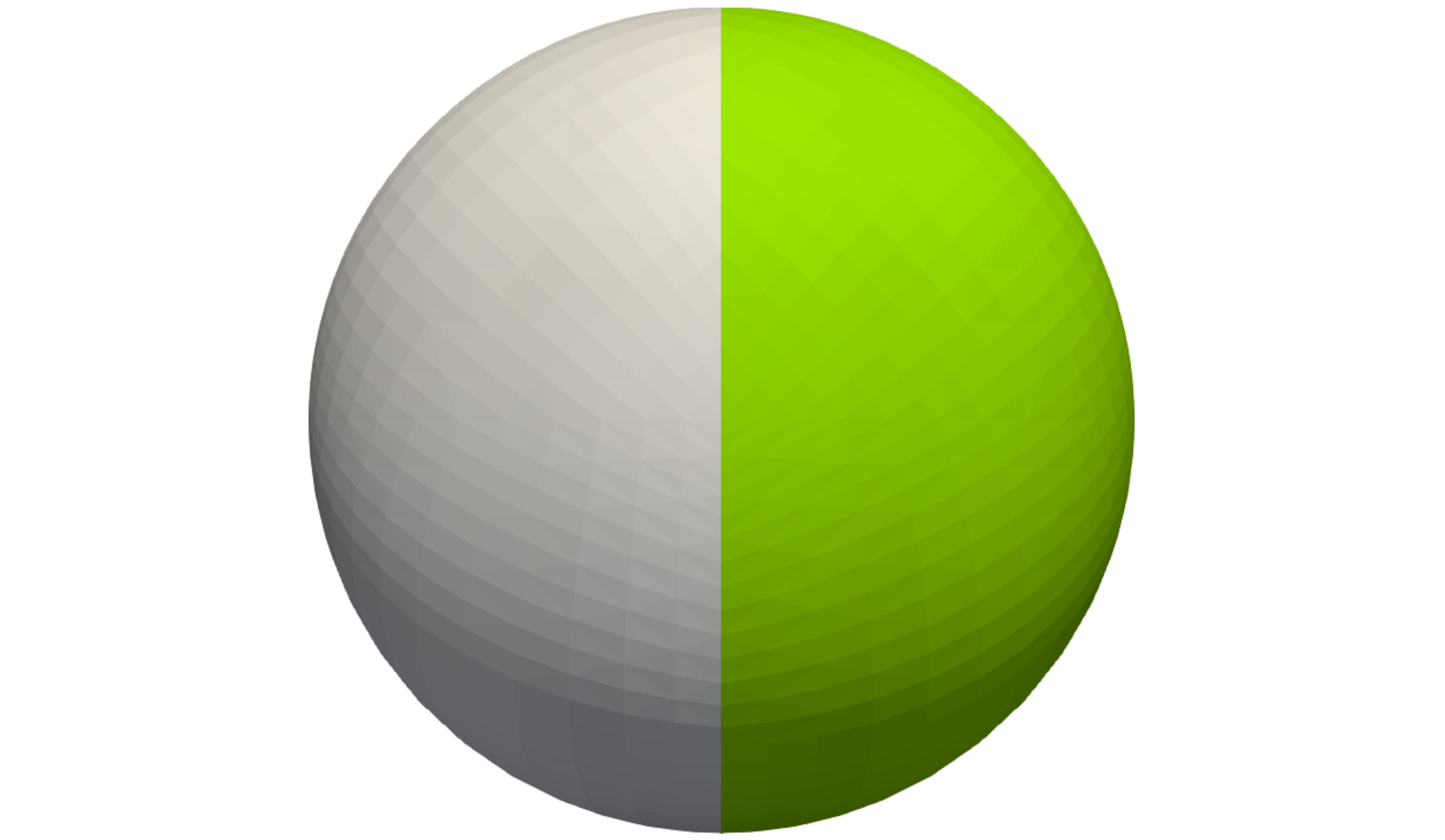}
\includegraphics[angle=-0,width=0.3\textwidth]{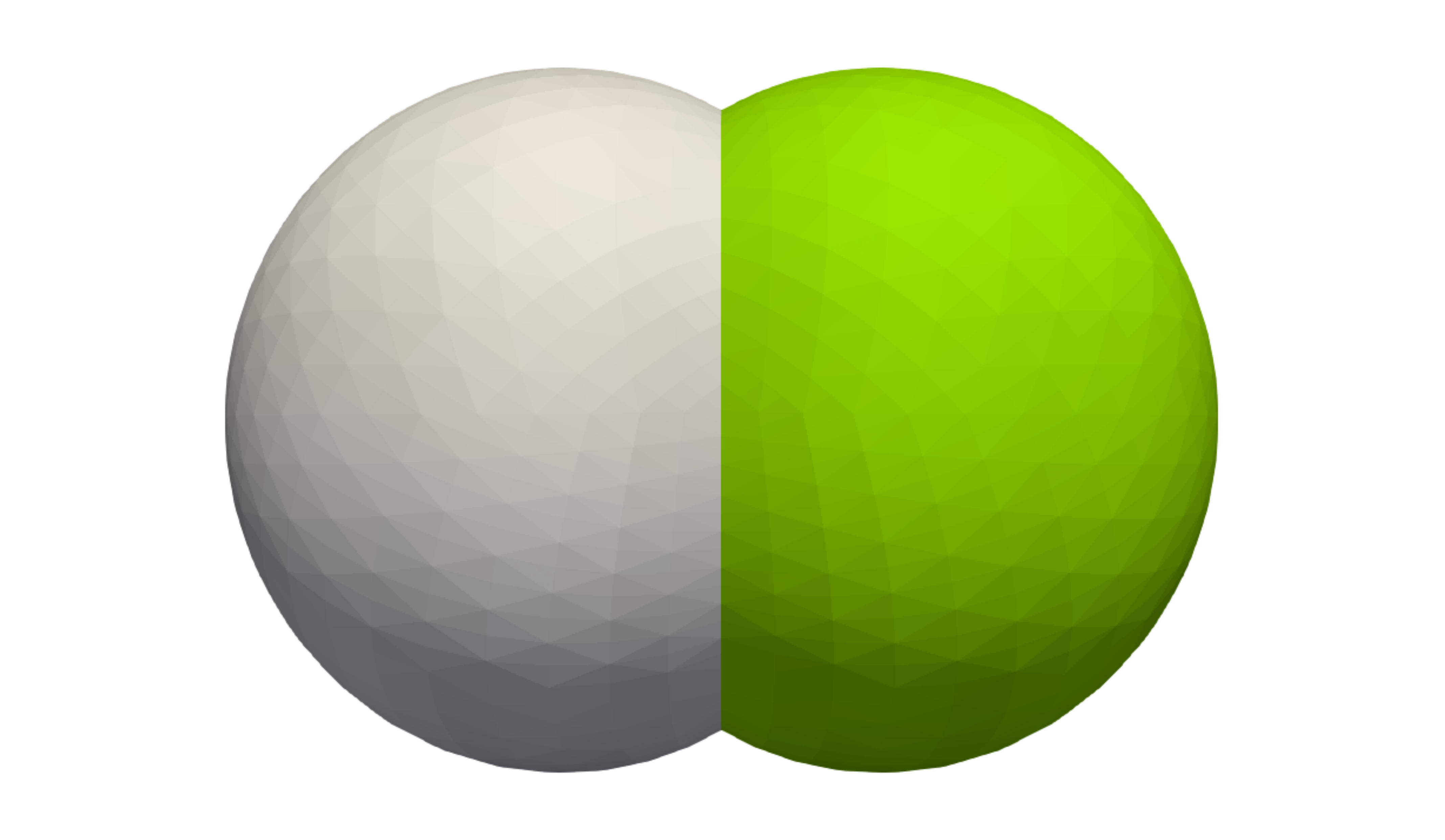}
\includegraphics[angle=-0,width=0.3\textwidth]{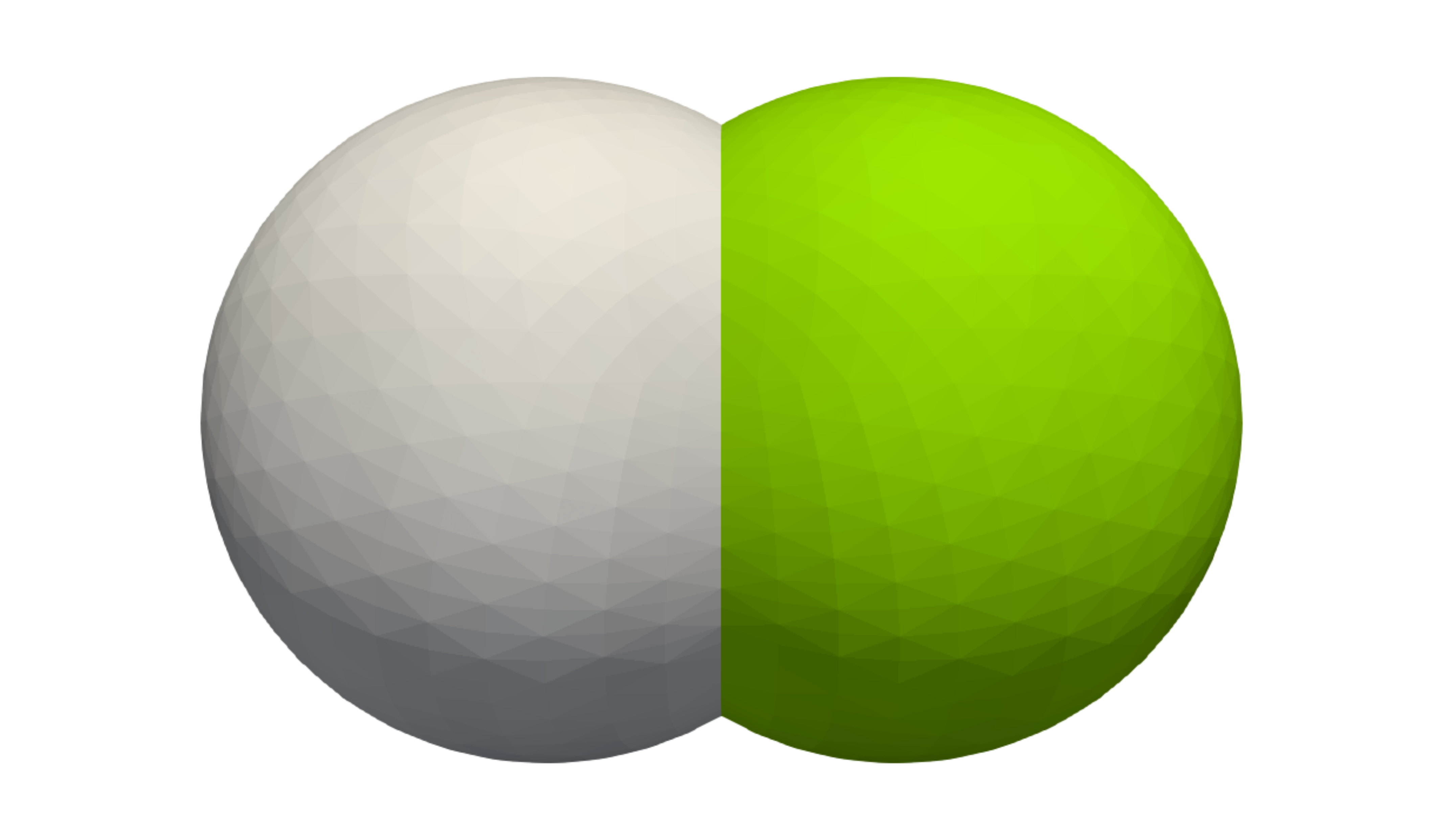}\\[0.5em]
\includegraphics[width=0.7\textwidth]{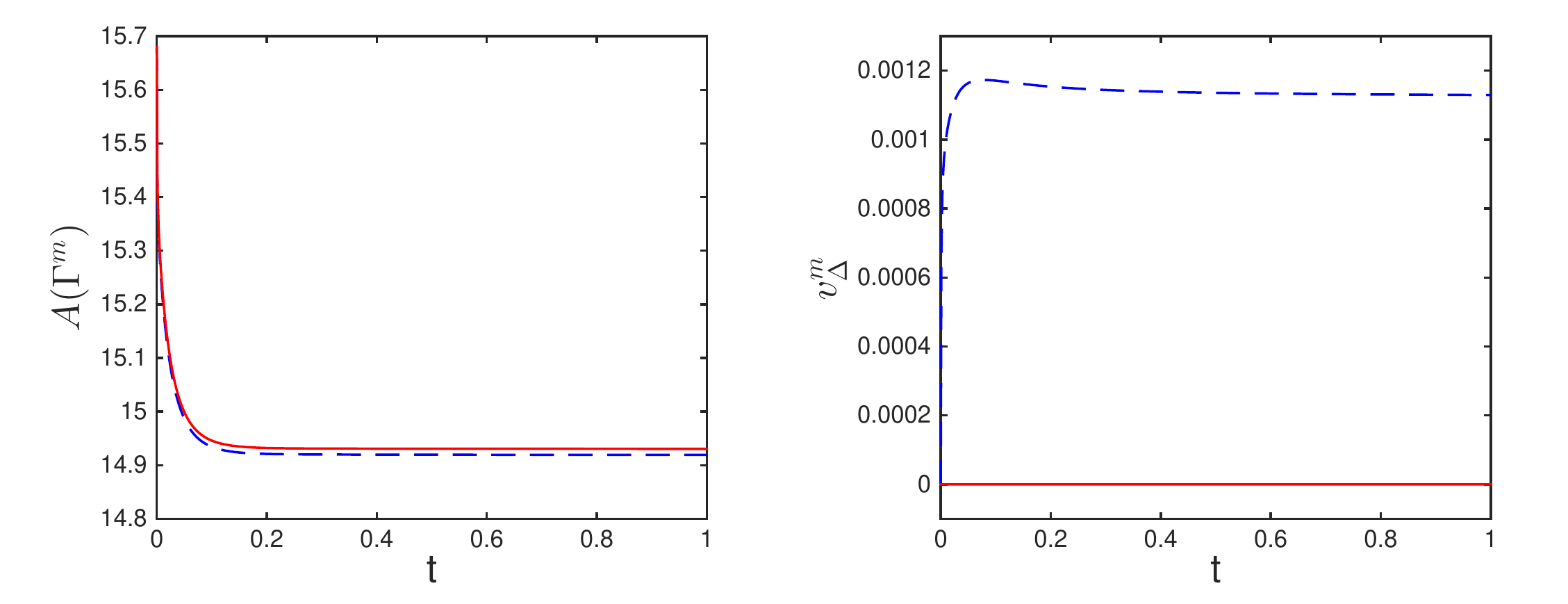}
\caption{
Evolution towards the 3d standard double bubble.
Plots of $\Gamma^m$ at times $t=0, 0.1, 1$.
We also show plots of the discrete energy $A(\Gamma^m)$ and 
the relative volume error $v_\Delta^m$ over time, where $K = 3267$ and $\ttau=10^{-3}$.
}
% ~/hpc_cluster/data/alberta/tjtrue/3d.db
\label{fig:3ddbw1}
\end{figure}%

\begin{figure}[!htb]
\center
\includegraphics[angle=-0,width=0.3\textwidth]{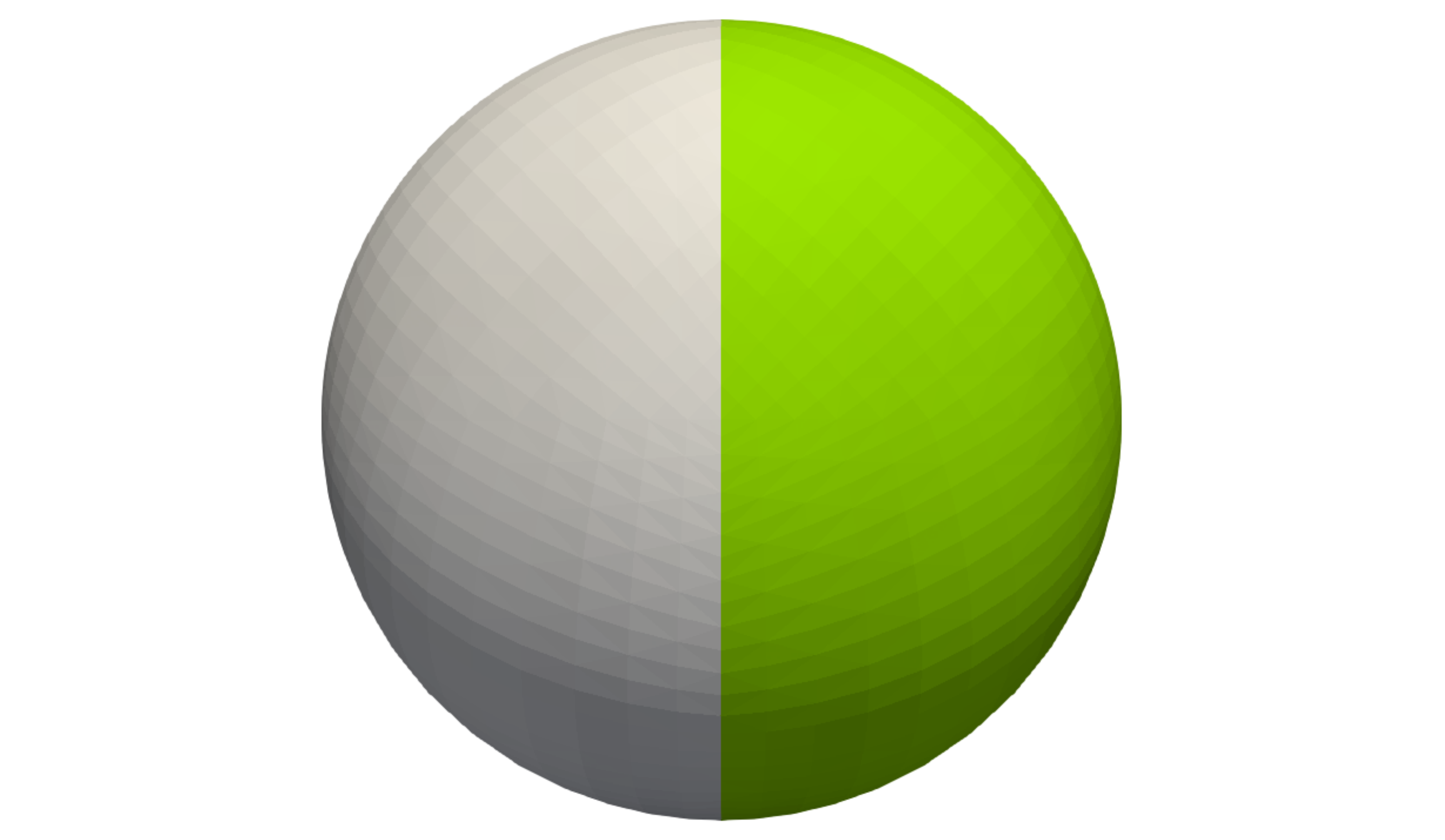}
\includegraphics[angle=-0,width=0.3\textwidth]{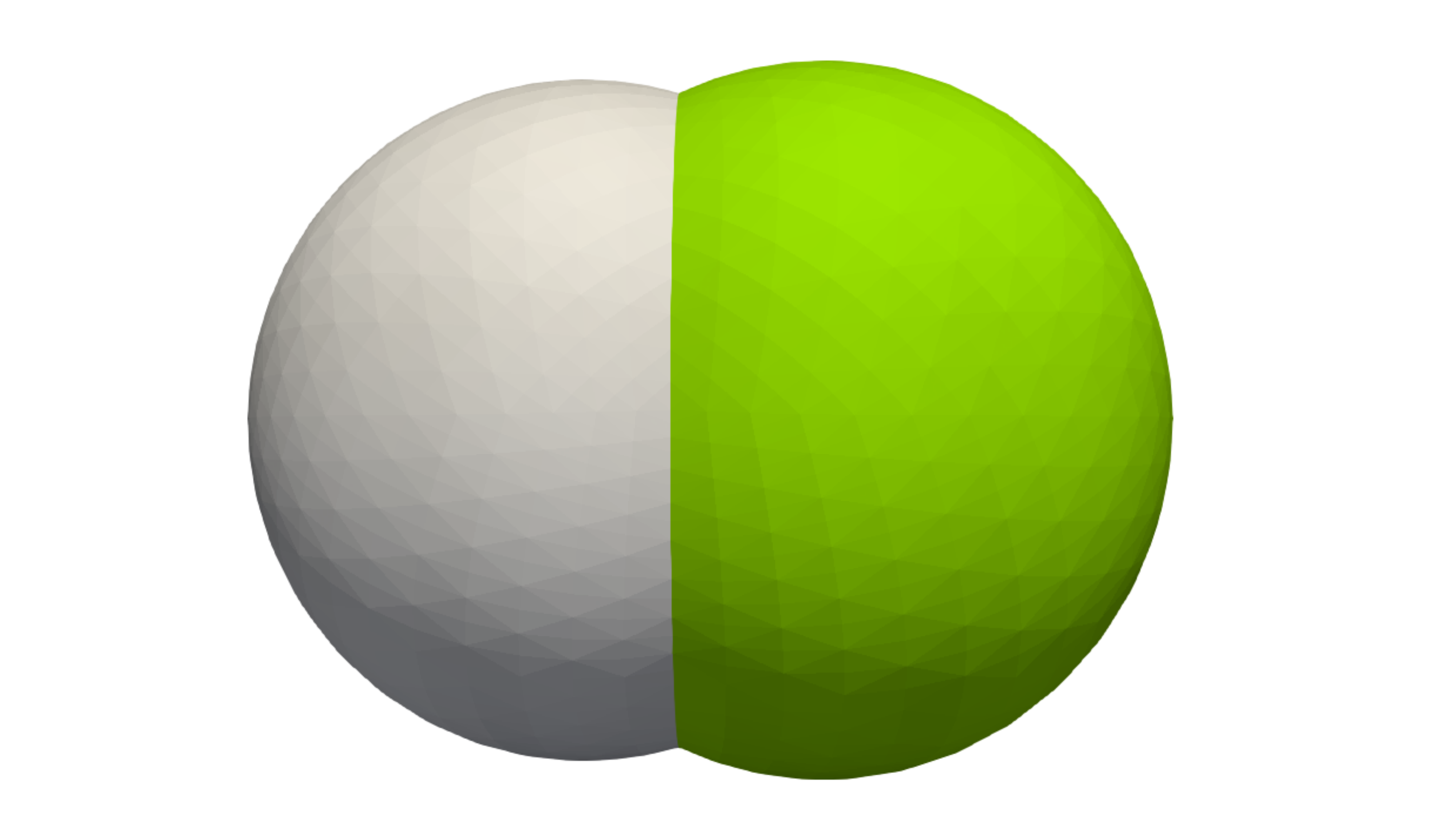}
\includegraphics[angle=-0,width=0.3\textwidth]{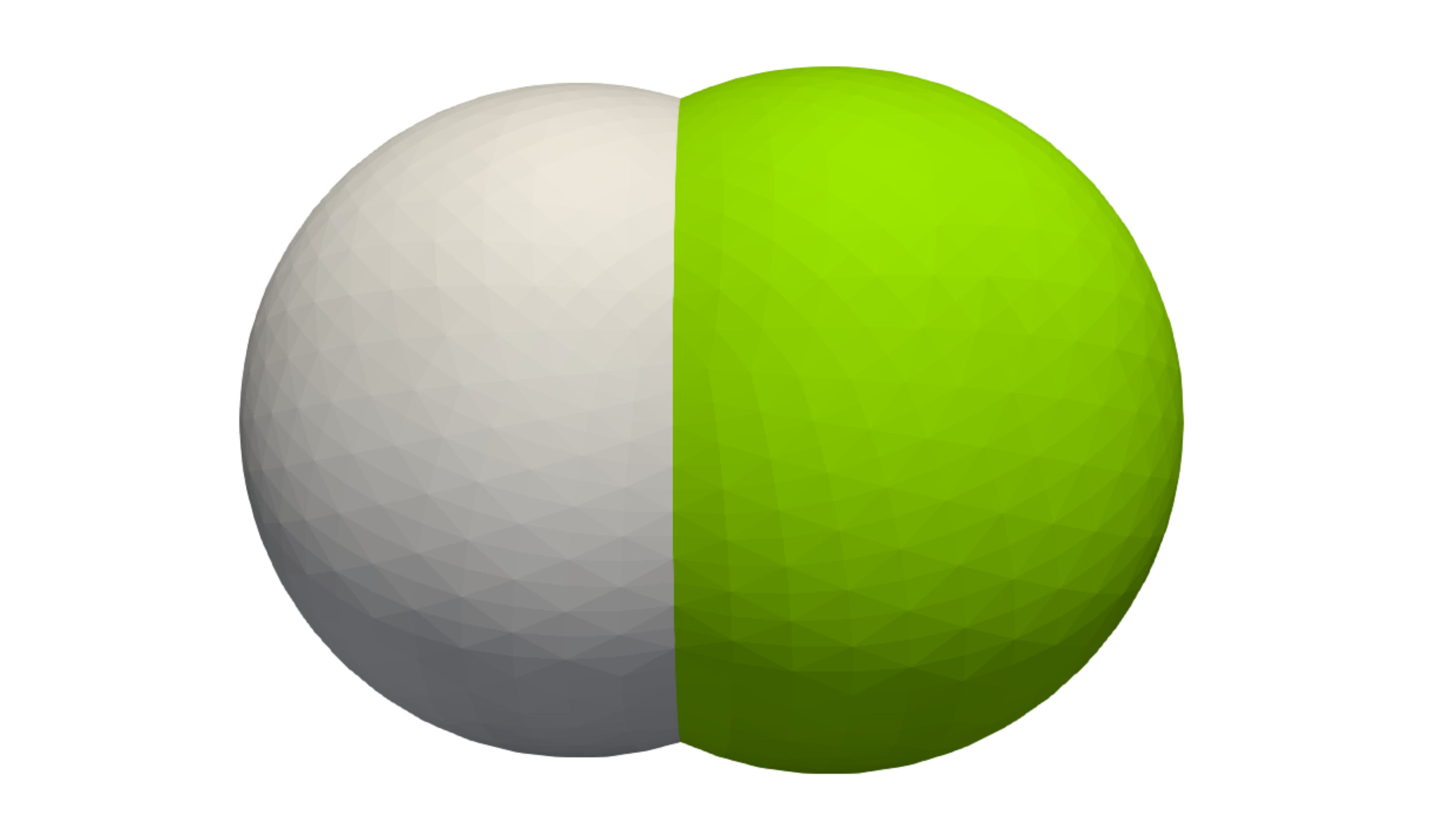}\\[0.5em]
\includegraphics[width=0.7\textwidth]{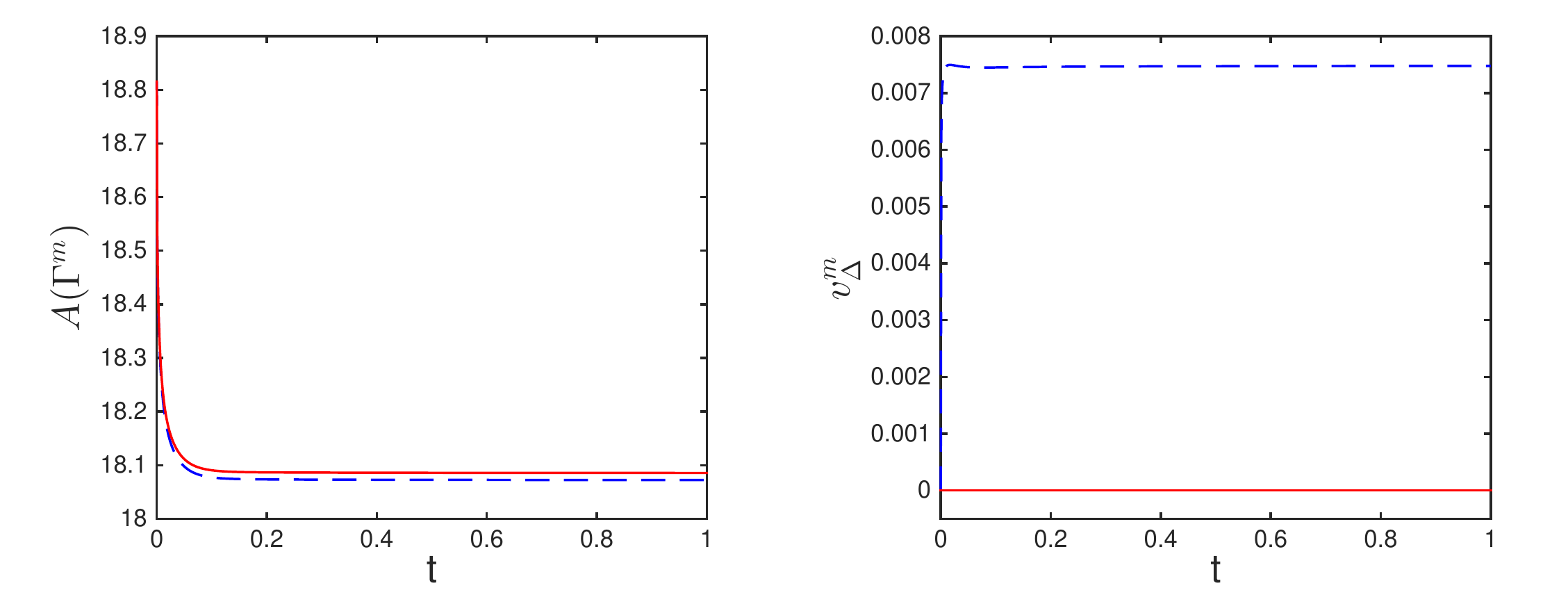}
\caption{
Evolution towards a 3d double bubble, with weightings $\sigma = (1.5,1,1)$.
Plots of $\Gamma^m$ at times $t=0, 0.1, 1$.
We also show plots of the discrete energy $A(\Gamma^m)$ and 
the relative volume error $v_\Delta^m$ over time, where $K = 3267$ and $\ttau=10^{-3}$.
}
% ~/hpc_cluster/data/alberta/tjtrue/3d.db_weight1.5
\label{fig:3ddbw1.5}
\end{figure}

\begin{figure}[!htb]
\center
\includegraphics[angle=-0,width=0.3\textwidth]{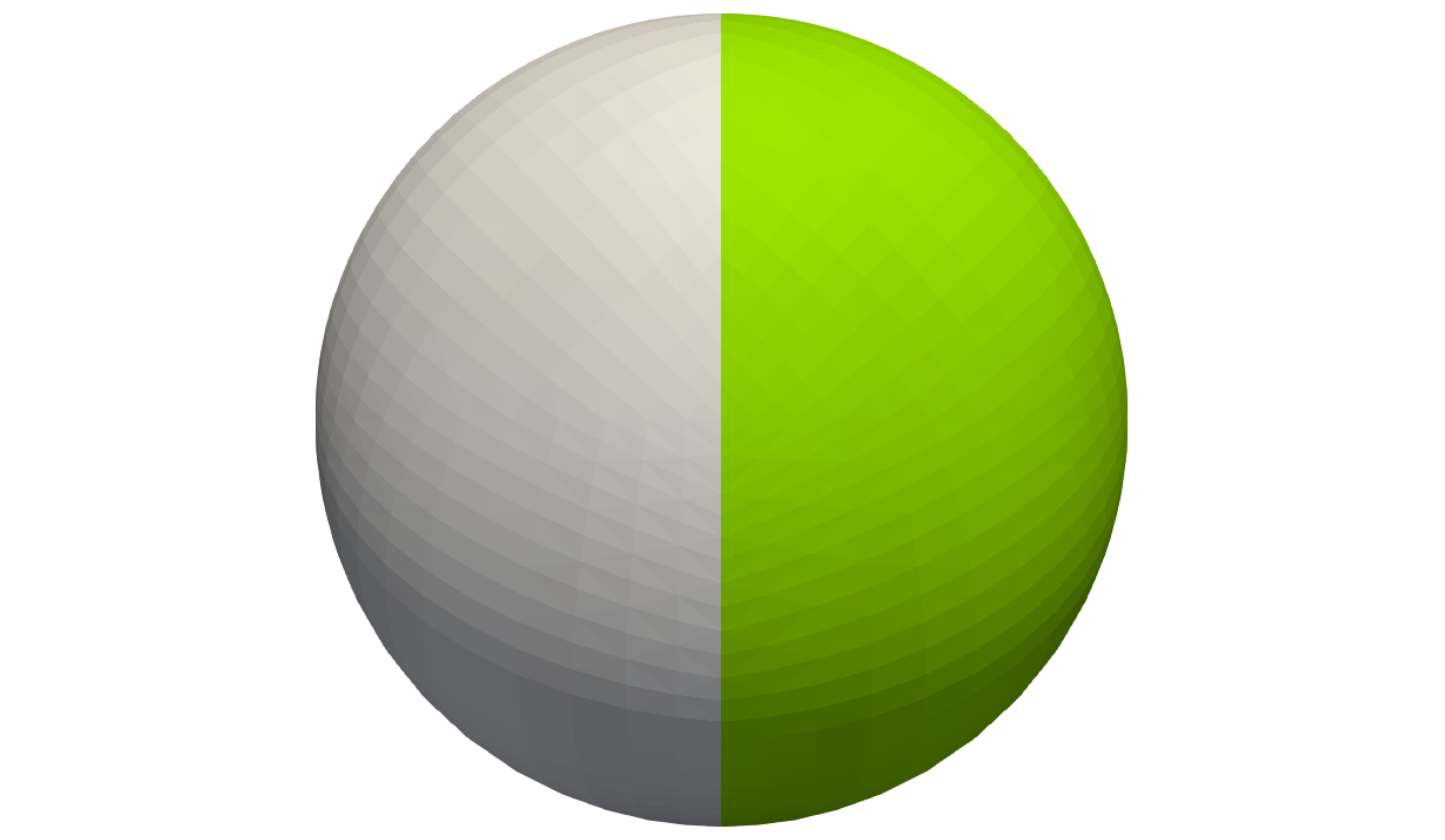}
\includegraphics[angle=-0,width=0.3\textwidth]{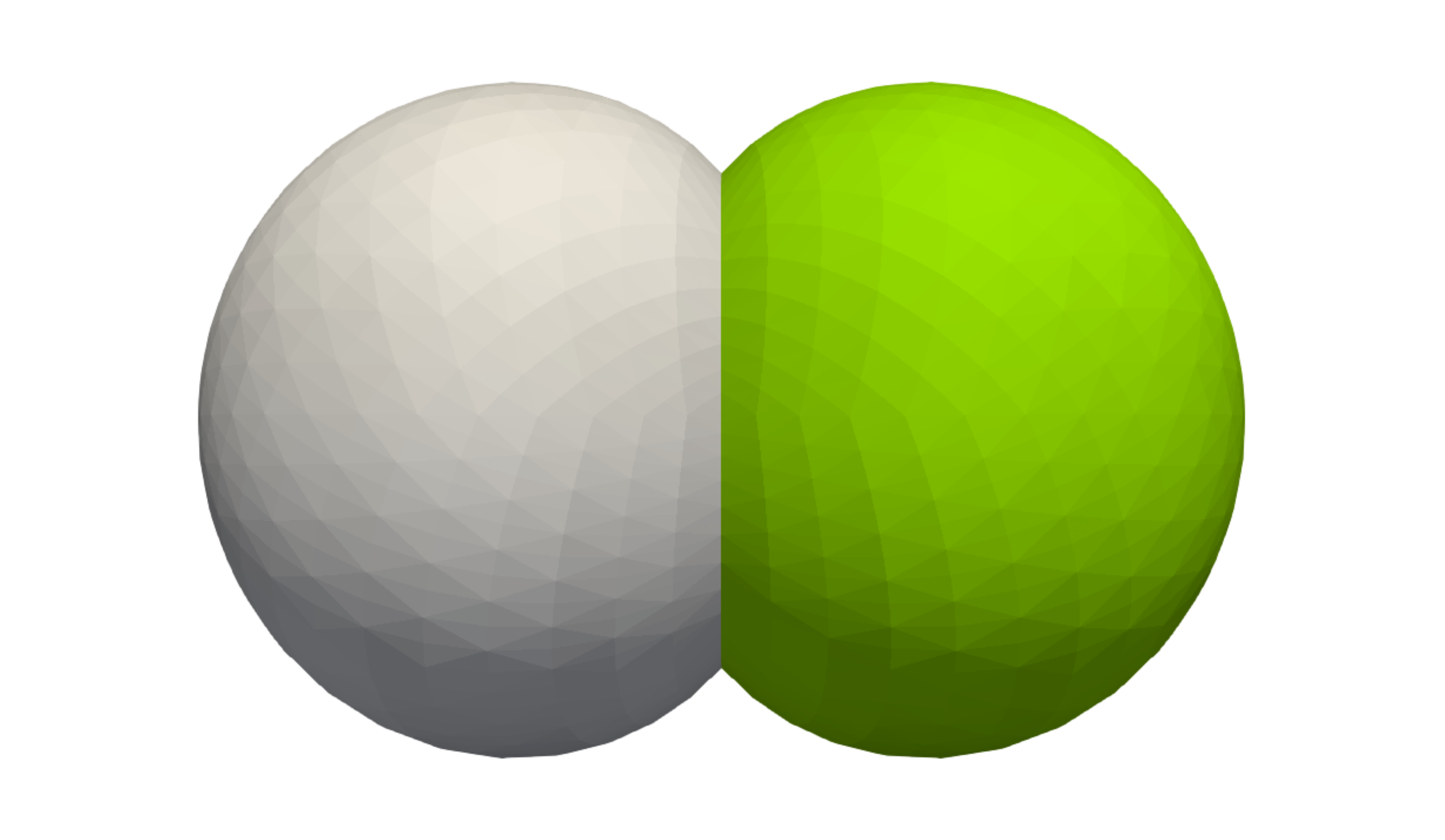}
\includegraphics[angle=-0,width=0.3\textwidth]{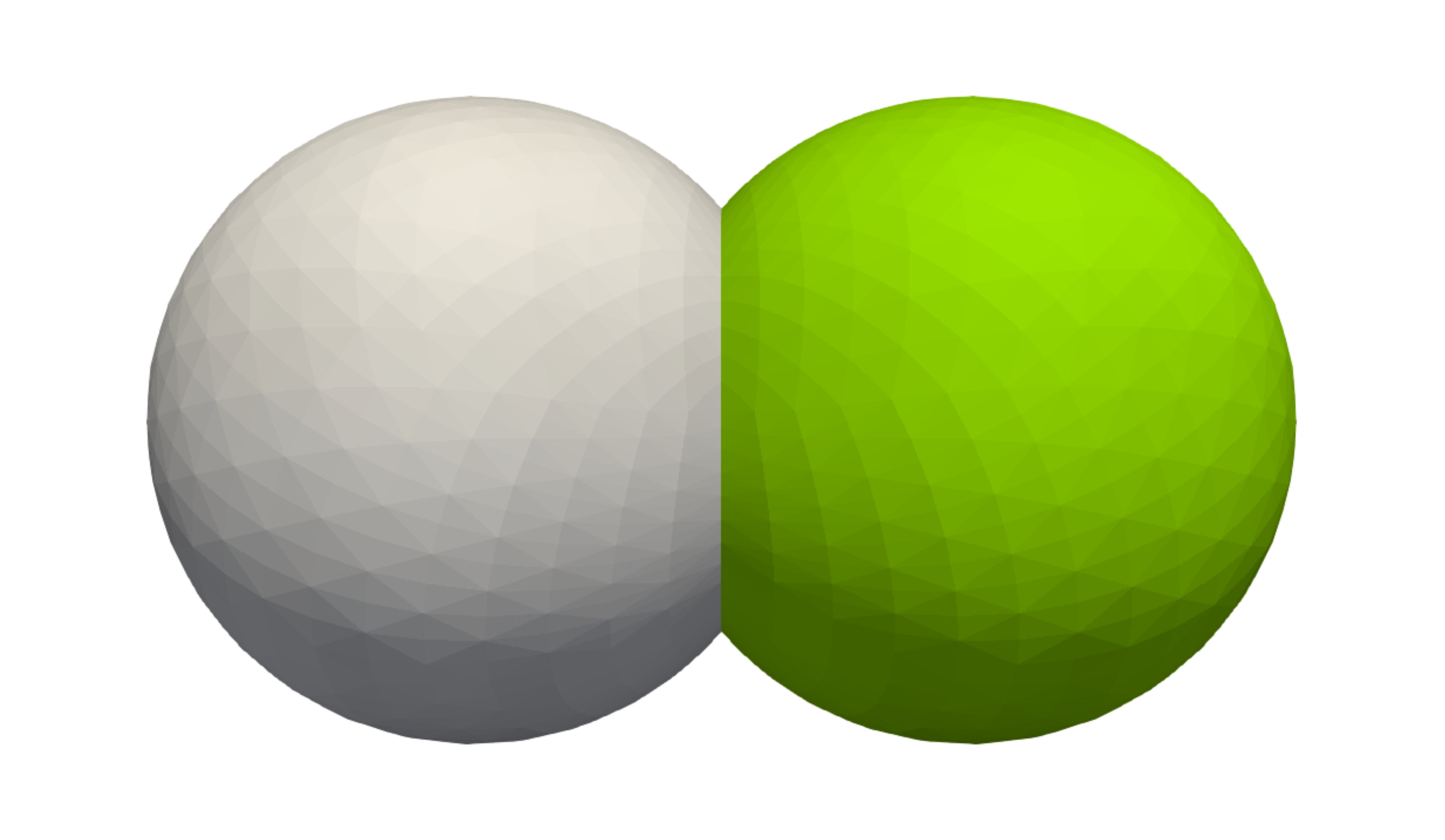}\\[0.5em]
\includegraphics[width=0.7\textwidth]{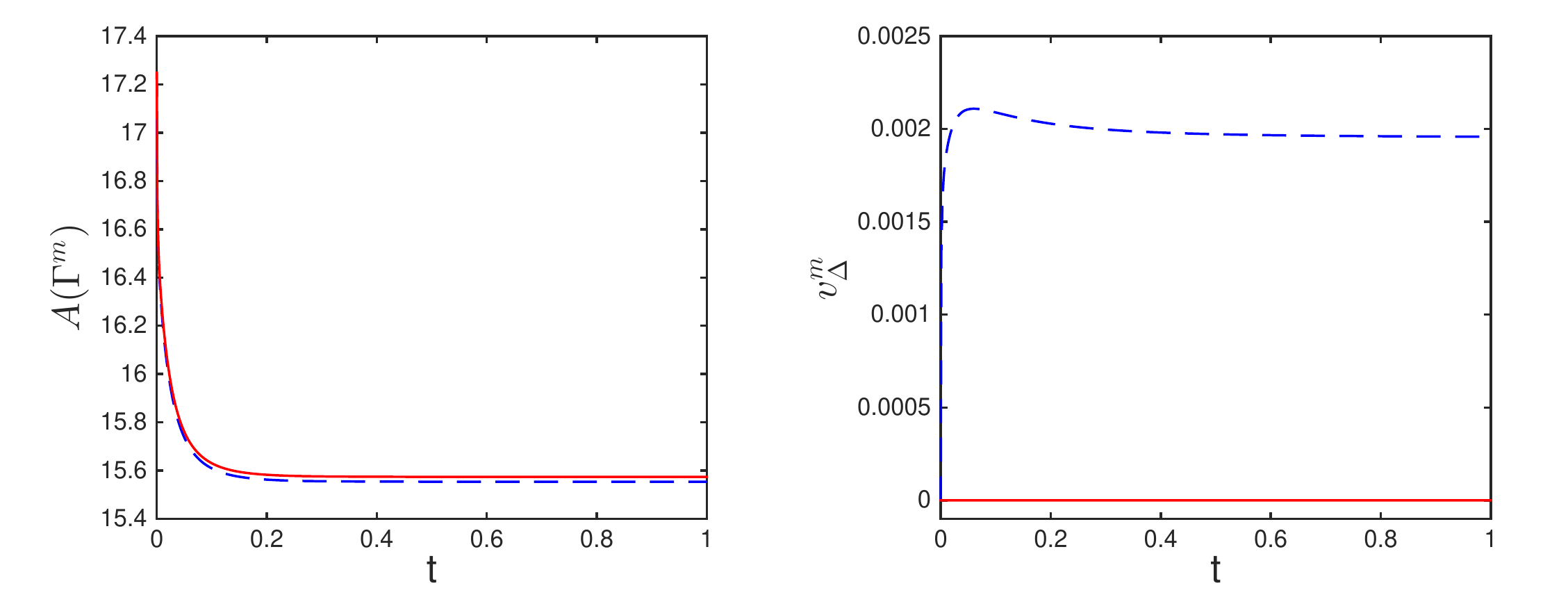}
\caption{
Evolution towards a 3d double bubble, with weightings $\sigma = (1,1.5,1)$.
Plots of $\Gamma^m$ at times $t=0, 0.1, 1$.
We also show plots of the discrete energy $A(\Gamma^m)$ and 
the relative volume error $v_\Delta^m$ over time, where $K = 3267$ and $\ttau=10^{-3}$.
}
% ~/hpc_cluster/data/alberta/tjtrue/3d.db_1weight1.5
\label{fig:3ddb1w1.5}
\end{figure}

\begin{figure}[!htb]
\center
\includegraphics[angle=-0,width=0.3\textwidth]{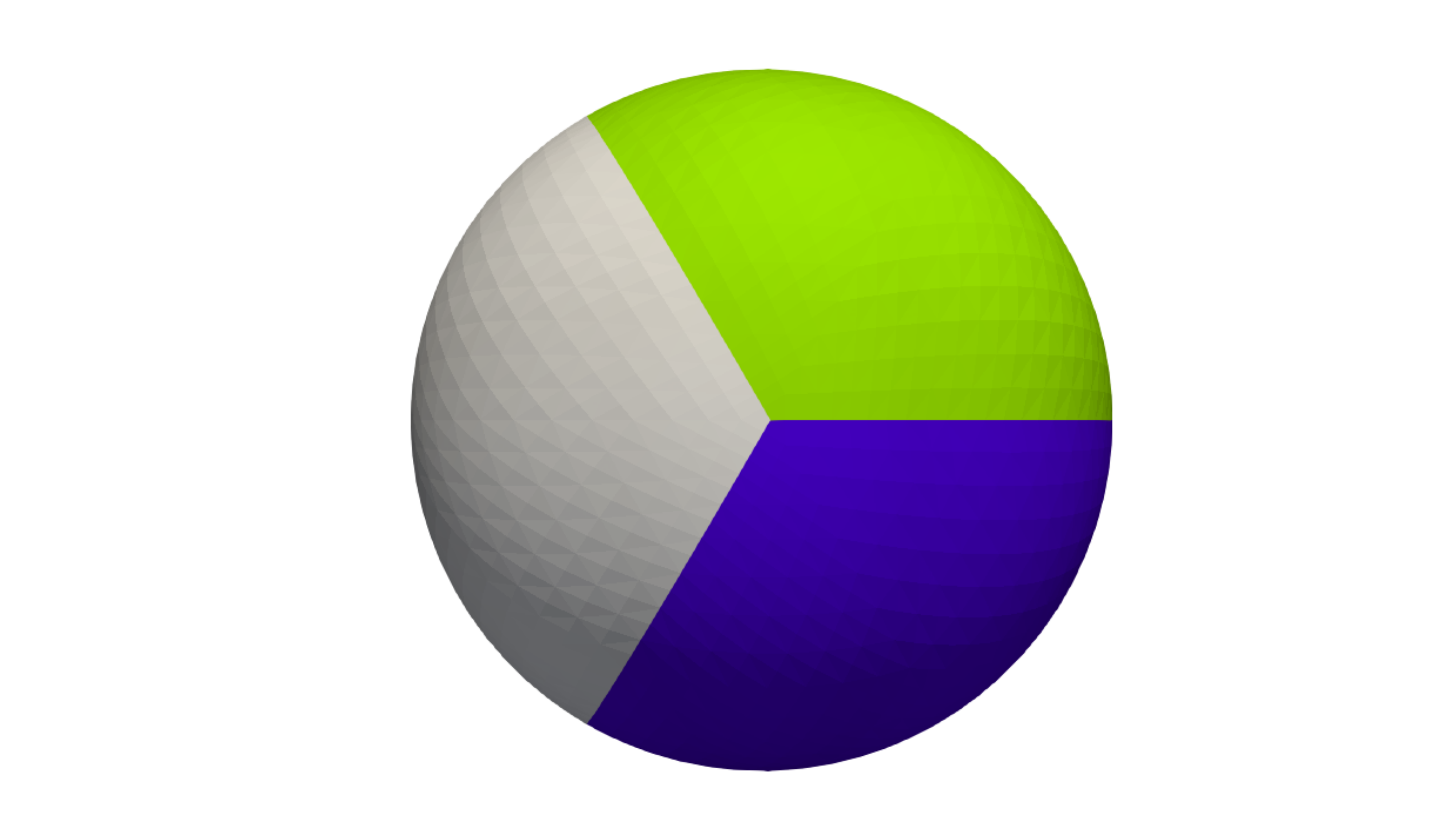}
\includegraphics[angle=-0,width=0.3\textwidth]{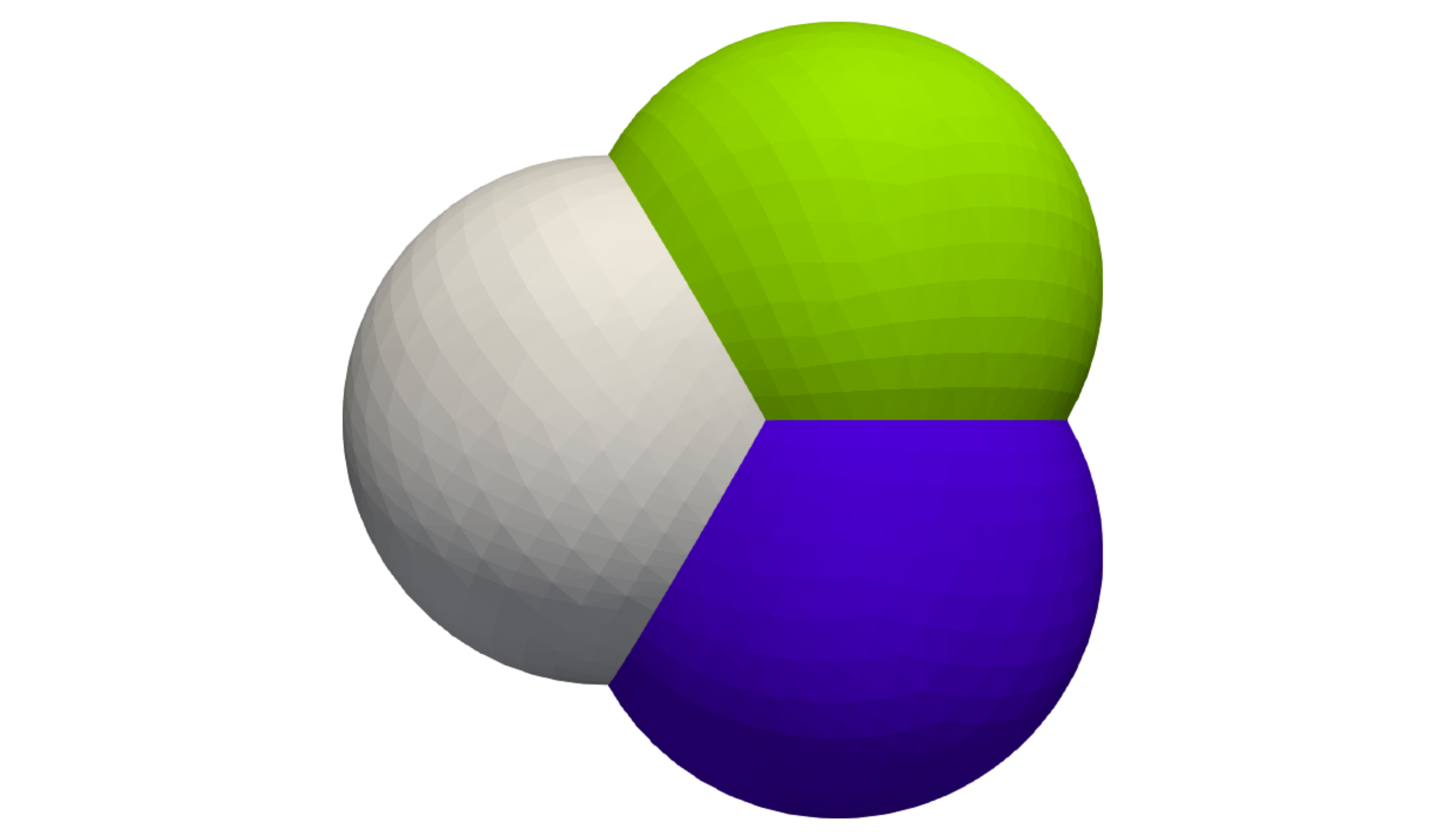}
\includegraphics[angle=-0,width=0.3\textwidth]{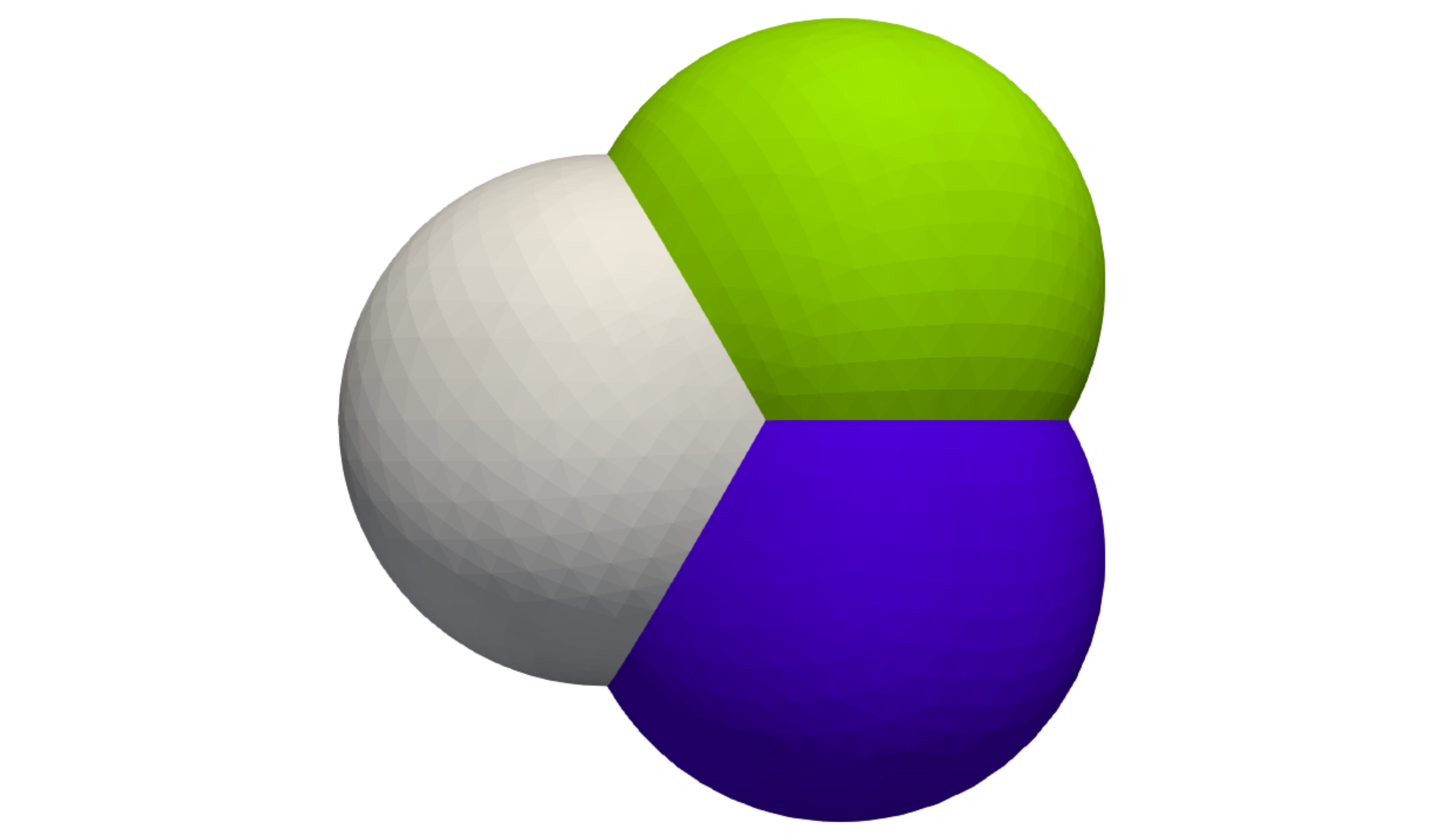}\\[0.5em]
\includegraphics[width=0.7\textwidth]{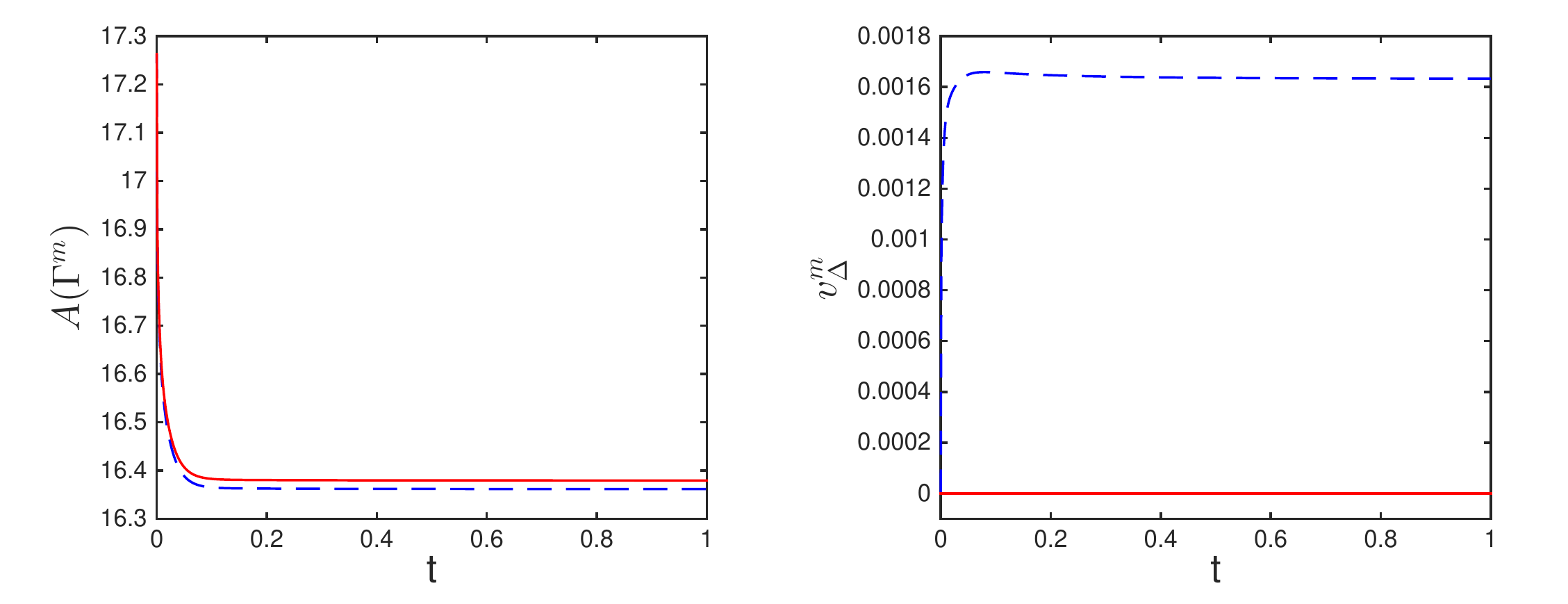}
\caption{
Evolution towards the standard 3d triple bubble.
Plots of $\Gamma^m$ at times $t=0, 0.1, 1$.
We also show plots of the discrete energy $A(\Gamma^m)$ and 
the relative volume error $v_\Delta^m$ over time, where $K = 6534$ and $\ttau=10^{-3}$.
}
\label{fig:3dtb}
\end{figure}

\begin{figure}[!htb]
\center
\includegraphics[angle=-0,width=0.3\textwidth]{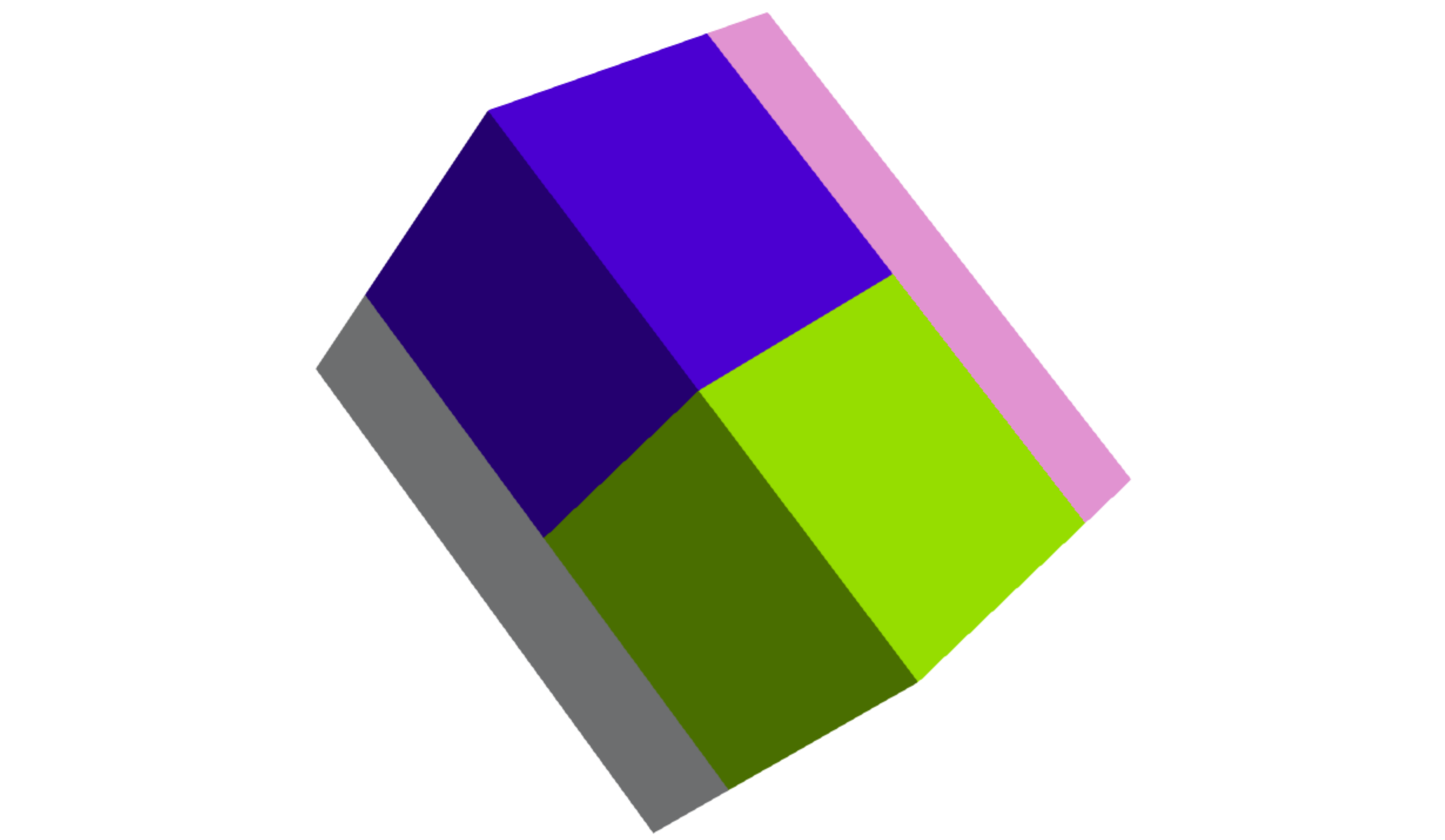}
\includegraphics[angle=-0,width=0.3\textwidth]{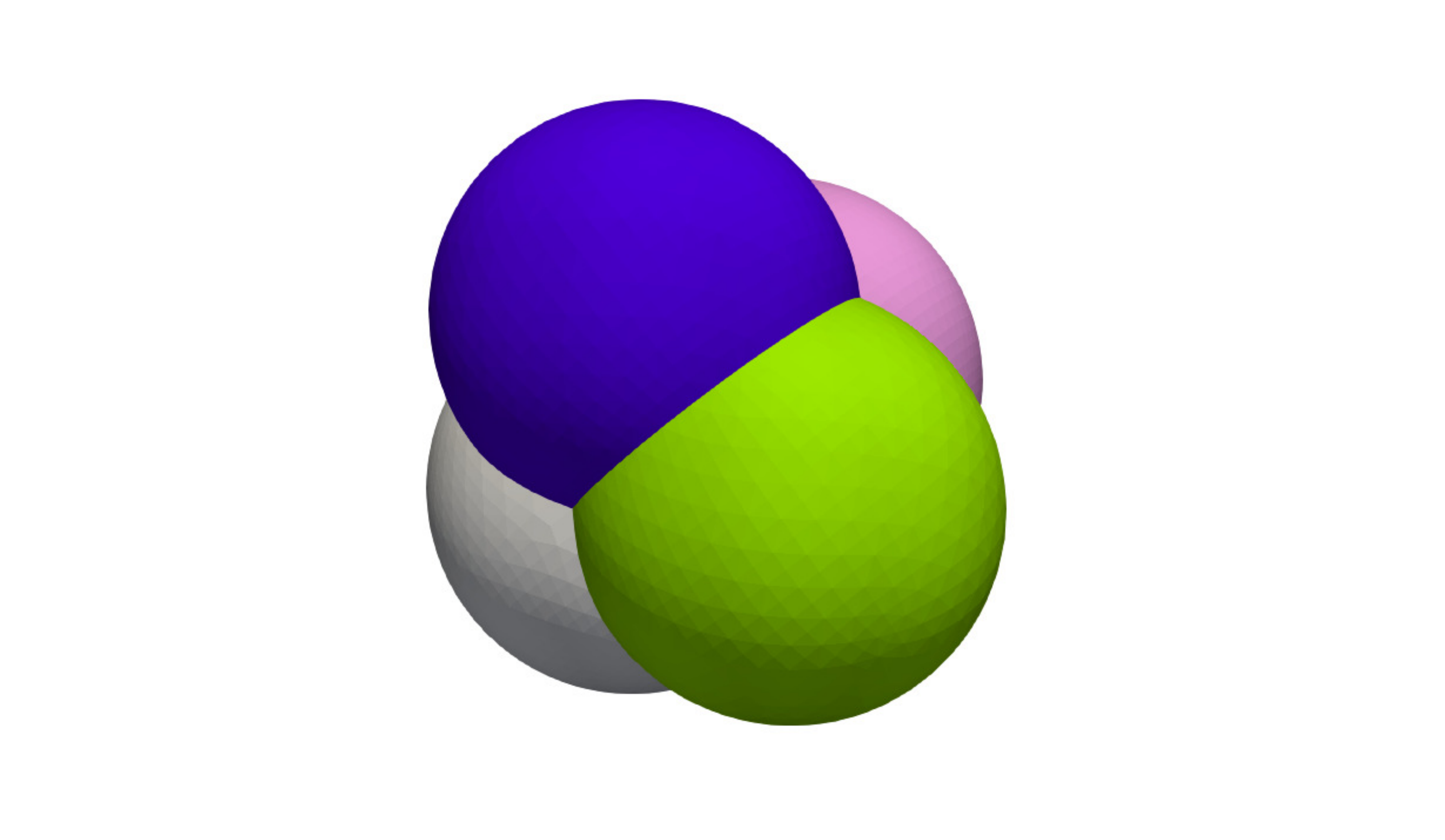}
\includegraphics[angle=-0,width=0.3\textwidth]{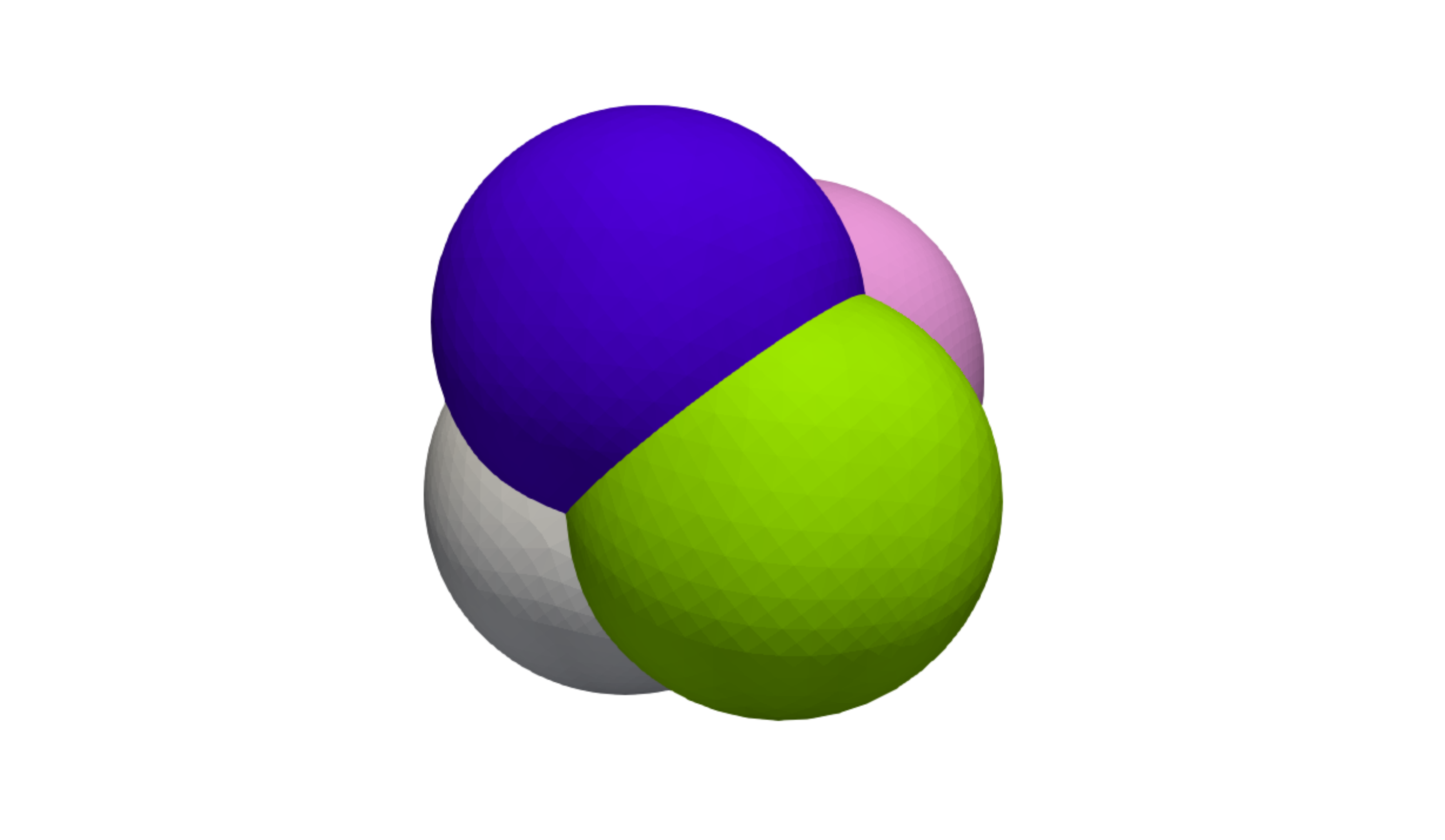}\\[0.5em]
\includegraphics[width=0.7\textwidth]{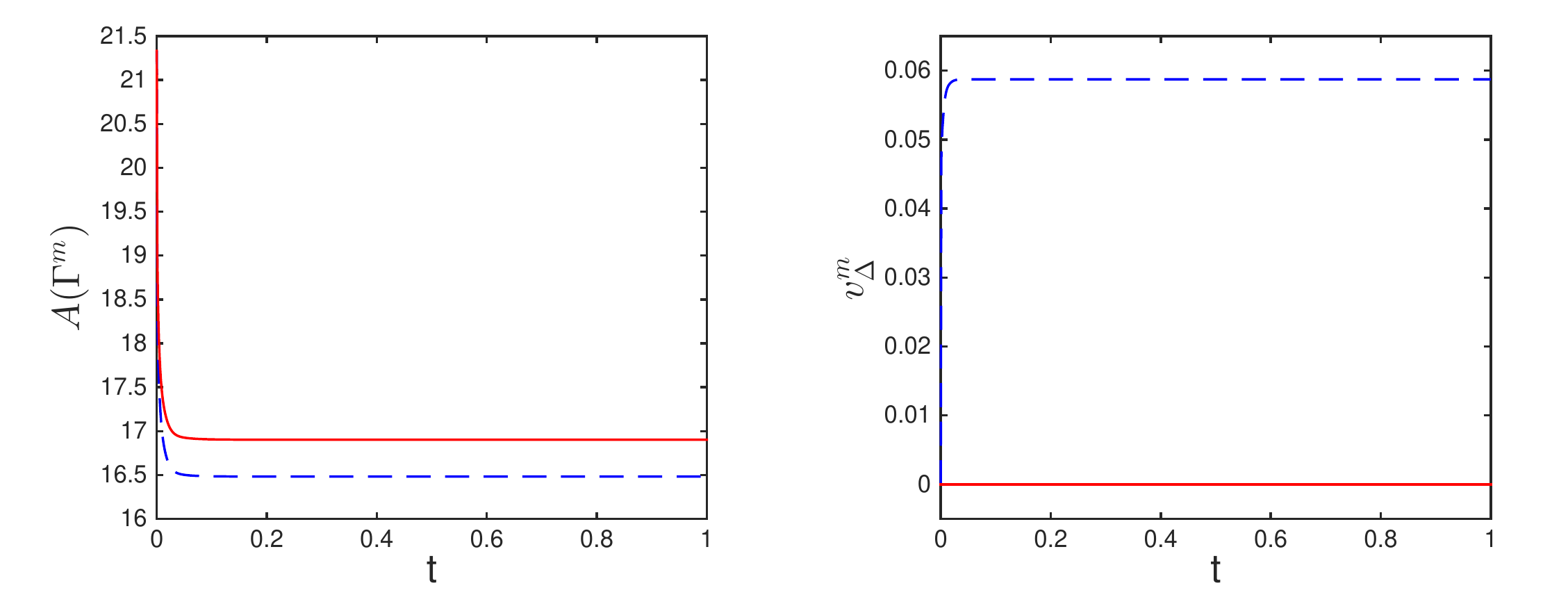}
\caption{
Evolution towards the standard 3d quadruple bubble.
Plots of $\Gamma^m$ at times $t=0, 0.1, 1$.
We also show plots of the discrete energy $A(\Gamma^m)$ and 
the relative volume error $v_\Delta^m$ over time, where $K = 8378$ and $\ttau=10^{-3}$.
}
\label{fig:3dqb}
\end{figure}

 We start with an initial surface cluster that is given by two halfspheres and a disk, meeting at a triple junction line. As shown in Fig.~\ref{fig:3ddbw1}, in the case of equal surface energy densities, we observe that the cluster evolves towards the symmetric standard double bubble, and the energy dissipation and volume conservation are well satisfied for the numerical solutions. We then use different weightings of surface energies, and the numerical results are reported in Fig.~\ref{fig:3ddbw1.5} and Fig.~\ref{fig:3ddb1w1.5}, respectively. We observe that the interface with higher weightings tends to shrink relative to the other two, thus leading to different triple junction angles. For example, in Fig.~\ref{fig:3ddb1w1.5}, the disk shrinks to form relatively large triple junction angles with the other two surfaces so that the contact angle conditions \eqref{eq:tj_b} are satisfied. Simulation results for the standard triple and quadruple bubbles are presented in Fig.~\ref{fig:3dtb} and \ref{fig:3dqb}, respectively. Regardless of the different setups, we can always observe the dissipation of the total surface area and the exact volume conservation for each enclosed bubble in these experiments.

\begin{figure}[!htb]
\center
\includegraphics[angle=-0,width=0.3\textwidth]{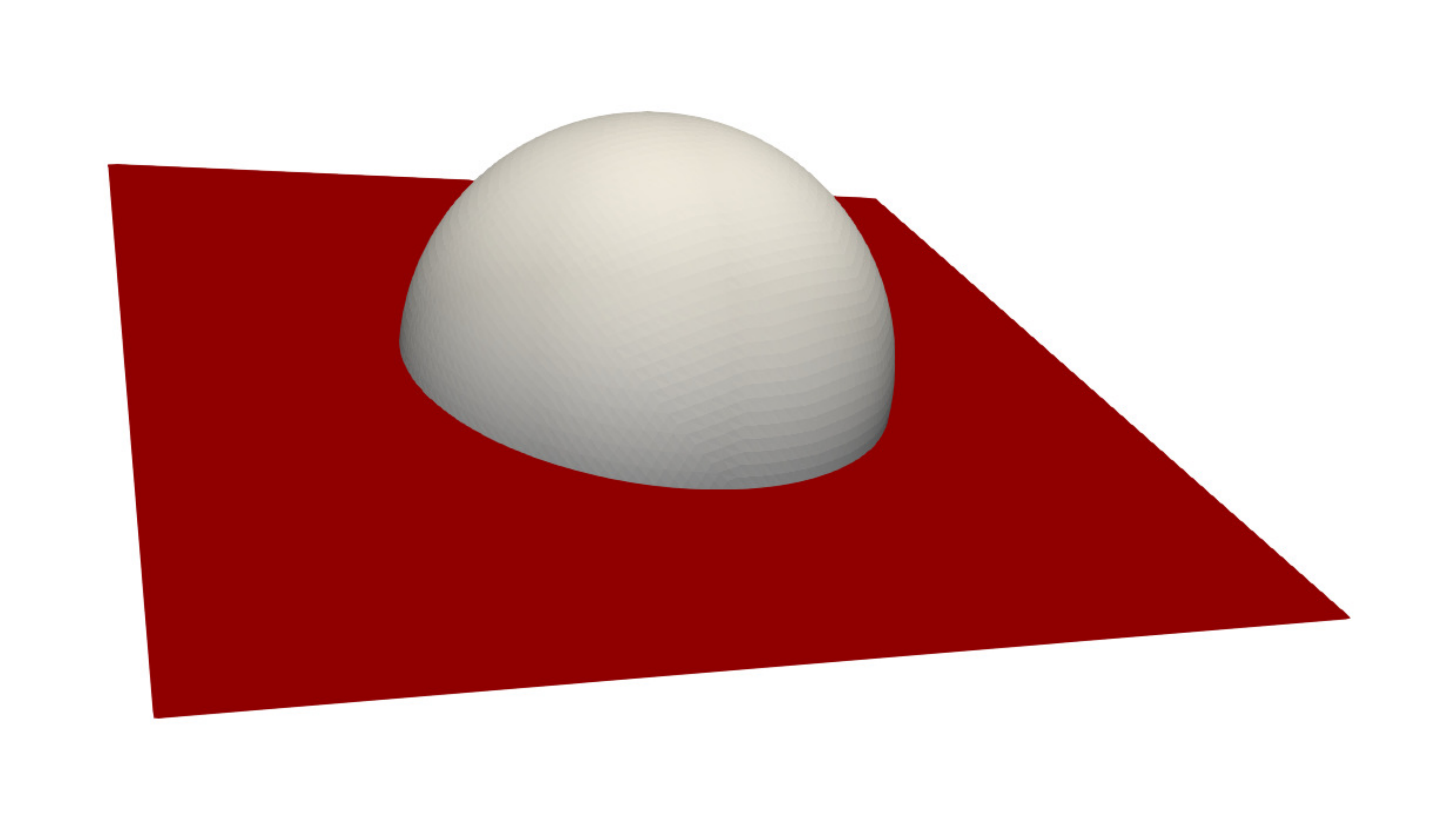}
\includegraphics[angle=-0,width=0.3\textwidth]{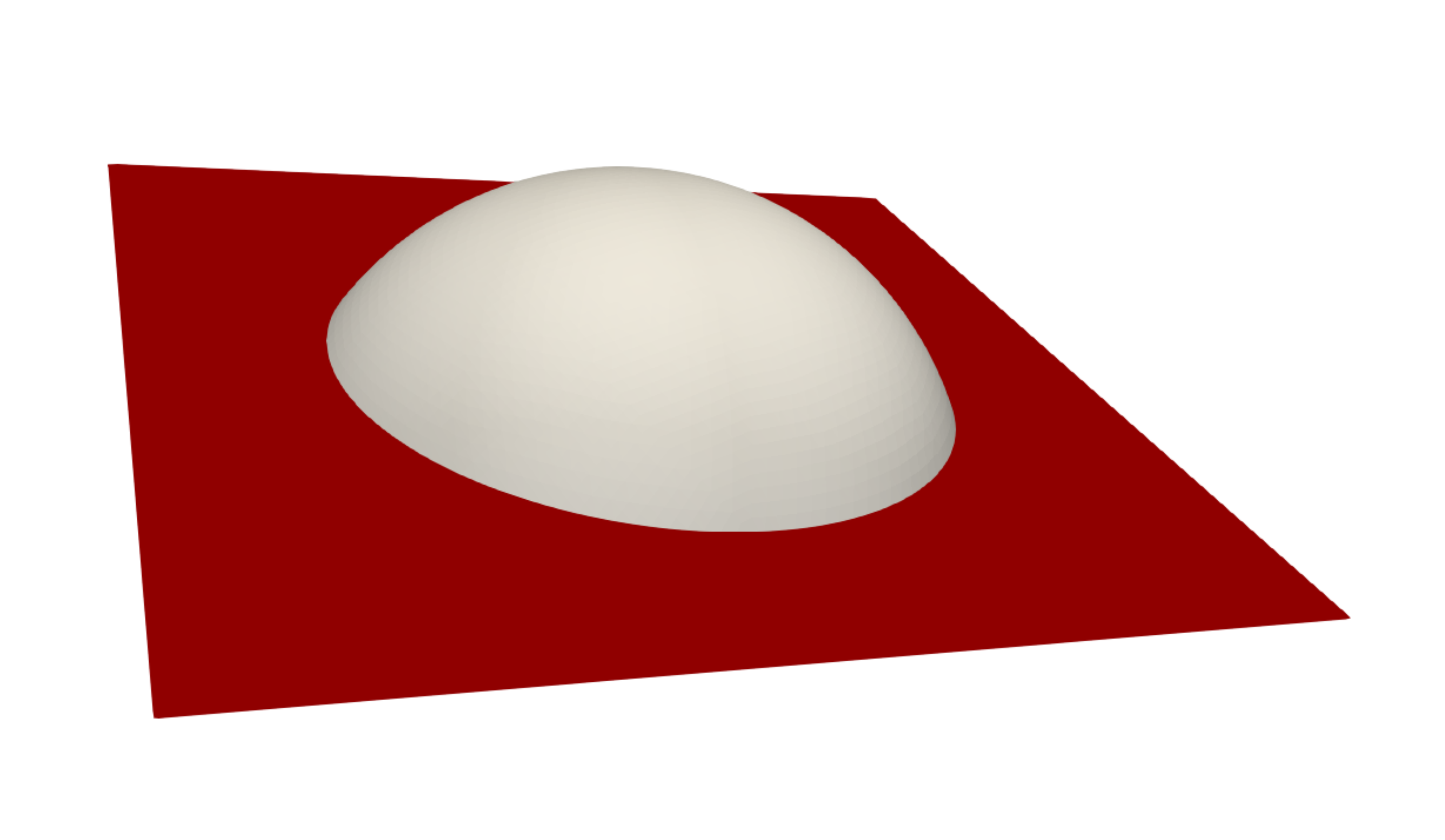} \\[0.5em]
\includegraphics[width=0.65\textwidth]{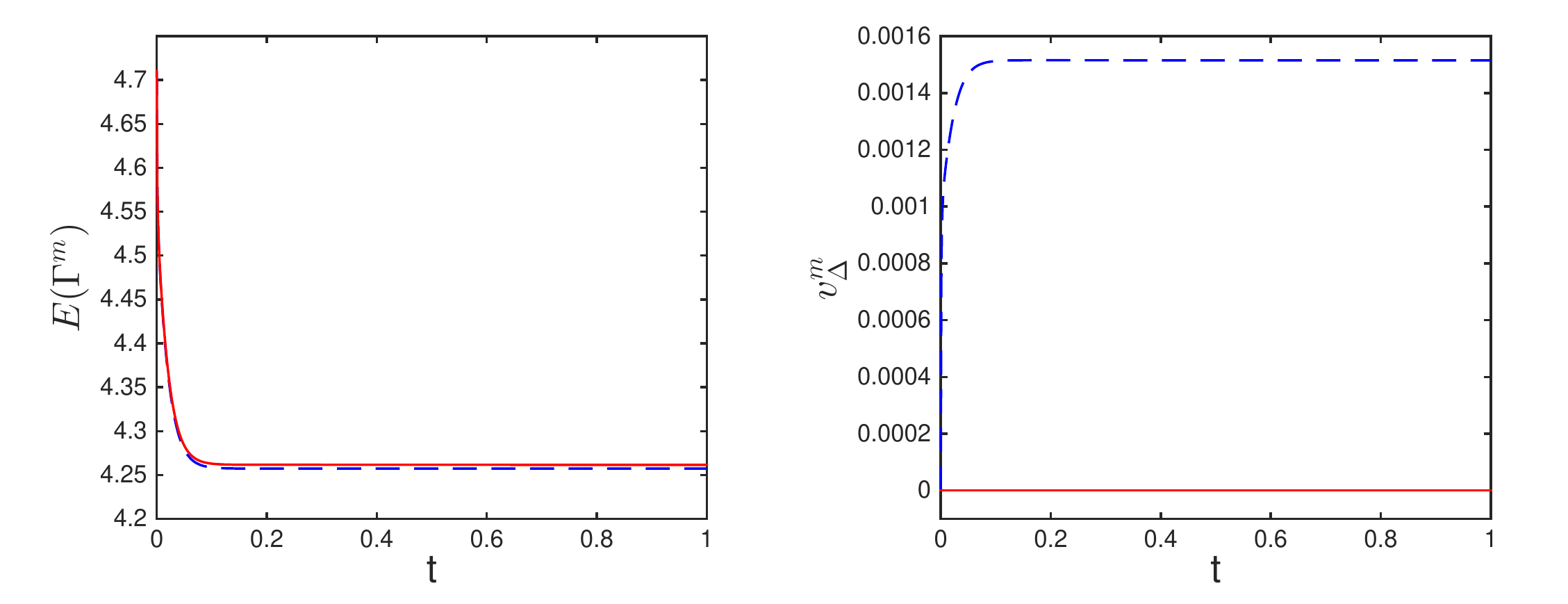}
\caption{
Evolution towards a drop on a substrate, with $\varrho=0.5$ so that
$\vartheta=60^\circ$.
Plots of $\Gamma^m$ at times $t=0, 1$.
We also show plots of the discrete energy $E(\Gamma^m)$ and 
the relative volume error $v_\Delta^m$ over time, where
$K = 4225$ and $\Delta t = 10^{-3}$.
}
\label{fig:3disodrop05}
\end{figure}%

\begin{figure}[!htb]
\center
\includegraphics[angle=-0,width=0.3\textwidth]{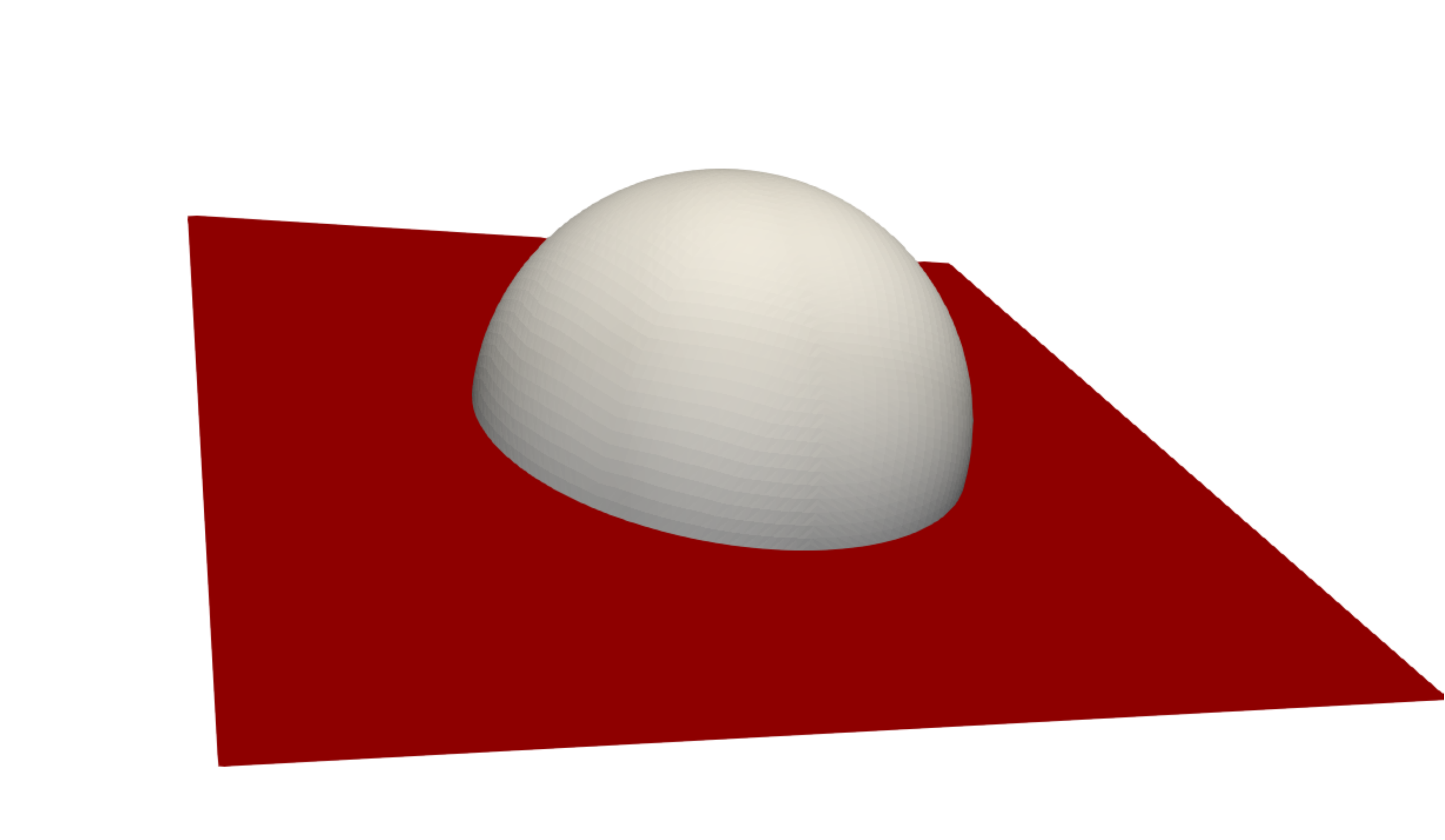}
\includegraphics[angle=-0,width=0.3\textwidth]{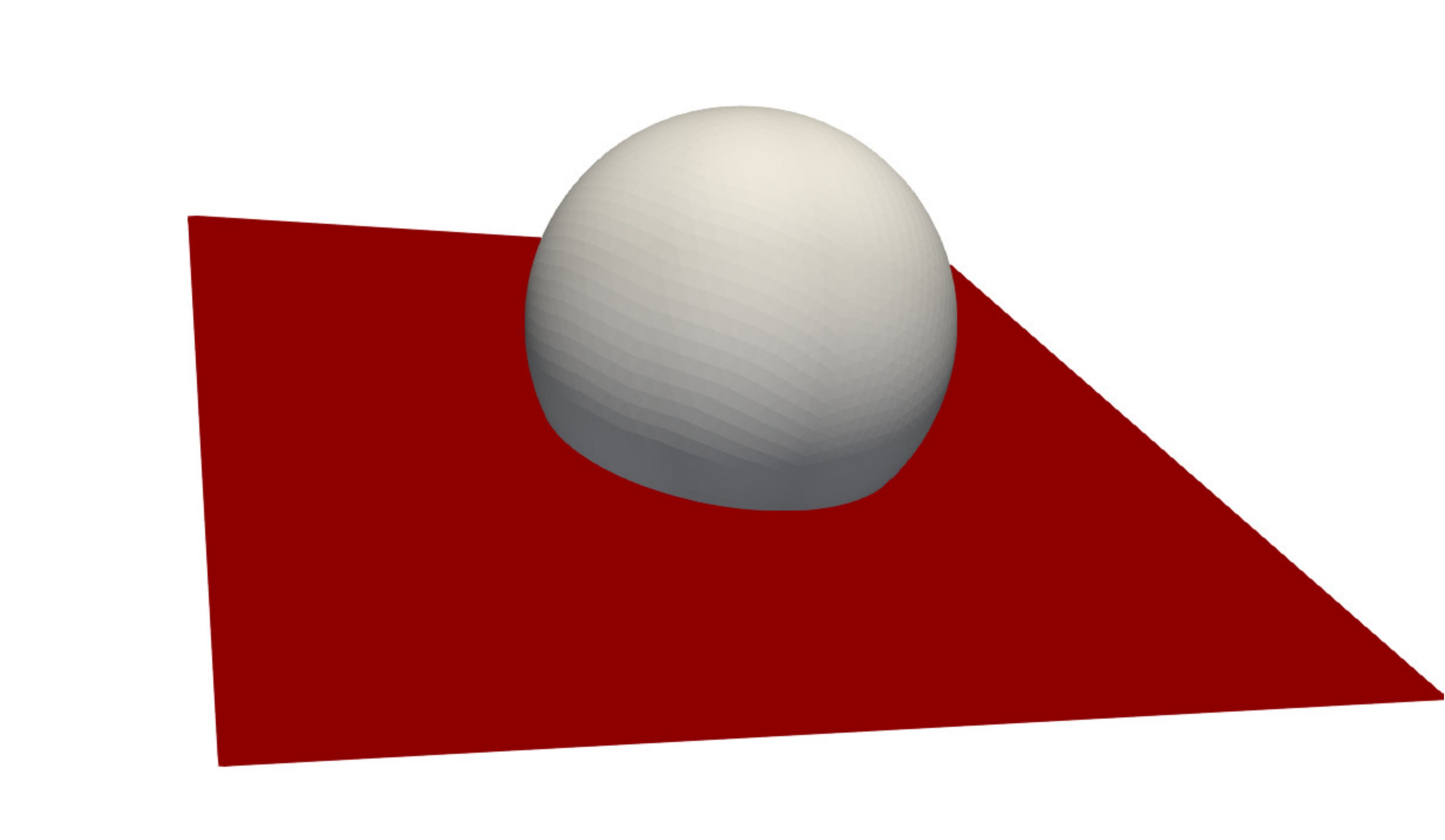} \\[0.5em]
\includegraphics[width=0.65\textwidth]{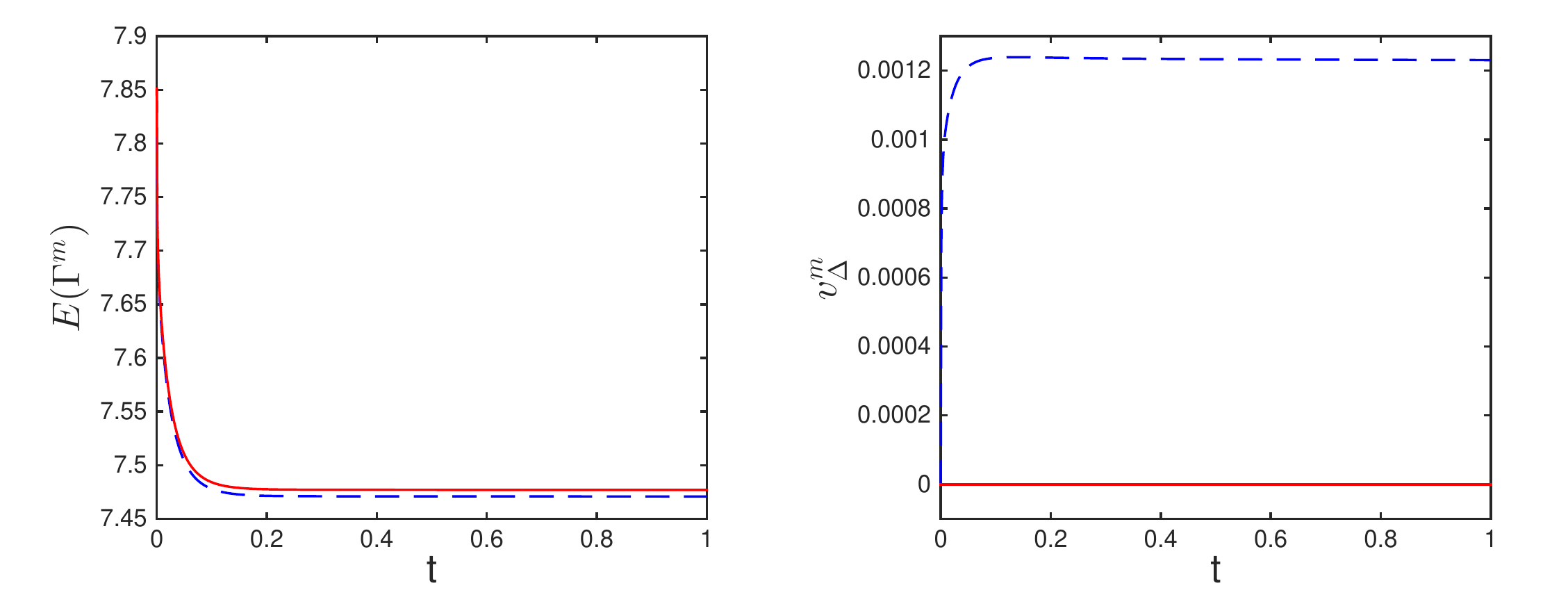}
\caption{
Evolution towards a drop on a substrate, with $\varrho=-0.5$ so that
$\vartheta=120^\circ$.
Plots of $\Gamma^m$ at times $t=0, 1$.
We also show plots of the discrete energy $E(\Gamma^m)$ and 
the relative volume error $v_\Delta^m$ over time, where
$K = 4225$ and $\Delta t = 10^{-3}$.
}
% ~/hpc_cluster/data/alberta/tjtrue/3d.isoangle_rho-05
\label{fig:3disodrop-05}
\end{figure}%

We then simulate the evolution of a single drop which is attached to a non-neutral substrate $\mathcal{D}_1=\{(q_1,q_2,q_3)\in\bR^3:q_3=0\}$, and initially the drop is chosen as a semisphere. The numerical results for $\varrho=0.5$ and $\varrho=-0.5$ are shown in
Figs.~\ref{fig:3disodrop05} and \ref{fig:3disodrop-05}, respectively. We can observe that the drop finally maintains the steady state with a contact angle 
of about $60^\circ$ when $\varrho=0.5$, and a contact angle of about
$120^\circ$ when $\varrho=-0.5$.

\begin{figure}[!htp]
\center
\includegraphics[angle=-0,width=0.3\textwidth]{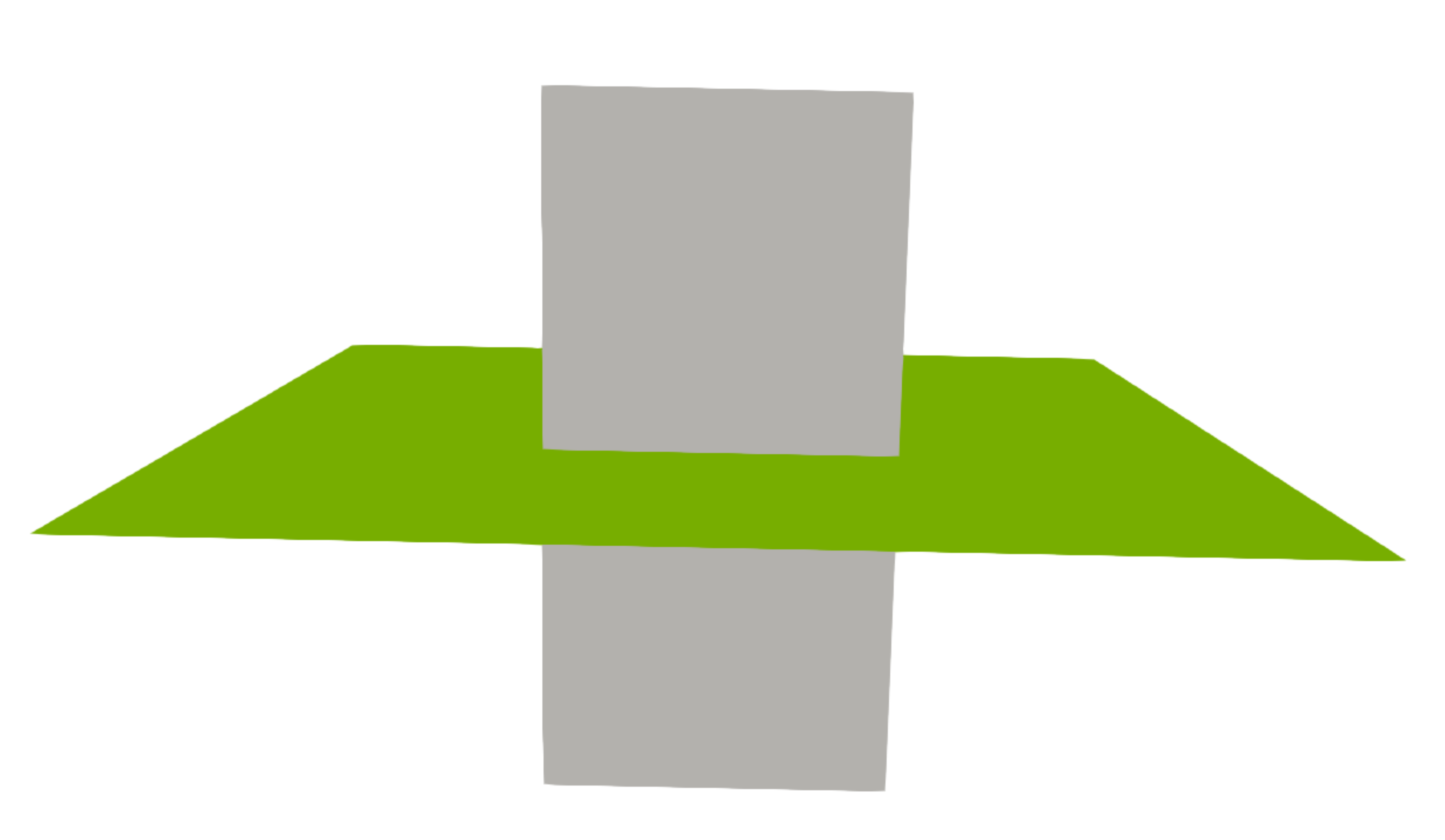}
\includegraphics[angle=-0,width=0.3\textwidth]{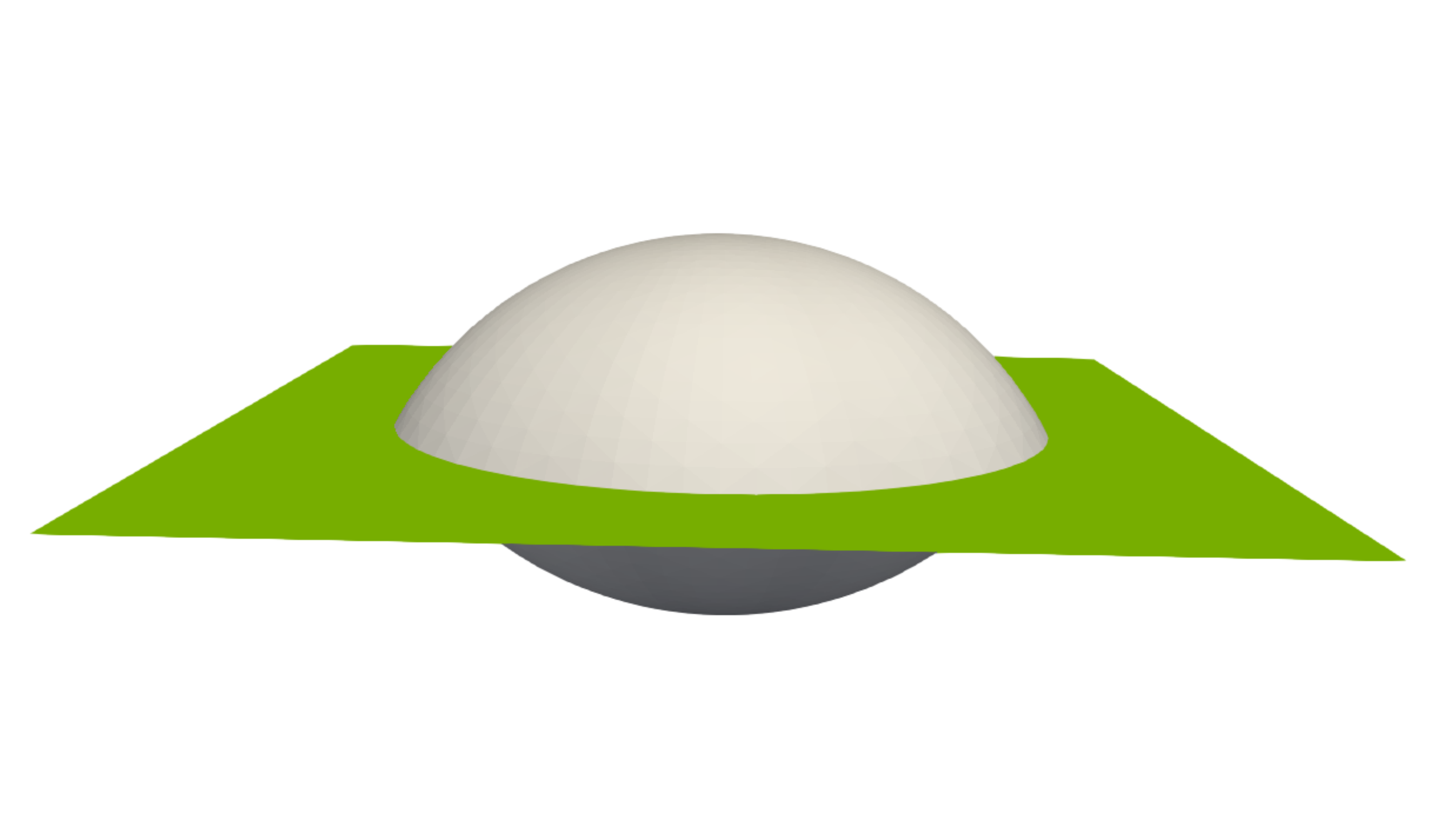} \\[0.5em]
\includegraphics[width=0.65\textwidth]{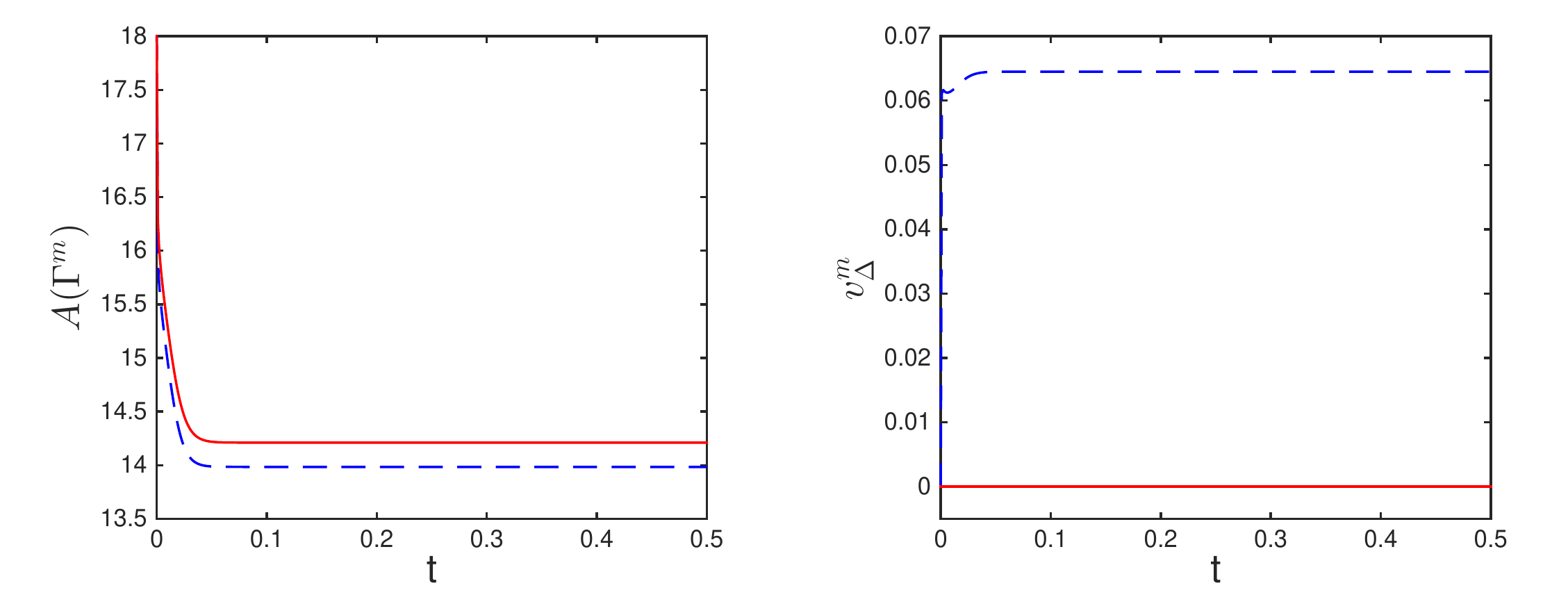}
\caption{
Plots of $\Gamma^m$ at times $t=0, 0.5$, with $\varrho=0$ so that
$\vartheta=90^\circ$.
We also show plots of the discrete energy $A(\Gamma^m)$ and 
the relative volume error $v_\Delta^m$ over time, where $K = 4802$ and $\Delta t = 10^{-3}$.
}
\label{fig:3dgrain}
\end{figure}

\begin{figure}[!htp]
\center
\includegraphics[angle=-0,width=0.3\textwidth]{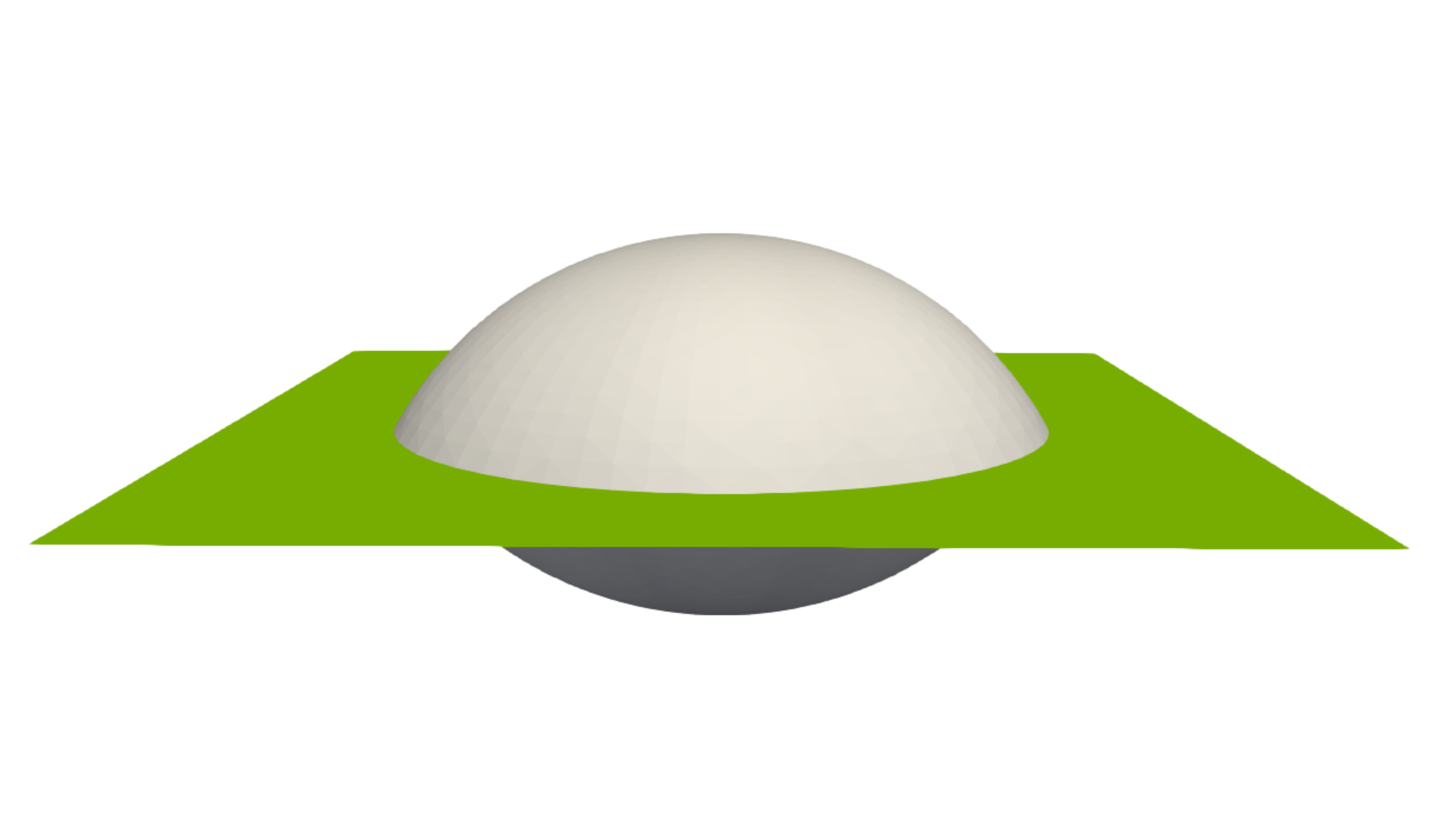}
\includegraphics[angle=-0,width=0.3\textwidth]{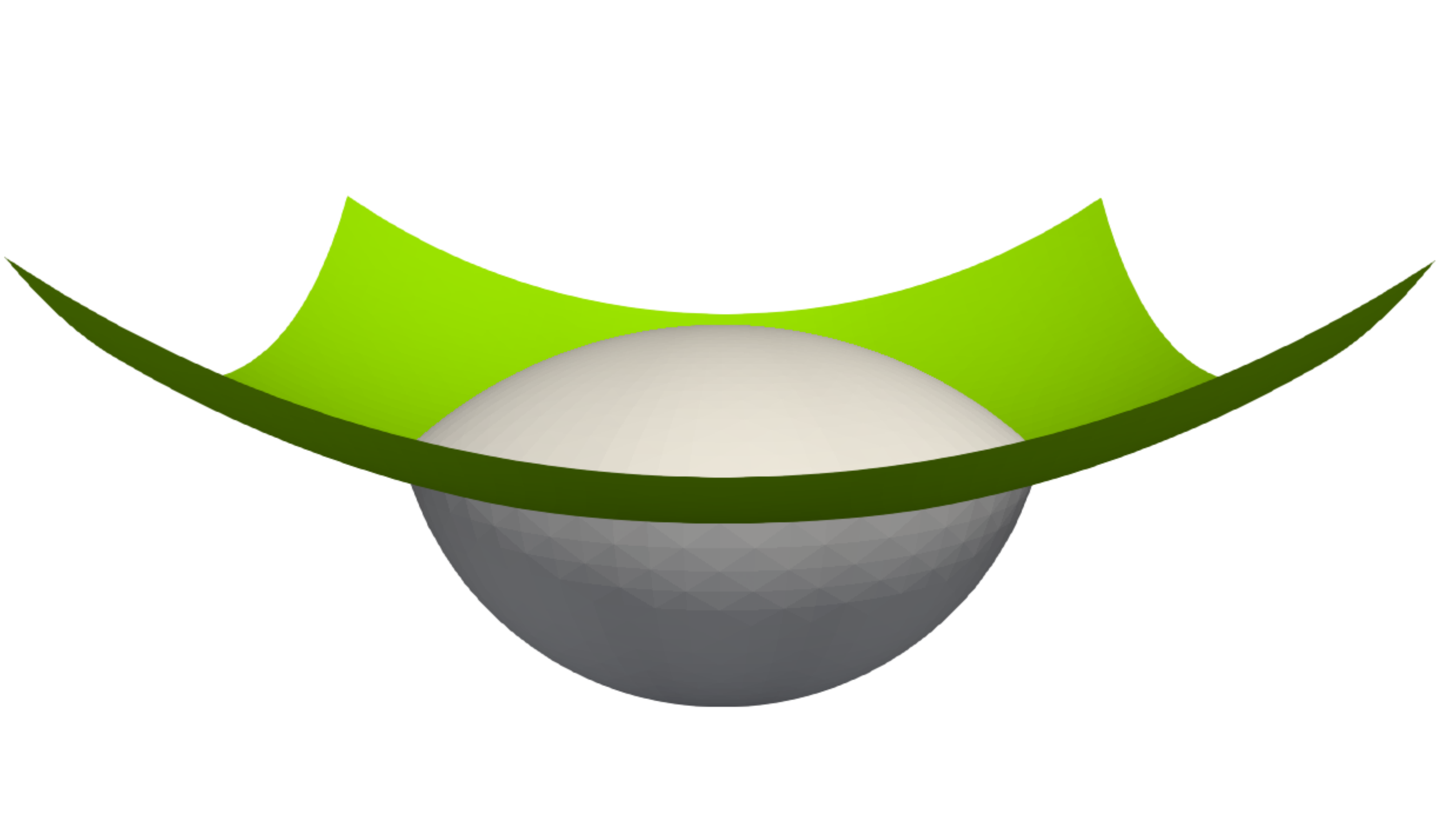} 
\caption{
Plots of $\Gamma^m$ at times $t=0, 0.5$, with $\varrho=0.5$ so that
$\vartheta=60^\circ$.
We also have $K = 4802$ and $\Delta t = 10^{-3}$.
}
\label{fig:3dgrainsmooth05}
\end{figure}

\begin{figure}[!htp]
\center
\includegraphics[angle=-0,width=0.4\textwidth]{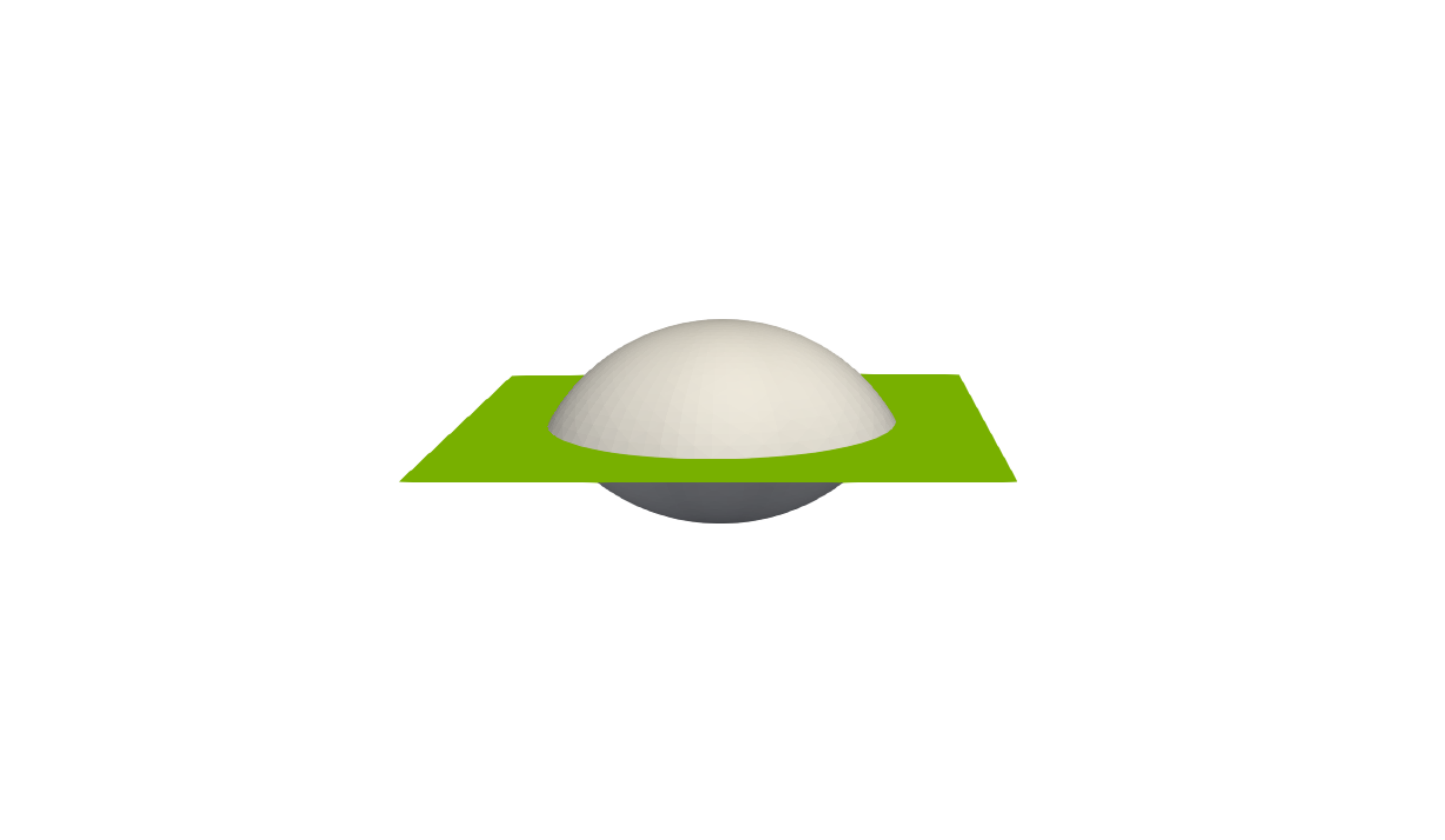}
\includegraphics[angle=-0,width=0.4\textwidth]{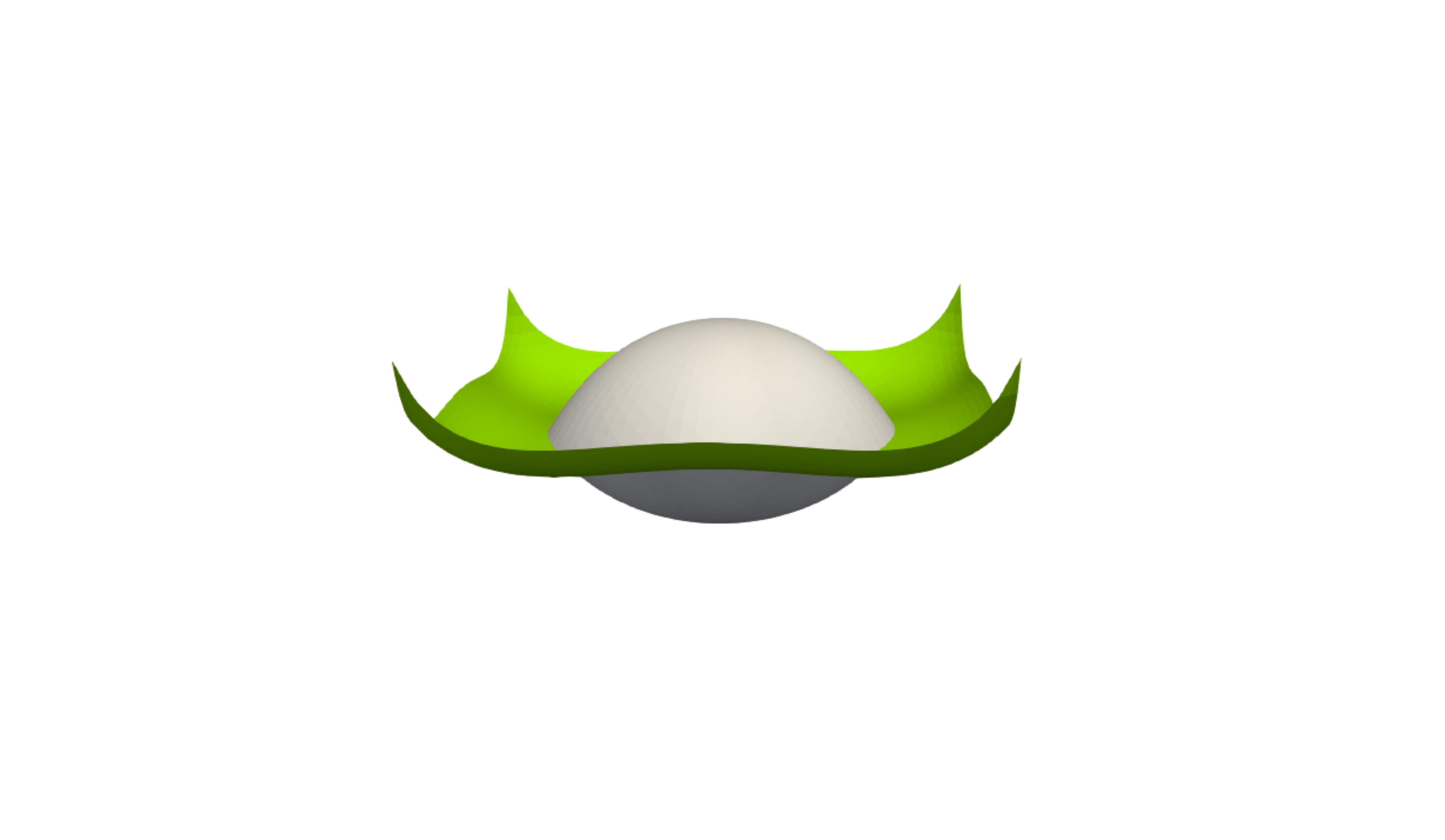}
\includegraphics[angle=-0,width=0.4\textwidth]{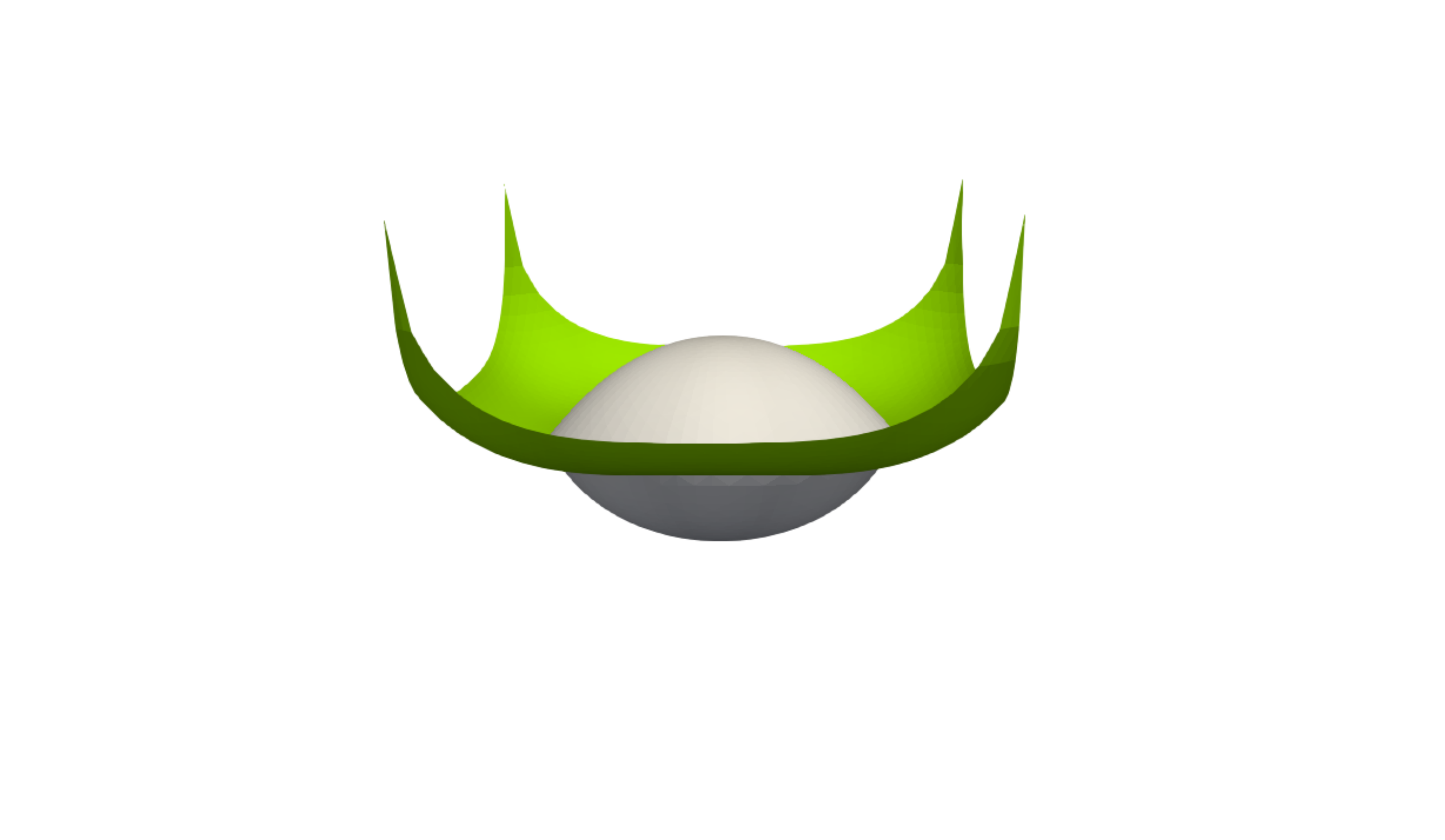}
\includegraphics[angle=-0,width=0.4\textwidth]{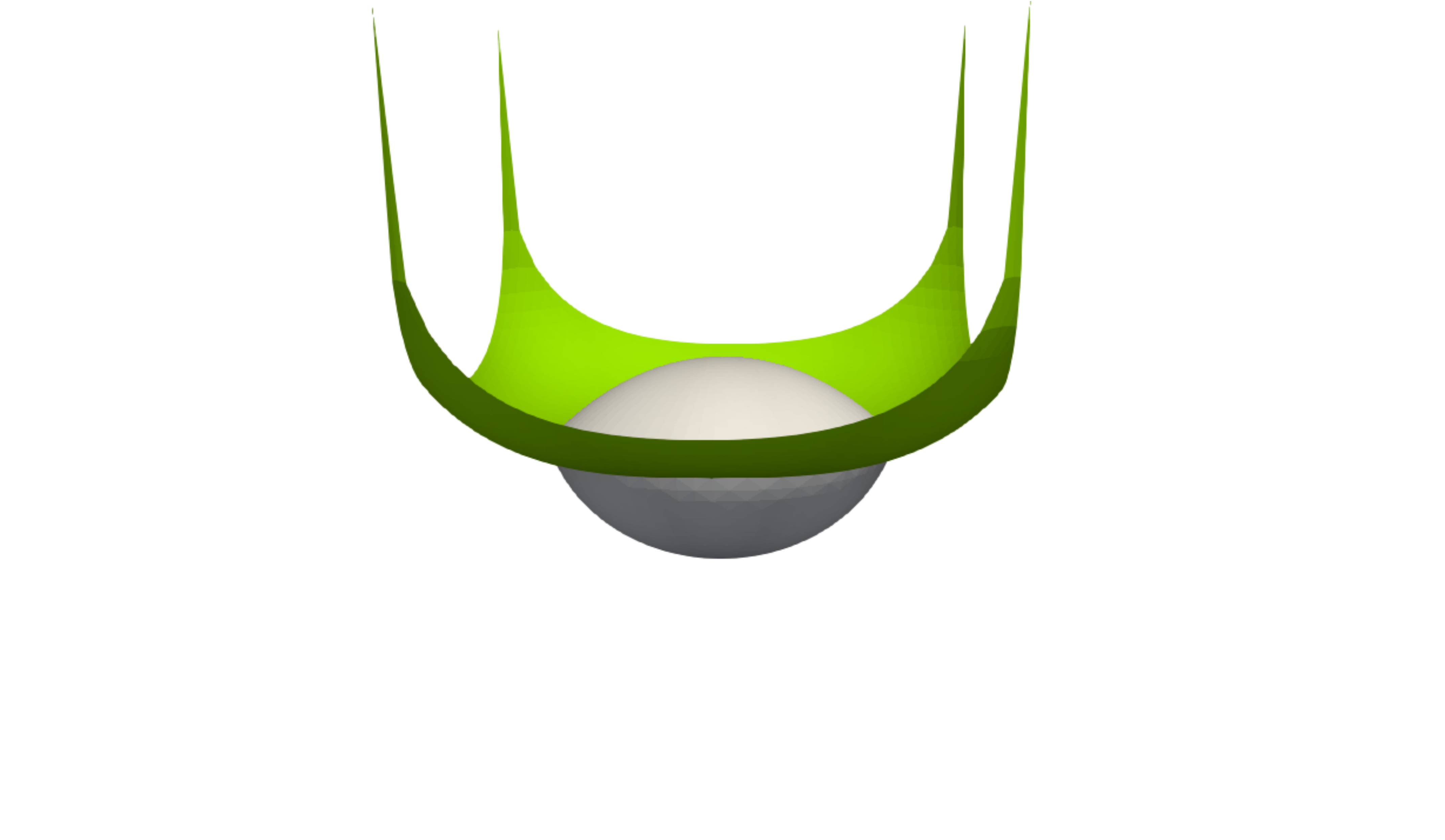} 
\caption{
Plots of $\Gamma^m$ at times $t=0, 0.01, 0.05, 0.1$, with $\varrho=0.75$, 
so that $\vartheta=41.4^\circ < 45^\circ$.
We also have $K = 4802$ and $\Delta t = 10^{-3}$.
}
\label{fig:3dgrainsmooth075}
\end{figure}

We next test the evolution of a surface cluster contained in a cylinder of square cross-section. As shown in Fig.~\ref{fig:3dgrain}, the cluster is made up of three surfaces, meeting at a triple junction line, and with one of the surfaces (coloured in green) attached to the external boundary of the cylinder 
$[-\frac{3}{2},\frac{3}{2}]^2\times\bR$. This gives rise to four boundary lines on the four planar boundaries. In the case when $\varrho=0$, we observe that the two surfaces of the initial cuboid remain symmetric and become spherical, and the third surface remains flat and attached orthogonally to the external boundaries.  We then start from the steady state in Fig.~\ref{fig:3dgrain} and consider different boundary energy contributions. 
When $\varrho=0.5$, as shown in Fig.~\ref{fig:3dgrainsmooth05}, the cluster forms a steady state with a contact angle of about $60^\circ$ at the external boundary. Observe that the central bubble is now no longer symmetric.
Increasing the value of the boundary energy contribution to $\varrho=0.75$
yields the results in Fig.~\ref{fig:3dgrainsmooth075}. Here we observe
an unbounded growth of the initially flat surface towards infinity, 
reminiscent of the NASA
experiments in zero gravity discussed in e.g.\ \cite{ConcusF74} 
and \cite[Chapter~6]{Finn86}. In fact, for the chosen value of $\varrho=0.75$,
the preferred contact angle is $41.4^\circ$, which is outside the range
$[45^\circ,135^\circ]$ for which it is known that a finite minimizer exists.

\subsection{Anisotropic numerical results in 3d} 

\begin{figure}[!htp]
\centering
\includegraphics[angle=-0,width=0.3\textwidth]{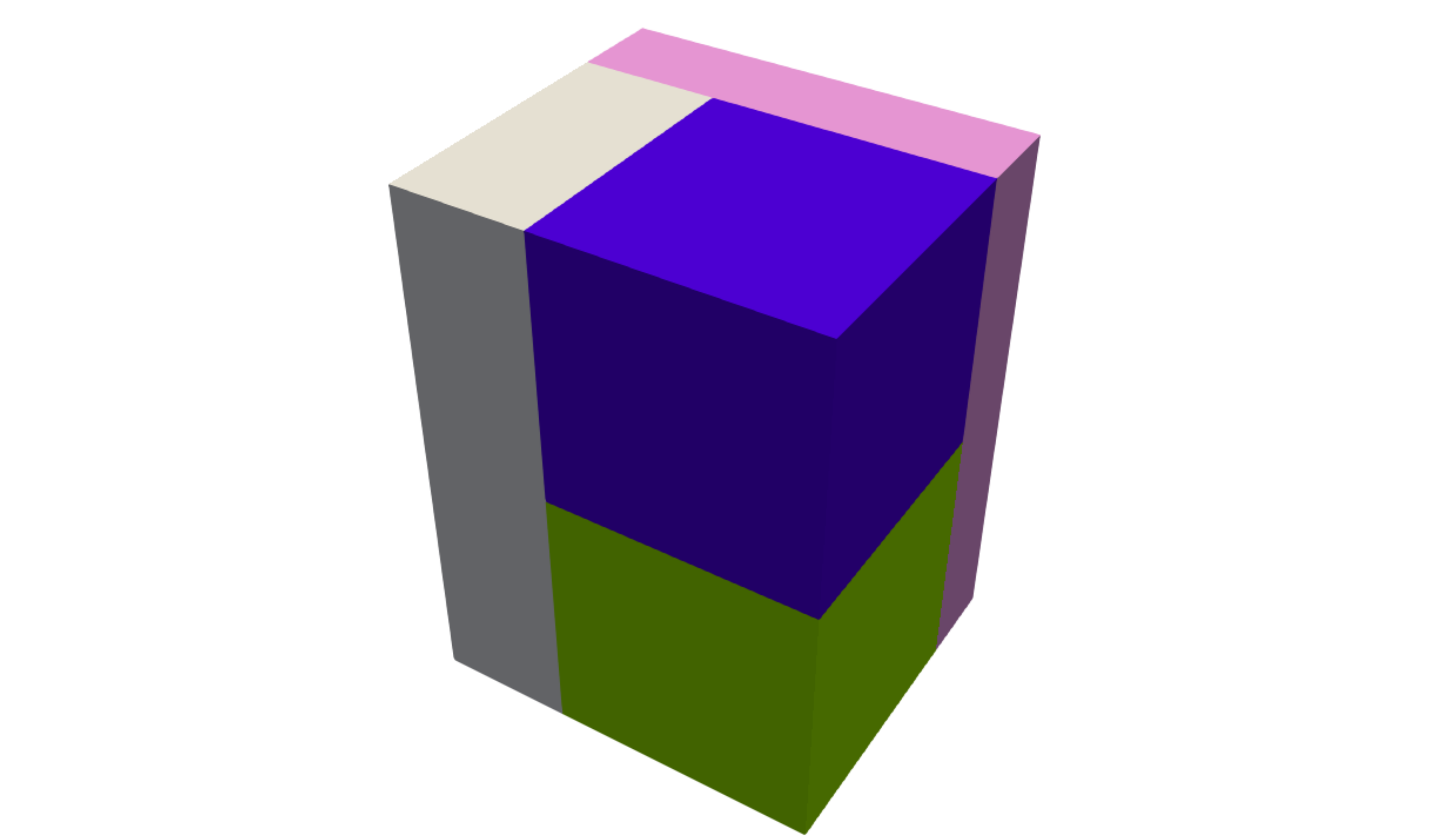}
\includegraphics[angle=-0,width=0.3\textwidth]{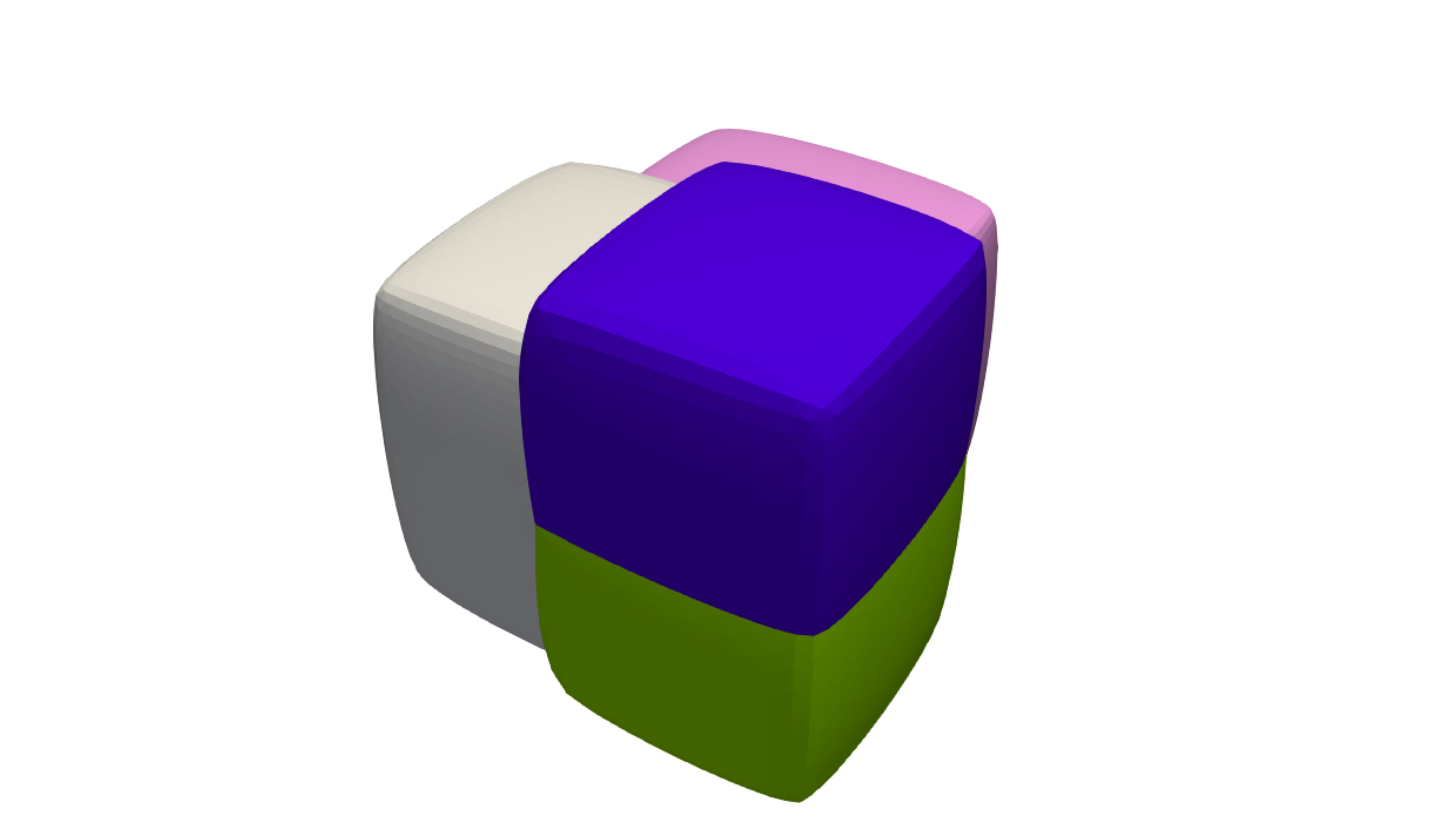}
\includegraphics[angle=-0,width=0.3\textwidth]{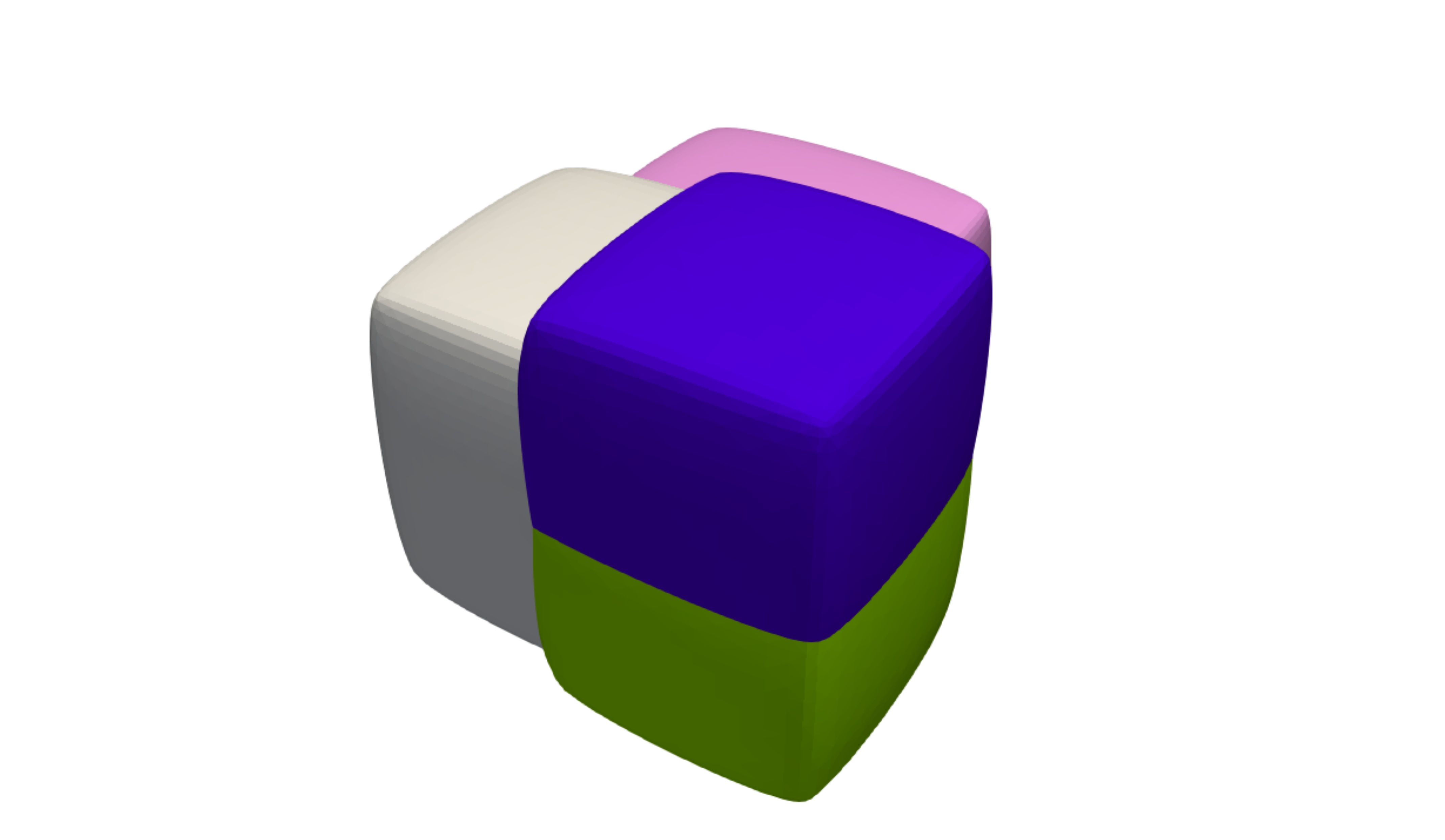}\\[0.5em]
\includegraphics[width=0.7\textwidth]{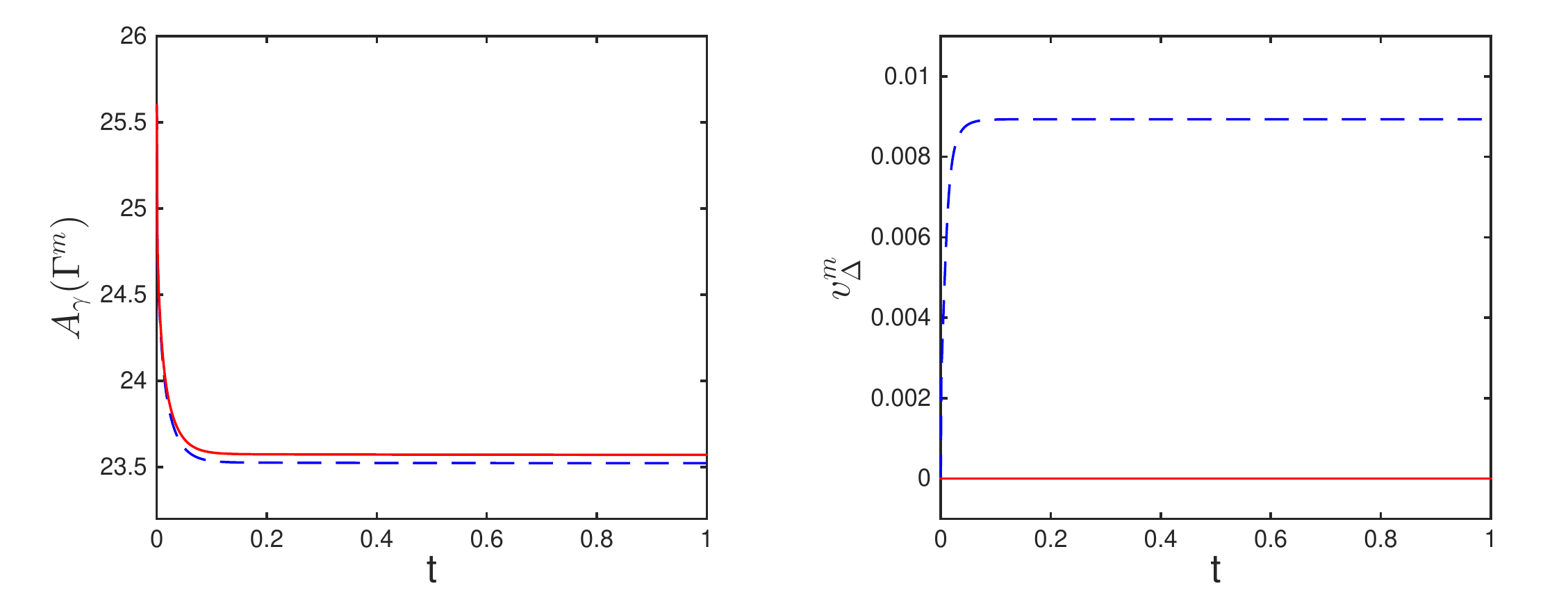}
\caption{
Evolution towards an anisotropic 3d quadruple bubble, for the 
anisotropy \eqref{eq:cuspgamma} with $L=3$, $r=1$ and $\epsilon=0.1$.
Plots of $\Gamma^m$ at times $t=0, 0.1, 1$.
We also show plots of the discrete energy $A_\gamma(\Gamma^m)$ and 
the relative volume error $v_\Delta^m$ over time, where $K = 8378$ and $\ttau = 10^{-3}$.
}
\label{fig:3dqbL3e01}
\end{figure}%

To observe the anisotropic effects, we repeat the experiment in Fig.~\ref{fig:3dqb} for the 3d quadruple bubble and use the smoothed $l^1$--norm anisotropy in \eqref{eq:cuspgamma} with $L=3, r=1$
and $\epsilon=0.1$. The numerical results are shown in Fig.~\ref{fig:3dqbL3e01}, where we find that the surfaces evolve into near cuboid shapes instead of spherical shapes as the steady state. During the simulations, the energy dissipation and volume conservation for the numerical solutions are observed as well.

\begin{figure}[!htp]
\center
\includegraphics[angle=-0,width=0.35\textwidth]{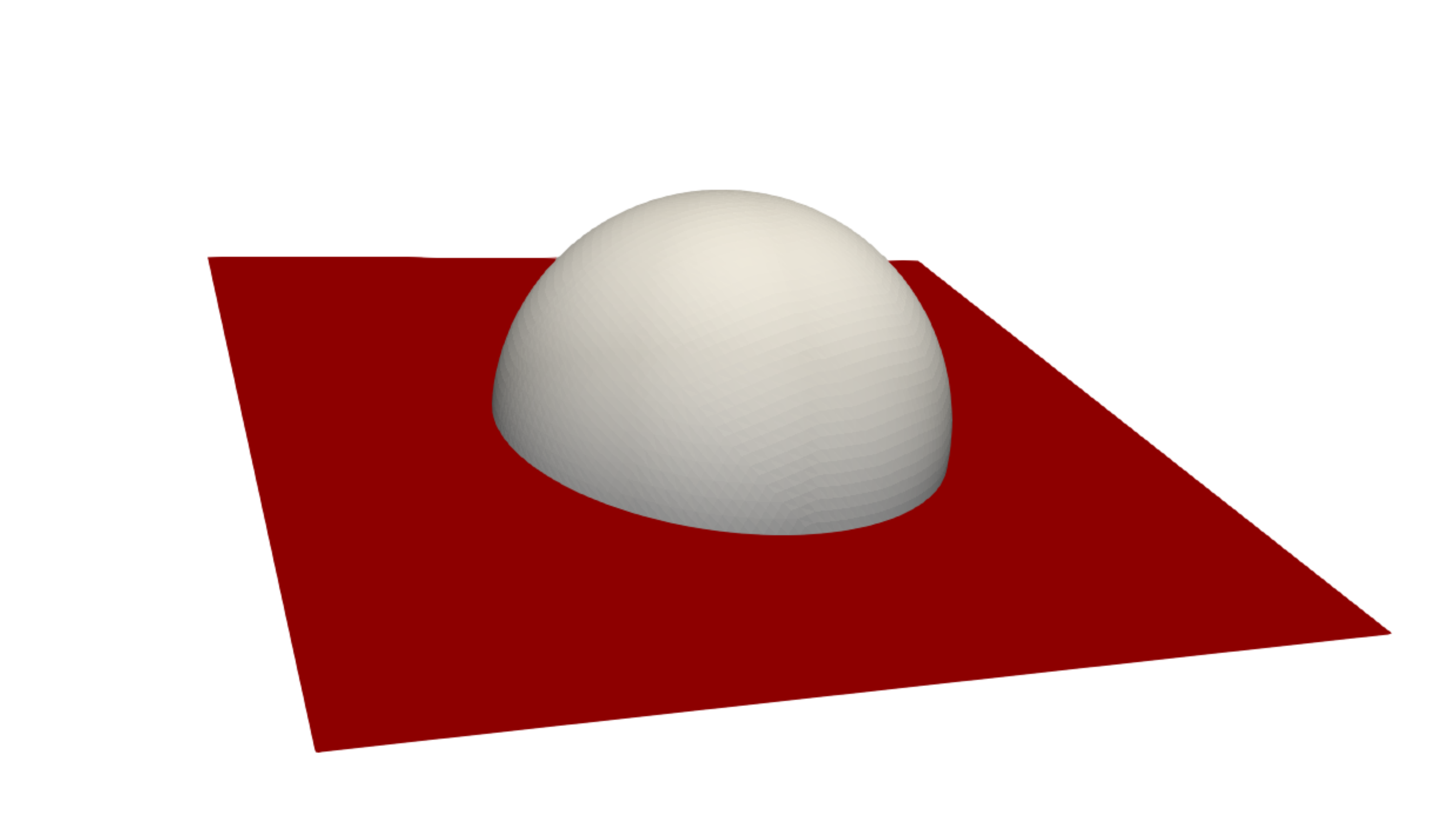}
\includegraphics[angle=-0,width=0.35\textwidth]{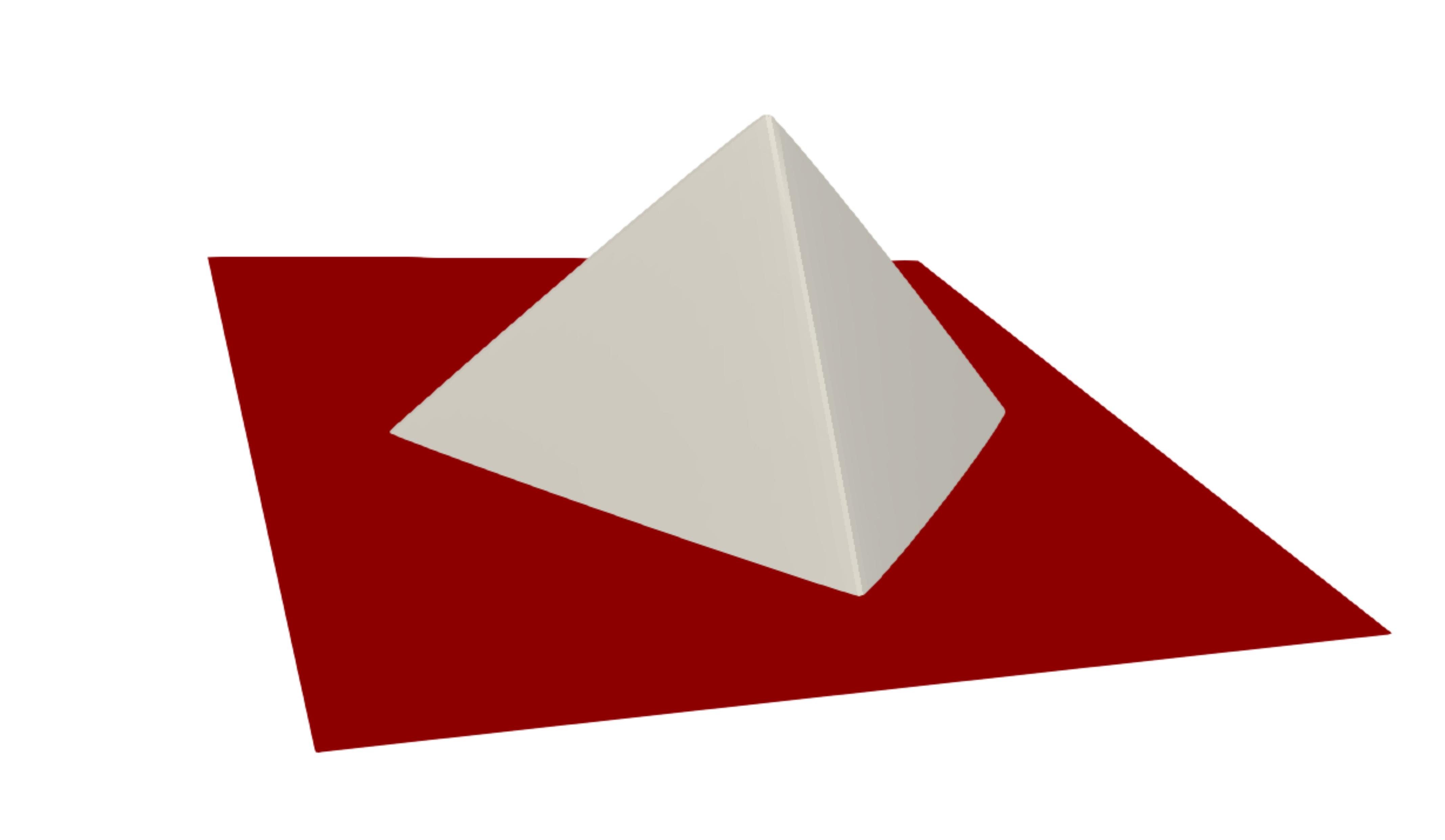} \\[0.5em]
\includegraphics[width=0.7\textwidth]{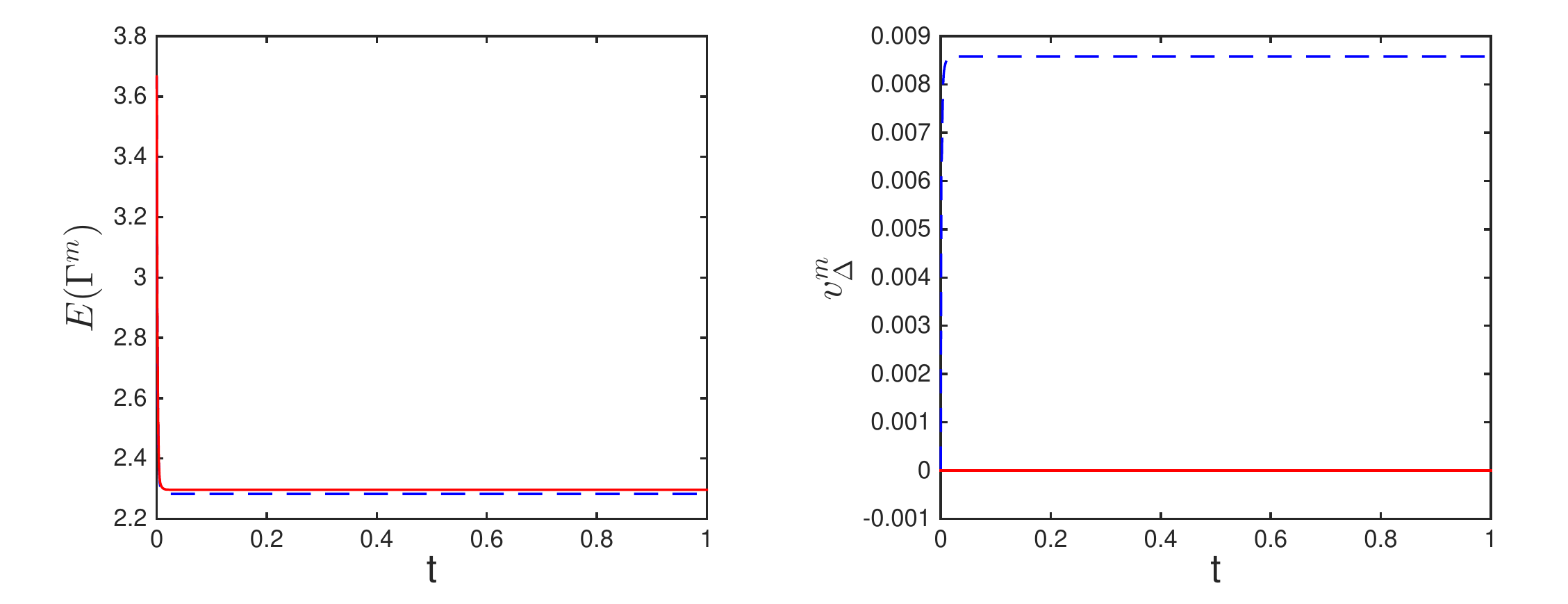}
\caption{
Evolution towards a drop on a substrate, with $\varrho=0.5$,
for the anisotropy \eqref{eq:cuspgamma} with $L=3$, $r=30$ and $\epsilon=0.1$.
Plots of $\Gamma^m$ at times $t=0, 1$.
We also show plots of the discrete energy $E(\Gamma^m)$ and 
the relative volume error $v_\Delta^m$ over time, where 
$K = 4225$ and $\Delta t = 10^{-3}$.
}
% ~/hpc_cluster/data/alberta/tjtrue/3d.anidrop_rho05
\label{fig:3danidrop05}
\end{figure}%

\begin{figure}[!htb]
\center
\includegraphics[angle=-0,width=0.35\textwidth]{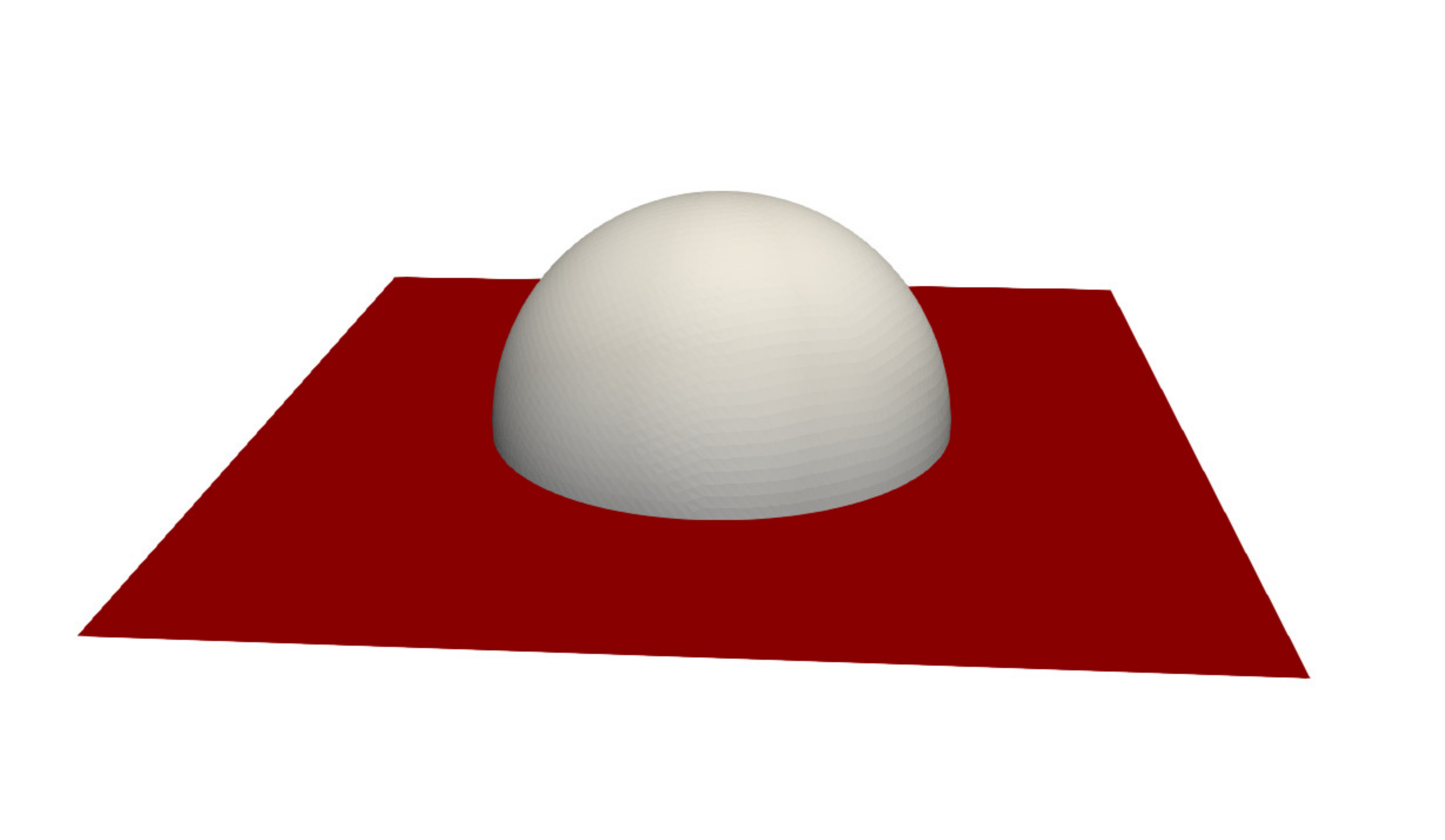}
\includegraphics[angle=-0,width=0.35\textwidth]{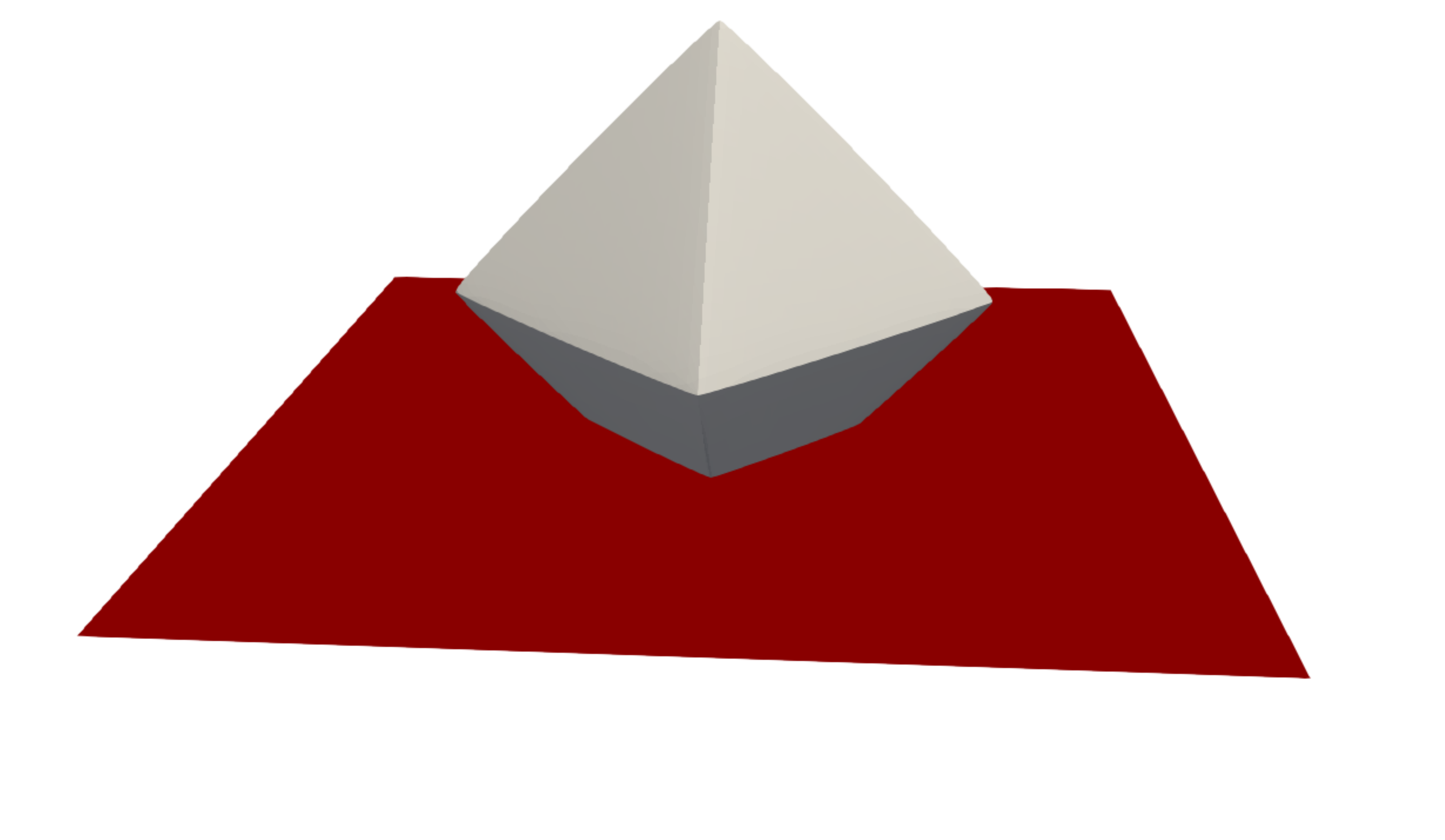} \\[0.5em]
\includegraphics[width=0.7\textwidth]{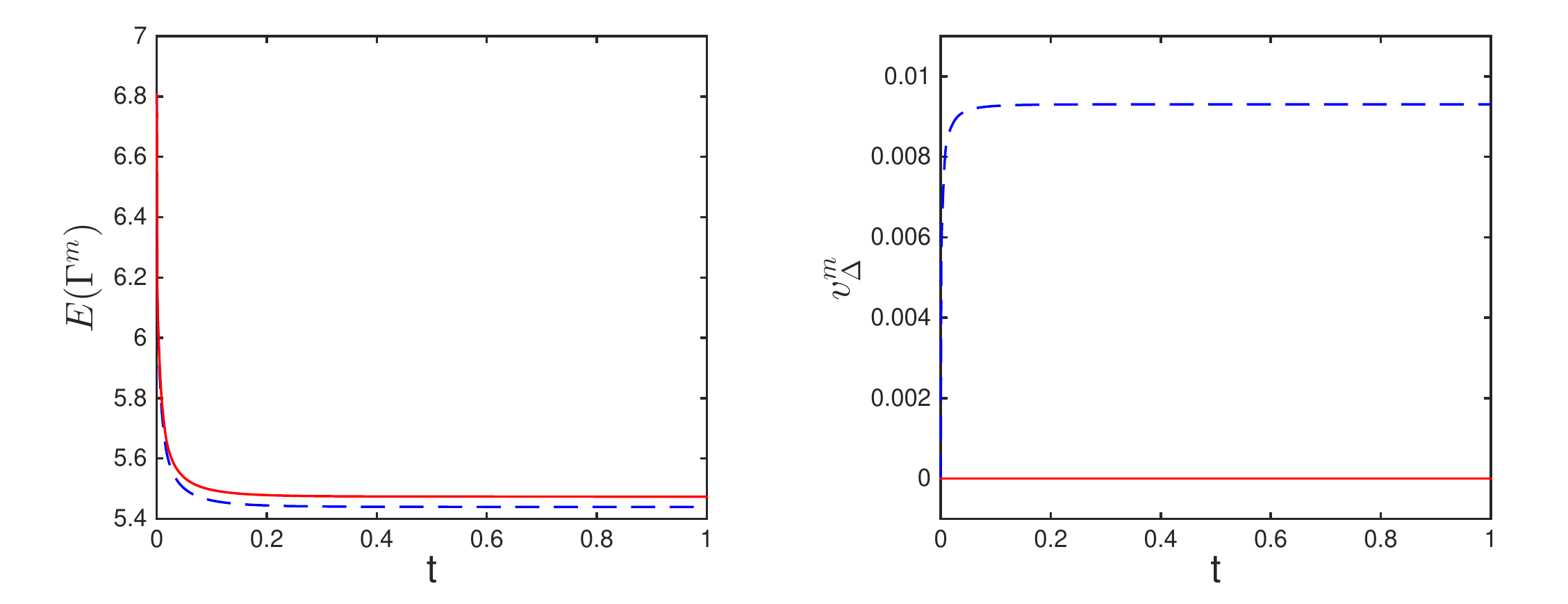}
\caption{
Evolution towards a drop on a substrate, with $\varrho=-0.5$,
for the anisotropy \eqref{eq:cuspgamma} with $L=3$, $r=30$ and $\epsilon=0.1$.
Plots of $\Gamma^m$ at times $t=0, 1$.
We also show plots of the discrete energy $E(\Gamma^m)$ and 
the relative volume error $v_\Delta^m$ over time, where 
$K = 4225$ and $\Delta t = 10^{-3}$.}
\label{fig:3danidrop-05}
\end{figure}%

Finally, we repeat the experiments in Figs.~\ref{fig:3disodrop05} and \ref{fig:3disodrop-05} but use the anisotropy in \eqref{eq:cuspgamma} with  $L=3$, $r=30$ and $\epsilon=0.1$. The simulation results are shown in Figs.~\ref{fig:3danidrop05} and \ref{fig:3danidrop-05}, where we observe the evolution of the drop is highly influenced by the chosen anisotropy $\gamma(\vec p)$ and the contact energy contribution parameter $\varrho$. 
We note that the numerical steady state for $\varrho=0.5$, which is 
visually nearly indistinguishable from the corresponding result for
$\varrho=0$, resembles the shapes of certain quantum dots, see e.g.\
\cite{AlshehriSAI19}.
%\footnote{Harald: Not sure if this is relevant. Add a sentence if so?}
Once again, we note that our numerical approximations exhibit the energy dissipation and volume conservation properties. 

\section{Conclusion} \label{sec:con}
In this work, we proposed a structure-preserving parametric finite element method for discretizing the surface diffusion of two-dimensional curve networks and three-dimensional surface clusters. The proposed method is based on an adaption of the BGN scheme from \cite{Barrett07,Barrett07Ani,Barrett10cluster,Barrett10} by using suitably time-weighted discrete normals,
and similarly appropriately weighted effective boundary velocity vectors,
instead of the conventional explicit treatment. As a consequence, the new method not only inherits the good mesh quality and the unconditional stability that the standard scheme enjoys, at least in the case of neutral external boundaries,
but also satisfies the exact volume conservation for each enclosed bubble in the system. In addition, the new scheme is also unconditionally stable in the case of non-neutral external boundaries.
These good properties were illustrated by numerical examples for the evolution of curve networks in 2d and surface clusters in 3d in the case of isotropic and anisotropic surface energies. Moreover, the reliability and applicability of the proposed scheme was demonstrated by comparing the numerical results with those of the standard BGN scheme. 

\section*{Acknowledgement}
The work of Bao was supported by the Ministry of Education of Singapore grant MOE2019-T2-1-063 (R-146-000-296-112). The work of Zhao was funded by the Alexander von Humboldt Foundation. 

\bibliographystyle{abbrv}
\bibliography{thebib}

\begin{thebibliography}{10}

\bibitem{AbelsAG15}
H.~Abels, N.~Arab, and H.~Garcke.
\newblock On convergence of solutions to equilibria for fully nonlinear
  parabolic systems with nonlinear boundary conditions.
\newblock {\em J. Evol. Equ.}, 15(4):913--959, 2015.

\bibitem{Abels2019}
H.~Abels, N.~Arab, and H.~Garcke.
\newblock Standard planar double bubbles are stable under surface diffusion
  flow.
\newblock {\em Commun. Anal. Geom.}, 29(5):1007--1060, 2021.

\bibitem{AlshehriSAI19}
K.~Alshehri, A.~Salhi, N.~Ahamad~Madhar, and B.~Ilahi.
\newblock Size and shape evolution of {GaAsSb-capped} {InAs/GaAs} quantum dots:
  Dependence on the {Sb} content.
\newblock {\em Crystals}, 9(10):530, 2019.

\bibitem{amilibia2001}
A.~M. Amilibia.
\newblock Existence and uniqueness of standard bubble clusters of given volumes
  in $\mathbb{R}^{N}$.
\newblock {\em Asian J. Math.}, 5(1):25--31, 2001.

\bibitem{AverbuchIR03}
A.~Averbuch, M.~Israeli, and I.~Ravve.
\newblock Electromigration of intergranular voids in metal films for
  microelectronic interconnects.
\newblock {\em J. Comput. Phys.}, 186:481--502, 2003.

\bibitem{Bansch05}
E.~B{\"a}nsch, P.~Morin, and R.~H. Nochetto.
\newblock A finite element method for surface diffusion: the parametric case.
\newblock {\em J. Comput. Phys.}, 203(1):321--343, 2005.

\bibitem{Bao17}
W.~Bao, W.~Jiang, Y.~Wang, and Q.~Zhao.
\newblock A parametric finite element method for solid-state dewetting problems
  with anisotropic surface energies.
\newblock {\em J. Comput. Phys.}, 330:380--400, 2017.

\bibitem{Zhao2021}
W.~Bao and Q.~Zhao.
\newblock A structure-preserving parametric finite element method for surface
  diffusion.
\newblock {\em SIAM J. Numer. Anal.}, 59(5):2775--2799, 2021.

\bibitem{BaoZ20preprint}
W.~Bao and Q.~Zhao.
\newblock An energy-stable parametric finite element method for simulating
  solid-state dewetting problems in three dimensions.
\newblock {\em J. Comput. Math.}, to appear, 2022.

\bibitem{Barrett07Ani}
J.~W. Barrett, H.~Garcke, and R.~N{\"u}rnberg.
\newblock Numerical approximation of anisotropic geometric evolution equations
  in the plane.
\newblock {\em IMA J. Numer. Anal.}, 28(2):292--330, 2007.

\bibitem{Barrett07b}
J.~W. Barrett, H.~Garcke, and R.~N{\"u}rnberg.
\newblock On the variational approximation of combined second and fourth order
  geometric evolution equations.
\newblock {\em SIAM J. Sci. Comput.}, 29(3):1006--1041, 2007.

\bibitem{Barrett07}
J.~W. Barrett, H.~Garcke, and R.~N{\"u}rnberg.
\newblock A parametric finite element method for fourth order geometric
  evolution equations.
\newblock {\em J. Comput. Phys.}, 222(1):441--467, 2007.

\bibitem{grain}
J.~W. Barrett, H.~Garcke, and R.~N\"urnberg.
\newblock A phase field model for the electromigration of intergranular voids.
\newblock {\em Interfaces Free Bound.}, 9(2):171--210, 2007.

\bibitem{Barrett08JCP}
J.~W. Barrett, H.~Garcke, and R.~N{\"u}rnberg.
\newblock On the parametric finite element approximation of evolving
  hypersurfaces in $\mathbb{R}^3$.
\newblock {\em J. Comput. Phys.}, 227(9):4281--4307, 2008.

\bibitem{Barrett08Ani}
J.~W. Barrett, H.~Garcke, and R.~N{\"u}rnberg.
\newblock A variational formulation of anisotropic geometric evolution
  equations in higher dimensions.
\newblock {\em Numer. Math.}, 109(1):1--44, 2008.

\bibitem{Barrett10}
J.~W. Barrett, H.~Garcke, and R.~N{\"u}rnberg.
\newblock Finite-element approximation of coupled surface and grain boundary
  motion with applications to thermal grooving and sintering.
\newblock {\em Eur. J. Appl. Math.}, 21(6):519--556, 2010.

\bibitem{Barrett10cluster}
J.~W. Barrett, H.~Garcke, and R.~N{\"u}rnberg.
\newblock Parametric approximation of surface clusters driven by isotropic and
  anisotropic surface energies.
\newblock {\em Interfaces Free Bound.}, 12(2):187--234, 2010.

\bibitem{Barrett2011}
J.~W. Barrett, H.~Garcke, and R.~N{\"u}rnberg.
\newblock The approximation of planar curve evolutions by stable fully implicit
  finite element schemes that equidistribute.
\newblock {\em Numer. Methods Partial Differ. Equ.}, 27(1):1--30, 2011.

\bibitem{Barrett20}
J.~W. Barrett, H.~Garcke, and R.~N{\"u}rnberg.
\newblock Parametric finite element approximations of curvature driven
  interface evolutions.
\newblock {\em Handb. Numer. Anal. (Andrea Bonito and Ricardo H. Nochetto,
  eds.)}, 21:275--423, 2020.

\bibitem{BowerC98}
A.~F. Bower and D.~Craft.
\newblock Analysis of failure mechanisms in the interconnect lines of
  microelectronic circuits.
\newblock {\em Fat. Frac. Eng. Mat. Struct.}, 21:611--630, 1998.

\bibitem{Brakke92}
K.~A. Brakke.
\newblock The surface evolver.
\newblock {\em Exp. Math.}, 1(2):141--165, 1992.

\bibitem{bronsard1995}
L.~Bronsard and B.~T. Wetton.
\newblock A numerical method for tracking curve networks moving with curvature
  motion.
\newblock {\em J. Comput. Phys.}, 120(1):66--87, 1995.

\bibitem{Cahn91}
J.~W. Cahn.
\newblock Stability, microstructural evolution, grain growth, and coarsening in
  a two-dimensional two-phase microstructure.
\newblock {\em Acta Metall.}, 39:2189--2199, 1991.

\bibitem{CahnH74}
J.~W. Cahn and D.~W. Hoffman.
\newblock A vector thermodynamics for anisotropic surfaces: {II.} {Curved} and
  faceted surfaces.
\newblock {\em Acta Metall.}, 22(10):1205--1214, 1974.

\bibitem{ConcusF74}
P.~Concus and R.~Finn.
\newblock On capillary free surfaces in the absence of gravity.
\newblock {\em Acta Math.}, 132(1):177--198, 1974.

\bibitem{CoxG04}
S.~J. Cox and F.~Graner.
\newblock Three-dimensional bubble clusters: {S}hape, packing, and growth rate.
\newblock {\em Phys. Rev. E}, 69(3):031409, 2004.

\bibitem{CoxGVM-PP03}
S.~J. Cox, F.~Graner, M.~F. Vaz, C.~Monnereau-Pittet, and N.~Pittet.
\newblock Minimal perimeter for {$N$} identical bubbles in two dimensions:
  calculations and simulations.
\newblock {\em Phil. Mag.}, 83(11):1393--1406, 2003.

\bibitem{CoxMG13}
S.~J. Cox, F.~Morgan, and F.~Graner.
\newblock Are large perimeter-minimizing two-dimensional clusters of equal-area
  bubbles hexagonal or circular?
\newblock {\em Proc. R. Soc. Lond. Ser. A Math. Phys. Eng. Sci.},
  469(2149):20120392, 10, 2013.

\bibitem{DaviG90}
F.~Davi and M.~E. Gurtin.
\newblock On the motion of a phase interface by surface diffusion.
\newblock {\em Z. Angew. Math. Phys.}, 41:782--811, 1990.

\bibitem{Davis04}
T.~A. Davis.
\newblock Algorithm 832: {UMFPACK} {V}4.3---an unsymmetric-pattern multifrontal
  method.
\newblock {\em ACM Trans. Math. Software}, 30(2):196--199, 2004.

\bibitem{DeckelnickDE05}
K.~Deckelnick, G.~Dziuk, and C.~M. Elliott.
\newblock Computation of geometric partial differential equations and mean
  curvature flow.
\newblock {\em Acta Numer.}, 14:139--232, 2005.

\bibitem{DepnerG13}
D.~Depner and H.~Garcke.
\newblock Linearized stability analysis of surface diffusion for hypersurfaces
  with triple lines.
\newblock {\em Hokkaido Math. J.}, 42(1):11--52, 2013.

\bibitem{ElliottG97a}
C.~M. Elliott and H.~Garcke.
\newblock Existence results for diffusive surface motion laws.
\newblock {\em Adv. Math. Sci. Appl.}, 7(1):465--488, 1997.

\bibitem{EscherMS98}
J.~Escher, U.~F. Mayer, and G.~Simonett.
\newblock The surface diffusion flow for immersed hypersurfaces.
\newblock {\em SIAM J. Math. Anal.}, 29(6):1419--1433, 1998.

\bibitem{Finn86}
R.~Finn.
\newblock {\em Equilibrium Capillary Surfaces}.
\newblock Grundlehren der Mathematischen Wissenschaften 284. Springer-Verlag,
  New York, 1986.

\bibitem{Foisy1993}
J.~Foisy, M.~Alfaro~Garcia, J.~Brock, N.~Hodges, and J.~Zimba.
\newblock The standard double soap bubble in {$R^2$} uniquely minimizes
  perimeter.
\newblock {\em Pac. J. Math.}, 159(1):47--59, 1993.

\bibitem{GarckeG20}
H.~Garcke and M.~G\"{o}{\ss}wein.
\newblock On the surface diffusion flow with triple junctions in higher space
  dimensions.
\newblock {\em Geom. Flows}, 5(1):1--39, 2020.

\bibitem{Garcke2021}
H.~Garcke and M.~G{\"o}{\ss}wein.
\newblock Non-linear stability of double bubbles under surface diffusion.
\newblock {\em J. Differ. Equ.}, 302:617--661, 2021.

\bibitem{GNSW}
H.~Garcke, B.~Nestler, B.~Stinner, and F.~Wendler.
\newblock {Allen-Cahn} systems with volume constraints.
\newblock {\em Math. Models Methods Appl. Sci.}, 18(08):1347--1381, 2008.

\bibitem{GarckeNS98}
H.~Garcke, B.~Nestler, and B.~Stoth.
\newblock On anisotropic order parameter models for multi-phase systems and
  their sharp interface limits.
\newblock {\em Physica D}, 115:87--108, 1998.

\bibitem{GarckeNC00}
H.~Garcke and A.~Novick-Cohen.
\newblock A singular limit for a system of degenerate {C}ahn--{H}illiard
  equations.
\newblock {\em Adv. Differential Equations}, 5(4-6):401--434, 2000.

\bibitem{Giga06}
Y.~Giga.
\newblock {\em Surface evolution equations}, volume~99 of {\em Monographs in
  Mathematics}.
\newblock Birkh\"{a}user, Basel, 2006.

\bibitem{Giga98}
Y.~Giga and K.~Ito.
\newblock On pinching of curves moved by surface diffusion.
\newblock {\em Commun. Appl. Anal.}, 2:393--405, 1998.

\bibitem{HausserV06a}
F.~Hau{\ss}er and A.~Voigt.
\newblock A discrete scheme for parametric anisotropic surface diffusion.
\newblock {\em J. Sci. Comput.}, 30(2):223--235, 2007.

\bibitem{Hoffman72}
D.~W. Hoffman and J.~W. Cahn.
\newblock A vector thermodynamics for anisotropic surfaces: {I. Fundamentals}
  and application to plane surface junctions.
\newblock {\em Surf. Sci.}, 31:368--388, 1972.

\bibitem{HMRR}
M.~Hutchings, F.~Morgan, M.~Ritor{\'e}, and A.~Ros.
\newblock Proof of the double bubble conjecture.
\newblock {\em Ann. of Math. (2)}, 155(2):459--489, 2002.

\bibitem{Jiang2021}
W.~Jiang and B.~Li.
\newblock A perimeter-decreasing and area-conserving algorithm for surface
  diffusion flow of curves.
\newblock {\em J. Comput. Phys.}, 443:110531, 2021.

\bibitem{JiangZB20}
W.~Jiang, Q.~Zhao, and W.~Bao.
\newblock Sharp-interface model for simulating solid-state dewetting in three
  dimensions.
\newblock {\em SIAM J. Appl. Math.}, 80(4):1654--1677, 2020.

\bibitem{Kovacs2020}
B.~Kov{\'a}cs, B.~Li, and C.~Lubich.
\newblock A convergent evolving finite element algorithm for {Willmore} flow of
  closed surfaces.
\newblock {\em Numer. Math.}, 149(3):595--643, 2021.

\bibitem{KraynikRS04}
A.~M. Kraynik, D.~A. Reinelt, and F.~van Swol.
\newblock Structure of random foam.
\newblock {\em Phys. Rev. Lett.}, 93(20):208301, 2004.

\bibitem{Li2020energy}
Y.~Li and W.~Bao.
\newblock An energy-stable parametric finite element method for anisotropic
  surface diffusion.
\newblock {\em J. Comput. Phys.}, 446:110658, 2021.

\bibitem{LiZG99}
Z.~Li, H.~Zhao, and H.~Gao.
\newblock A numerical study of electro-migration voiding by evolving level set
  functions on a fixed cartesian grid.
\newblock {\em J. Comput. Phys.}, 152:281--304, 1999.

\bibitem{MerrimanBO94}
B.~Merriman, J.~K. Bence, and S.~J. Osher.
\newblock Motion of multiple functions: a level set approach.
\newblock {\em J. Comput. Phys.}, 112(2):334--363, 1994.

\bibitem{Morgan2007}
F.~Morgan.
\newblock Colloquium: Soap bubble clusters.
\newblock {\em Rev. Mod. Phys.}, 79(3):821, 2007.

\bibitem{Morgan1998wulff}
F.~Morgan, C.~French, and S.~Greenleaf.
\newblock Wulff clusters in {$R^2$}.
\newblock {\em J. Geom. Anal.}, 8(1):97--115, 1998.

\bibitem{Mullins57}
W.~W. Mullins.
\newblock Theory of thermal grooving.
\newblock {\em J. Appl. Phys.}, 28(3):333--339, 1957.

\bibitem{Mullins58}
W.~W. Mullins.
\newblock The effect of thermal grooving on grain boundary motion.
\newblock {\em Acta Metall.}, 6(6):414--427, 1958.

\bibitem{GNSSW}
B.~Nestler, F.~Wendler, M.~Selzer, B.~Stinner, and H.~Garcke.
\newblock Phase-field model for multiphase systems with preserved volume
  fractions.
\newblock {\em Phys. Rev. E}, 78(1):011604, 2008.

\bibitem{Neubauer2002}
R.~Neubauer.
\newblock Ein {F}initeelementeansatz f{\"u}r {K}r{\"u}mmungsflu{\ss} von unter
  {T}ripelpunktbedingungen verbundenen {K}urven.
\newblock Master's thesis, University Bonn, Bonn, 2002.

\bibitem{Robert09}
R.~N\"urnberg.
\newblock Numerical simulations of immiscible fluid clusters.
\newblock {\em Appl. Numer. Math.}, 59:1612--1628, 2009.

\bibitem{mullins}
R.~N\"urnberg.
\newblock A structure preserving front tracking finite element method for the
  {M}ullins--{S}ekerka problem.
\newblock arXiv: 2111.15418, 2021.

\bibitem{Pan03}
J.~Pan.
\newblock Modelling sintering at different length scales.
\newblock {\em Int. Mater. Rev.}, 48(2):69--85, 2003.

\bibitem{PanW08}
Z.~Pan and B.~Wetton.
\newblock A numerical method for coupled surface and grain boundary motion.
\newblock {\em European J. Appl. Math.}, 19(3):311--327, 2008.

\bibitem{PaoliniT20}
E.~Paolini and V.~M. Tortorelli.
\newblock The quadruple planar bubble enclosing equal areas is symmetric.
\newblock {\em Calc. Var. Partial Differential Equations}, 59(1):20, 2020.

\bibitem{Ruuth1998}
S.~J. Ruuth.
\newblock Efficient algorithms for diffusion-generated motion by mean
  curvature.
\newblock {\em J. Comput. Phys.}, 144(2):603--625, 1998.

\bibitem{Alberta}
A.~Schmidt and K.~G. Siebert.
\newblock {\em Design of Adaptive Finite Element Software: The Finite Element
  Toolbox {ALBERTA}}, volume~42 of {\em Lecture Notes in Computational Science
  and Engineering}.
\newblock Springer-Verlag, Berlin, 2005.

\bibitem{SmithSC02}
K.~A. Smith, F.~J. Solis, and D.~L. Chopp.
\newblock A projection method for motion of triple junctions by levels sets.
\newblock {\em Interfaces Free Bound.}, 4(3):263--276, 2002.

\bibitem{sullivan1996}
J.~M. Sullivan and F.~Morgan.
\newblock Open problems in soap bubble geometry.
\newblock {\em Int. J. Math.}, 7(06):833--842, 1996.

\bibitem{Taylor76}
J.~E. Taylor.
\newblock The structure of singularities in soap-bubble-like and soap-film-like
  minimal surfaces.
\newblock {\em Ann. of Math. (2)}, 103(3):489--539, 1976.

\bibitem{Taylor99}
J.~E. Taylor.
\newblock A variational approach to crystalline triple-junction motion.
\newblock {\em J. Stat. Phys.}, 95(5):1221--1244, 1999.

\bibitem{Thaddey99}
B.~Thaddey.
\newblock {N}umerik f\"ur die {E}volution von {K}urven mit {T}ripelpunkt.
\newblock Master's thesis, University Freiburg, Freiburg, 1999.

\bibitem{Wecht2000double}
B.~Wecht, M.~Barber, and J.~Tice.
\newblock Double crystals.
\newblock {\em Acta Crystallographica, Sect. A}, 56(1):92--95, 2000.

\bibitem{Wichiramala04}
W.~Wichiramala.
\newblock Proof of the planar triple bubble conjecture.
\newblock {\em J. Reine Angew. Math.}, 567:1--49, 2004.

\bibitem{ZhaoMOW98}
H.-K. Zhao, B.~Merriman, S.~Osher, and L.~Wang.
\newblock Capturing the behavior of bubbles and drops using the variational
  level set approach.
\newblock {\em J. Comput. Phys.}, 143(2):495--518, 1998.

\bibitem{Zhao2020parametric}
Q.~Zhao, W.~Jiang, and W.~Bao.
\newblock A parametric finite element method for solid-state dewetting problems
  in three dimensions.
\newblock {\em SIAM J. Sci. Comput.}, 42(1):B327--B352, 2020.

\bibitem{Zhao20}
Q.~Zhao, W.~Jiang, and W.~Bao.
\newblock An energy-stable parametric finite element method for simulating
  solid-state dewetting.
\newblock {\em IMA J. Numer. Anal.}, 41(3):2026--2055, 2021.

\end{thebibliography}
\end{document}